\documentclass[3p,times]{elsarticle}

\usepackage{latexsym}

\usepackage{indentfirst}
\usepackage[english]{babel}
\usepackage{cancel}

\usepackage{lipsum}
\usepackage{array}
\usepackage{color} 
\usepackage{tabularx}
\usepackage{graphicx} 
\usepackage{amsmath}
\usepackage{amssymb}
\usepackage{amsthm}
\newtheorem{proposition}{Proposition}
\usepackage{amsfonts}	
\usepackage{moreverb}
\usepackage{dsfont}
\usepackage{tipa}
\usepackage{upgreek}
\usepackage{bm}
\usepackage{multirow}
\usepackage{soul}
\usepackage{ textcomp }
\usepackage[dvipsnames]{xcolor}
\usepackage{mathtools}
\usepackage{setspace}
\usepackage{algorithm}
\usepackage{algorithmic}
\allowdisplaybreaks

\newcommand\BibTeX{{\rmfamily B\kern-.05em \textsc{i\kern-.025em b}\kern-.08em
		T\kern-.1667em\lower.7ex\hbox{E}\kern-.125emX}}

\usepackage{hyperref} 
\hypersetup{
	colorlinks=true,                          
	linkcolor=blue, 
	citecolor=red, 
	urlcolor=blue  } 

\newcommand{\x}{\mathbf{x}}
\newcommand{\X}{\mathbf{X}}

\newcommand{\mbf}[1]{\mathbf{#1}}			%

\newcommand{\Q}{\mathbf{Q}}

\newcommand{\w}{\mathbf{w}}

\newcommand{\q}{\mathbf{q}}
\newcommand{\F}{\mathbf{F}}
\newcommand{\f}{\mathbf{f}}

\newcommand{\n}{\mathbf{n}}
\newcommand{\de}[2]{\frac {\partial #1}{\partial#2}}

\renewcommand{\v}{\mathbf{v}}

\newcommand{\RIIcolor}[1]{{\leavevmode\color{black} #1}}
\newcommand{\RIcolor}[1]{{\leavevmode\color{black} #1}}

\newcommand{\K}{\mathbf{K}}

\newcommand{\G}{\mathbf{G}}
\newcommand{\A}{\AAA}

\newcommand{\halb}{\frac{1}{2}}

\newcommand{\bpsi}{\boldsymbol{\psi}}
\newcommand{\bxi}{\boldsymbol{\xi}}
\newcommand{\diag}{\textnormal{diag}}

\newcommand{\bdm}{\begin{displaymath}}
	\newcommand{\edm}{\end{displaymath}}

\newcommand{\bea}{\begin{eqnarray} }
	\newcommand{\eea}{\end{eqnarray} }

\newcommand{\AAA}{{\boldsymbol{A}}}

\newcommand{\JJ}{{\mathbf{J}}}

\newcommand{\Popt}{\mathbf{P}_{\mathrm{opt}}}
\newcommand{\xbin}{\mbf{x}_c}

\newfont{\numerikEleven}{ecrm1000}
\newfont{\numerikTen}{cmss10}
\newfont{\numerikNine}{cmss9}
\newfont{\numerikEight}{cmss8}

\journal{Journal of Computational Physics}

\begin{document} 
	\begin{frontmatter}
		\title{On general and complete multidimensional Riemann solvers for nonlinear systems of hyperbolic conservation laws}   
		\author[UniVR]{Elena Gaburro\corref{cor1}}
		\ead{elena.gaburro@univr.it}
		\cortext[cor1]{Corresponding author}
		
		\author[INRIA]{Mario Ricchiuto}
		\ead{mario.ricchiuto@inria.fr}
		
		\author[UniTN]{Michael Dumbser}
		\ead{michael.dumbser@unitn.it}
		
		\address[UniVR]{Department of Computer Science, University of Verona, Strada le Grazie 15, Verona, 37134, Italy}		
		
		\address[INRIA]{Inria, Univ. Bordeaux, CNRS, Bordeaux INP, IMB, UMR 5251, 200 Avenue de la Vieille Tour, 33405 Talence, France}		
		
		\address[UniTN]{Laboratory of Applied Mathematics, DICAM, University of Trento, via Mesiano 77, 38123 Trento, Italy} 
		
		
		\begin{abstract} 
			In this work, we introduce  a framework to design  multidimensional Riemann solvers for nonlinear systems of hyperbolic conservation laws on general unstructured polygonal Voronoi-like tessellations. 
			In this framework
			we propose two simple but complete solvers.
			The first method is a direct extension of the Osher-Solomon Riemann solver to multiple space dimensions. Here, the multidimensional numerical dissipation is obtained by integrating the absolute value of the flux Jacobians over a dual triangular mesh around each node of the primal polygonal grid. The required nodal gradient is then evaluated on a local computational simplex involving the  $d+1$
			neighbors meeting at each corner.
			The second method is a genuinely multidimensional upwind flux. 
			By   introducing a fluctuation form
			of finite volume methods with corner fluxes, 
			we  show an equivalence  with  residual distribution schemes (RD).
			This   naturally allows to construct genuinely multidimensional upwind corner fluxes starting from existing and well-known RD fluctuations.  
			In particular, we explore the use of corner fluxes obtained from 
			the so-called N scheme, i.e. the Roe's original optimal multidimensional upwind advection scheme. 
			
			Both methods use the full eigenstructure of the underlying hyperbolic system and are therefore complete by construction.  
			A simple higher order extension up to fourth order in space and time is then introduced at the aid of a CWENO reconstruction in space and an ADER approach in time, leading to a family of high order accurate fully-discrete one-step schemes based on genuinely multidimensional Riemann solvers.  
			
			We present applications of our new numerical schemes to several test problems for the compressible Euler equations of gas-dyanamics. In particular, we show that the proposed schemes are at the same time carbuncle-free and preserve certain stationary shear waves exactly.

		\end{abstract}
		
		\begin{keyword}
			finite volume (FV) and residual distribution (RD) schemes for hyperbolic conservation laws \sep
			genuinely multidimensional Riemann solvers \sep 
			multidimensional Osher scheme \sep
			N scheme applied in the context of FV schemes \sep 
			unstructured polygonal meshes \sep
			compressible Euler equations \sep
		\end{keyword}
	\end{frontmatter}
	


	\section{Introduction} 
	\label{sec.introduction}
	
	The groundbreaking work of Godunov~\cite{godunov} on finite volume (FV) schemes in 1959 \RIIcolor{made it possible} to solve nonlinear hyperbolic systems of \RIIcolor{partial differential equations} (PDEs) with discontinuous solutions numerically\RIIcolor{.
		This was achieved by making use,} for the first time, of the exact solution of the Riemann problem (RP) between two adjacent states at a cell interface, in order to obtain a stable numerical flux for systems of nonlinear hyperbolic conservation laws via nonlinear upwinding.
	Later, Godunov-type finite volume schemes were developed that allowed also the use of only approximate solutions of the RP to get the numerical fluxes, so-called approximate Riemann solvers, see e.g.~\cite{roe1d,hll,hllem,osherandsolomon,Torohllc,torobook}. Since then, several attempts have been made to extend the idea of nonlinear upwinding via the exact or approximate solution of Riemann problems also to the multidimensional case. First, developments based on genuinely multidimensional Riemann solvers have been made by Roe \textit{et al.} in~\cite{RoeMultiD,RumseyLeerRoe}, Colella~\cite{ColellaMultiD} and Saltzman~\cite{SaltzmanMultiD} and later also in the context of finite volume evolution of Galerkin methods~\cite{FVEG}. Some thorough studies of the solution of two-dimensional Riemann problems have been provided in the work of Schulz-Rinne \textit{et al.} in~\cite{schulzrinne,RinneCollinsGlaz}. 
	
	A significant step forward in the development of genuinely multidimensional Riemann solvers for high order Godunov-type  FV schemes was made by Balsara \textit{et al.} in a series of fundamental papers~\cite{balsarahlle2d,balsarahllc2d,BalsaraMultiDRS,MUSIC1,MUSIC2,ADERdivB}, where a new family of genuinely multidimensional Riemann solvers of the HLL-type has been developed, including also multidimensional HLLC-type solvers. In all the aforementioned \RIIcolor{works}, the multidimensional Riemann solvers were not only interesting from a mathematical and \RIIcolor{conceptual} point of view, but also a necessary key building block for the construction of exactly divergence-free numerical schemes for magneto-hydrodynamics (MHD). 
	Along these lines, later, multidimensional Riemann solvers have also been successfully employed in the context of semi-implicit schemes for divergence-free MHD, see e.g.~\cite{SIMHD,HybridHexa1,CompatibleMHD}. We also would like to mention a very recent contribution made in~\cite{PHRaph1,PHRaph2}, \RIIcolor{extending}, for the first time ever, commonly used nodal solvers from the context of Lagrangian hydrodynamics~\cite{Despres2005,Despres2009,Maire2007,Maire2008,Maire2009,Maire2010,Maire2011} (where the use of multidimensional node solvers is mandatory to guarantee geometric compatibility on moving unstructured meshes), to the Eulerian description of compressible gas dynamics on fixed grids.  
	
	\bigskip 
	\RIIcolor{In paralle} to the efforts  to formulate  genuinely multidimensional Riemann solvers, \RIIcolor{there have been many works over the years dedicated} to 
	devise formulations which dot not directly rely on the Riemann Problem, and  may allow more naturally to go beyond the use of one dimensional wave propagation ideas. 
	These works have led first to the design of the schemes now  known as  residual distribution (RD)~\cite{abgrall2022hyperbolic}.
	These methods started as generalizations  to multiple dimensions P.L. Roe's upwind splitting~\cite{roe1d,Roe82}. 
	Over a decade, joint efforts of groups revolving around   P.L. Roe and H. Deconinck have managed to propose 
	numerical methods for the Euler equations  combining a multidimensional generalization
	of Roe's linearization,  a   decomposition of the steady state equations in scalar components, which allowed
	a unique and  efficient  approximation of steady   two-dimensional  flows~\cite{Roe:87,Roe:90,Roe92,Roe93,Dec93,Deconinck1993,hiro3}.
	The methods proposed relied on splitting  flux integrals to cell nodes.   In this new setting upwinding is not anymore  related to choosing information  between  two  states, 
	but allows to select   as target for the flux  update among  three states  at once on a triangle, four on a tetrahedron or quadrangle, etc.
	Among these methods, the so called N scheme emerged as the optimal  positive coefficient first order scheme:  featuring the minimal cross-wind numerical dissipation  
	among positivity preserving first order upwind schemes, and providing exact preservation of steady data  for mesh aligned advection~\cite{Roe:87,Roe:90,Roe92}.
	
	The methods have seen since then several evolutions. Recasting system  equations as coupled scalar  wave equations has been often replaced
	with a full matrix   generalization only requiring the knowledge of the eigenstructure of the  Jacobian matrices of the underlying hyperbolic equations~\cite{vanderweide}.
	The requirement of an exact linearization  of the    latter matrices  to evaluate multidimensional upwind distribution, required e.g. for the N scheme,
	was removed in~\cite{crd02,rcd05} allowing the use of a more practical  direct quadrature of the fluxes,  instead of the quadrature of the  linearized quasi-linear form.
	Ad-hoc fixes for unwanted features as Carbuncle's and expansion shocks arising in high speed applications have been
	also proposed~\cite{GARICANOMENA201643,rd-efix}. 
	\RIIcolor{The appearance of such anomalies is partly a consequence of the low cross-stream dissipation of the N scheme, which, on the other hand, is appealing for the approximation of complex flows with contact waves and boundary layers on coarse meshes~\cite{CiCP-13-479}. 
		Higher-order extensions for steady-state hyperbolic equations have been proposed in~\cite{abg2003,alr11}, and the high-order extension to time-dependent problems has been clarified in~\cite{abg2001c,rcd05,rad07,ra10,abg2017}.}
	Recently,  the RD setting has been shown to allow naturally the satisfaction of several additional conservation constraints~\cite{abgrall2018general,abgrall2022reinterpretation,Abgrall23}.
	Many of these contributions are   summarized and compared in the reviews~\cite{RD-ency,mr11,RD-ency2,abgrall2022hyperbolic}.  These developments have required  somewhat steering away   from some of the initial ideas which motivated the design of the methods.  This has
	led   to more recent work on a different setting called Active Flux methods (see e.g. the recent lecture~\cite{Roe24} and~\cite{barsukow2019active,abgrall2023extensions,abg_bar_klin_2025semi}).\\

	In this paper we proceed differently and  combine    ideas  from  the multidimensional upwind RD  and multidimensional Riemann solver  frameworks. This provides a new setting to design 
	finite volume methods  
	with   multiple state corner numerical fluxes, and  genuinely multidimensional dissipation. 
	Within this setting we propose and test two specific definitions of the corner fluxes. The first  is a direct extension of the Osher-Solomon Riemann solver to multiple space dimensions. Here, the multidimensional numerical dissipation is obtained by integrating the absolute value of the flux Jacobians over a dual triangular mesh around each node of the primal polygonal grid. The required nodal gradient is then evaluated on a local computational simplex involving the  $d+1$
	neighbors meeting at each corner. The second method is a genuinely multidimensional upwind flux,
	obtained by an equivalence with RD schemes on the mentioned  local computational simplex. 
	In particular, the corner flux is defined as the internal normal flux plus a fluctuation  obtained from  the   N scheme.\\

	The rest of the paper is organized as follows. In Section~\ref{sec.setting} we briefly introduce the overall setting of the paper, including an introduction to the notation and the class of governing PDEs under consideration. Subsequently, in Section~\ref{sec.osher.o1} we present a new genuinely multidimensional Osher-type Riemann solver, which is a straightforward extension of~\cite{OsherNC,OsherUniversal} to the multidimensional case. A peculiar feature of the new multidimensional Osher-type flux proposed in this paper is the fact that it provides an expression for the entire multidimensional numerical flux tensor, and not just a numerical flux projected into a particular normal direction. 
	In Section~\ref{sec.rd.o1} we present a fluctuation formulation of
	multidimensional Riemann solvers with corner fluxes, allowing to establish
	an analogy with the so-called  residual distribution (RD) framework
	This provides a whole new setting to design  corner fluxes  for finite volume (FV) schemes, embedding   a genuinely multidimensional upwind flavor.  
	As an example we will use in the tests numerical fluxes based on a generalization of the  N scheme.
	All multidimensional fluxes presented before can then be naturally extended to higher order in space and time via a CWENO reconstruction in space~\cite{CWENO1,CWENO2,CWENO3,cravero2018cweno,ADERCWENO,gaburro2020high} and a fully-discrete one-step ADER discretization in time~\cite{toro2005ader,toro3,dumbser2008unified,frontiers2020}, see Section~\ref{sec.highorder}. Comparative computational results are presented for the different schemes in Section~\ref{sec.tests}. The paper closes with some concluding remarks and an outlook to future research given in Section~\ref{sec.conclusions}.  
	
	\section{Governing PDEs and domain discretization}  
	\label{sec.setting}
	
	We aim at solving general nonlinear hyperbolic systems of conservation laws in multiple space dimensions that read 
	\begin{equation}
		\label{eq.pde1}
		\partial_t \Q + \nabla \cdot \F(\Q) = 0 \quad \textnormal{on} \quad 
		\Omega \times [0,t_f] \subset \mathbb{R}^d \times\mathbb{R}^+
	\end{equation}
	with $\Q \in \Omega_\Q \subset \mathbb{R}^m$ the vector of $m$ conservative variables, $\F=\F(\Q)=\left( \f_1(\Q),\f_2(\Q),\dots,\f_d(\Q) \right)$ the flux tensor
	and $d$ the number of space dimensions. 
	The system can be also written in the following quasi-linear form as 
	\RIIcolor{
		\begin{equation}
			\label{eq.pde2}
			\partial_t \Q +   \A(\Q) \cdot \nabla \Q = 0, \quad \text{or equivalently}  \quad 
			\partial_t \Q +  \sum_{j=1}^d \A_j(\Q) \, \frac{\partial \Q}{\partial x_j}  = 0,
	\end{equation}}
	with $\A(\Q) = \left( \A_1(\Q),\A_2(\Q),\dots,\A_d(\Q) \right)$
	\RIIcolor{the rank 3 tensor gathering} the Jacobian matrices of the fluxes $\A_i(\Q) = \partial \f_i(\Q)/\partial \Q$ \RIIcolor{and 
		$\nabla \Q$ made by the $d$ columns of the Jacobian of $Q$}.
	
	\RIIcolor{We assume the above system to verify all} the classical properties of a hyperbolic system, i.e. it is endowed with a (convex) entropy pair, the components  of \RIIcolor{$\K_{n} := \A(\Q) \cdot \n $  (defined here below in details)} have a full set of real eigenvalues and linearly independent eigenvectors for all normal vectors $\n \neq 0$. 
	For later use we define, given a vector $ \n = (n_1, n_2, \dots, n_d)  \in\mathbb{R}^d$ the projected  matrix  
	\begin{equation}\label{eq.pde3}
		\K_{n}:= \A(\Q) \cdot  \n = \sum_{i=1}^d \A_i(\Q) \, n_i.
	\end{equation}
	Since the system is hyperbolic $\K_{n}$ also has a full set of real eigenvalues and linearly independent eigenvectors. 
	
	The computational domain $\Omega$ is discretized via a set of Voronoi-like polygons/polyhedra $\Omega_c$. 
	To construct the polygonal/polyhedral tessellation we first cover our domain with a set of generator points which are connected via a Delaunay triangulation in 2D or tetrahedralization in 3D; 
	then, each polygon/polyhedron is constructed around a generator point by connecting the centroids (i.e. barycenters) of the dual triangles/tetrahedra sharing that generator. 
	Since we employ the centroids of the triangles/tetrahedra, and not the circumcenters as in the classical definition of Voronoi elements, we should always refer to our elements as Voronoi-like polygons/polyhedra, but sometimes, to lighten the notation, we will use just the word Voronoi. 
	The use of centroid-based Voronoi-like elements avoids small edges/faces and in particular zero length/area ones, hence it increases the homogeneity and the quality of the obtained tessellation. 
	Moreover, it guarantees that each vertex of the polygons (respectively polyhedra), a part for those on a boundary, is attached to exactly 3 Delaunay triangles (respectively 4 Delaunay tetrahedra). 
	
	In what follows, we restrict our notation for the geometry to the two dimensional case $d=2$ and we will make explicit the three-dimensional $d=3$ extension only in the multidimensional flux definitions.

	A considered generator, which is also a vertex of a Delaunay triangle, is indicated with the letter ${c}$, and eventually, if needed, with the letters ${a,b}$. It is depicted with a point symbol in the Figures and its coordinates are indicated with $\x_c$ (respectively $\x_a, \x_b$). 
	Correspondingly, the Voronoi element around this point is indicated with $\Omega_c$ (respectively $\Omega_a, \Omega_b$). 
	
	Two elements $\Omega_a$ and $\Omega_c$ which share the same edge/face $\partial \Omega_{ac} = \Omega_c \cap \Omega_a$ are called neighbors. 
	The set of neighbors of the polygonal element $\Omega_c$ is denoted by $\mathcal{N}_c$ and contains all polygons that share an edge with $\Omega_c$. The \textit{unit} normal vector pointing from $\Omega_c$ to $\Omega_a$ is denoted by $\hat{\n}_{ac}$, while the vector including the length 
	$$|\partial \Omega_{ac}| = \int \limits_{\partial \Omega_{ac}} dS$$ 
	is denoted by $\n_{ac} = \hat{\n}_{ac} |\partial \Omega_{ac}|$. 
	
	A generic vertex of the Voronoi tessellation is denoted by $p$ (and eventually, if needed, with the letters ${r,q}$). It is depicted with a star symbol in the Figures and its coordinates are denoted with $\x_p$.
	These points are indeed the barycenters of the Delaunay triangles: we denote with $T_p$ the Delaunay triangle around each~$p$ whose vertexes consequently are the generators of the $3$ polygons which are attached to the point~$p$. 
	Furthermore, the set of polygonal cells around a point $p$ is denoted by $\mathcal{C}_p$ and the set of vertexes $p$ that compose a polygon $\Omega_c$ is denoted by $\mathcal{P}_c$.
	
\begin{figure}[h]
		\begin{minipage}{0.5\textwidth}
			\centering\includegraphics[width=0.8\textwidth]{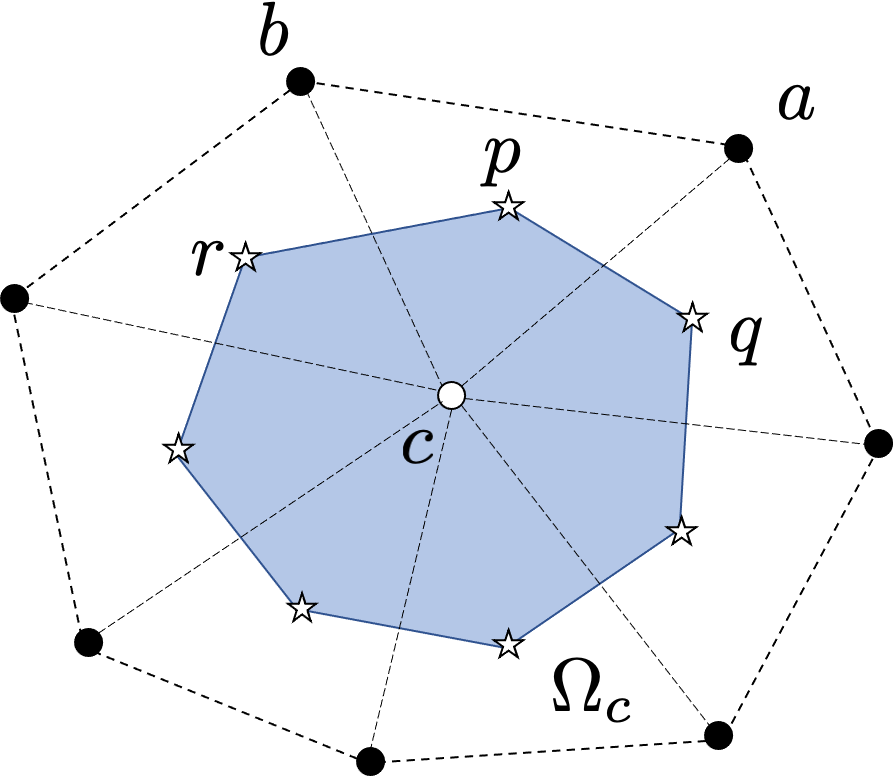}
		\end{minipage}\hfill
		\begin{minipage}{0.5\textwidth}
			\centering\includegraphics[width=0.7\textwidth]{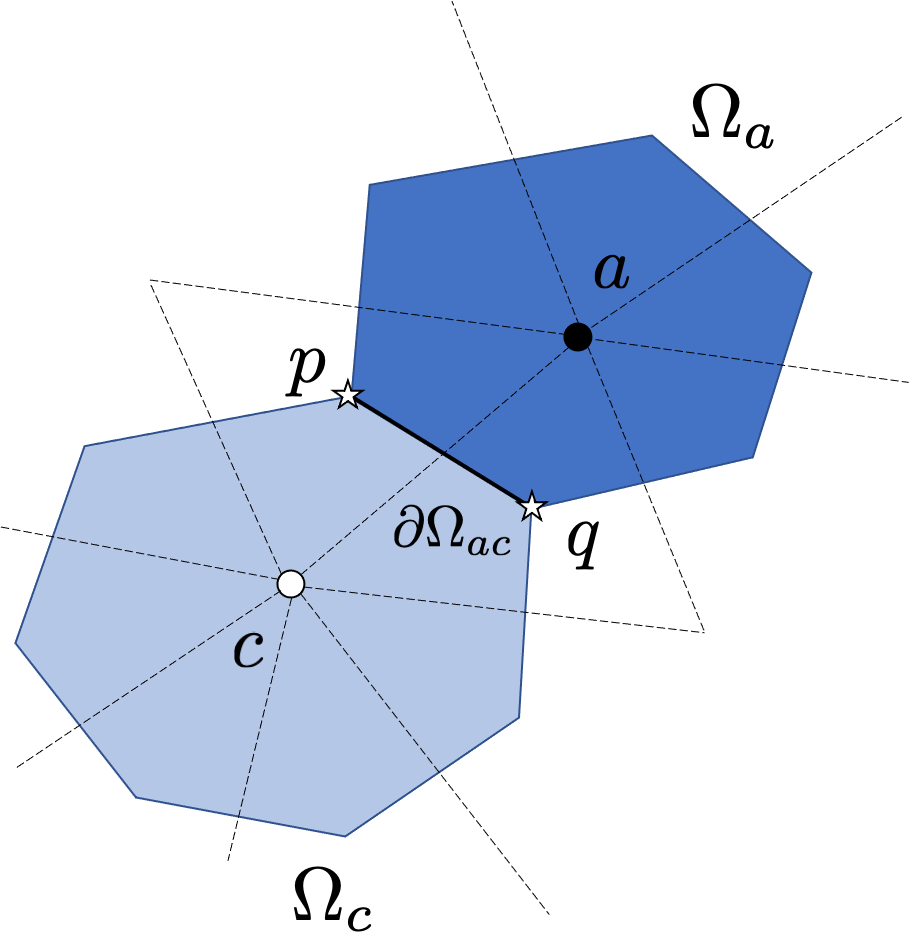}
		\end{minipage}
		\caption{In this Figure we depict the polygonal Voronoi-like cell $\Omega_c$ and, on the right, we highlight the edge $\partial\Omega_{ac}$ between the neighboring Voronoi cells  
			$\Omega_c$ (colored in light blue) and $\Omega_a$ (colored in blue). Note that the vertexes of the Delaunay mesh, i.e. the Voronoi generators ${a,b,c}$, are depicted with point symbols and the vertexes of the Voronoi ${p,q,r}$, i.e. the barycenters of the dual triangles, are depicted with star symbols. 
			\label{fig:voronoi-geo}}
	\end{figure}

\begin{figure}[h]	\centering\includegraphics[width=0.4\textwidth]{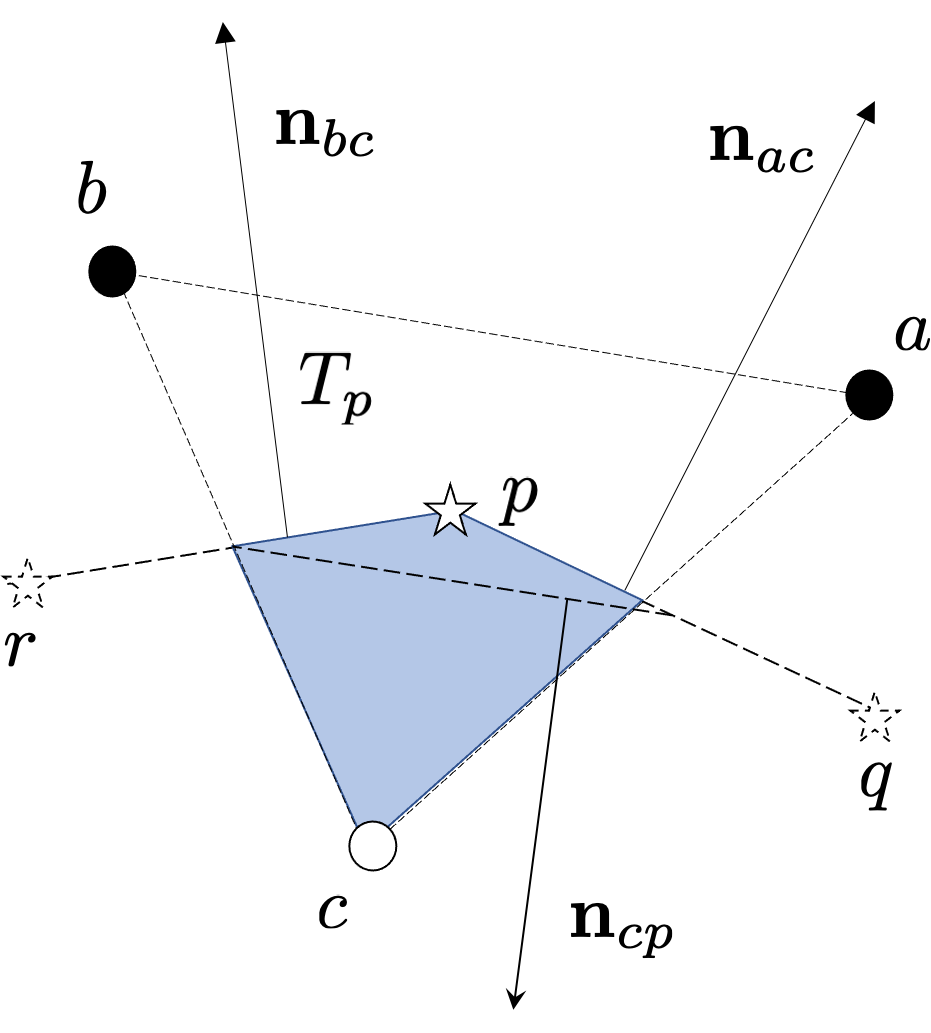}
		\caption{In this Figure we sketch the edge normal vectors $\n_{ac}$ and $\n_{bc}$ and the corner normal pointing from cell $\Omega_c$ to point $p$ which is denoted by $\n_{pc}$. 
			Note that in general the intersection of edge $\overline{ac}$ with $\overline{pq}$ is \textit{not} in the mid point of the edge $\overline{ac}$. Likewise for the intersection of the edge $\overline{bc}$ with $\overline{pr}$.}
		\label{fig:ncp-geo}
	\end{figure}

	A sketch of the notation used in this paper is provided in Figures~\ref{fig:voronoi-geo} and~\ref{fig:ncp-geo} where the set of elements attached to node $p$ corresponds to $\mathcal{C}_p=\left\{ a, b, c \right\}$ (which also defines the triangle $T_p$) and the set of points which are vertexes of $\Omega_c$ is $\mathcal{P}_c = \left\{ p, q, r, \dots  \right\}$.
	
	Finally, we define the corner normal pointing from cell $\Omega_c$ to point $p$ as 
	\begin{equation}
		\n_{pc} = \halb \left( \n_{ac} + \n_{bc} \right) = -\n_{cp}, 
		\label{eqn.cornernormal} 
	\end{equation}
	and we refer to Figure~\ref{fig:ncp-geo} for a visual representation. 
	
	\section{First order finite volume scheme using multidimensional corner fluxes}
	\label{sec.osher.o1}
	
	As usual in cell-centered finite volume schemes, the discrete data $\Q_c^n$ are assumed to represent the cell averages of the conserved variables $\Q$ 
	\begin{equation}
		\Q_c^n = \frac{1}{|\Omega_c|} \int \limits_{\Omega_c} \Q(\x,t^n) d\x, 
		\qquad 
		\textnormal{with}
		\qquad 
		|\Omega_c| = \int \limits_{\Omega_c}  d\x
		\label{eqn.cellaverage} 
	\end{equation}
	which we discretize on the primal Voronoi tessellation made by the polygons $\Omega_c$.
	Moreover, a classical first order finite volume scheme for~\eqref{eq.pde1} reads 
	\begin{equation}
		\Q_c^{n+1} = \Q_c^n - \frac{\Delta t}{|\Omega_c|} \sum_{\Omega_a \in \mathcal{N}_c} |\partial \Omega_{ac}| \, \hat{\f}_{ac},
		\label{eqn.fv} 
	\end{equation}
	where, 
	in addition to the previously defined notation, 
	$\hat{\f}_{ac}$ is the numerical flux in the normal direction $\hat{\n}_{ac}$ that approximates the integral
	\begin{equation}
		\hat{\f}_{ac} \approx \frac{1}{\Delta t \, |\partial \Omega_{ac}|} \int \limits_{t^n}^{t^{n+1}} \int \limits_{\partial \Omega_{ac}} \F(\Q(\x,t)) \cdot \hat{\n}_{ac} \, dS dt.  
		\label{eqn.fluxdef} 
	\end{equation} 
	In the context of classical edge-based \textit{two-point fluxes} and first order explicit Euler time integration the numerical flux in normal direction would take the simple form
	\begin{equation}
		\hat{\f}_{ac} = \hat{\f}_{ac}(\Q_c^n, \Q_a^n, \hat{\n}_{ac}),
		\label{eq.fv.2states}
	\end{equation}
	where as flux function $\hat{\f}_{ac}$ any classical 1d Riemann solver can be used (see e.g.~\cite{torobook}) and it is evaluated 
	at the midpoint of the edge $\partial\Omega_{ac}$.
	
	In this paper, instead, we employ \textit{point fluxes}, i.e. genuinely multidimensional Riemann solvers, as in~\cite{balsarahlle2d,balsarahllc2d,BalsaraMultiDRS,schneider2021multidimensional,MUSIC1,MUSIC2} and thus the use of the midpoint rule is no longer convenient. 
	Thus, in order to approximate the integral in~\eqref{eqn.fluxdef}, 
	one can employ the trapezoidal rule in such a way to include data at the vertexes $p$ of the control volume $\Omega_c$. 
	Thus, the flux across the edge becomes 
	\begin{equation}
		\hat{\f}_{ac} = \halb \left( \hat{\F}_{ac}^- + \hat{\F}_{ac}^+  \right) \cdot \hat{\n}_{ac} = \halb \left( \hat{\F}_{q}  + \hat{\F}_{p}  \right) \cdot \hat{\n}_{ac}, 
	\end{equation}
	with $\hat{\F}_{ac}^- = \hat{\F}_{q}$ and $\hat{\F}_{ac}^+ = \hat{\F}_{p}$ the numerical flux tensors evaluated in the two end points $q$ and $p$ of the edge $\partial \Omega_{ac}$. \RIIcolor{The corner flux $\hat{\F}_{p} = \hat{\F}_{p}(\Q_a^n, \Q_b^n, \Q_c^n)$ depends on the cell averages of the three elements $\Omega_a$, $\Omega_b$ and $\Omega_c$ surrounding node $p$ and its precise expression still needs to be determined. It will be detailed below in the following Sections \ref{sec.MultiDOsher} and \ref{sec.rd.o1}, which are the core of this paper.}

	Using the corner normal defined by~\eqref{eqn.cornernormal},  
	we can rewrite the finite volume scheme~\eqref{eqn.fv} also using the corner normals and the corner fluxes 
	instead of the edge normals and the edge fluxes as
	\begin{equation}
		\Q_c^{n+1} = \Q_c^n - \frac{\Delta t}{|\Omega_c|} \sum_{\Omega_a \in \mathcal{N}_c} |\partial \Omega_{ac}| \, \halb \left( \hat{\F}_{ac}^- + \hat{\F}_{ac}^+  \right) \cdot \n_{ac} = 
		\Q_c^n - \frac{\Delta t}{|\Omega_c|} \sum_{p \in \mathcal{P}_c}  \hat{\F}_{p}\cdot \n_{pc}, 
		\label{eqn.fv.corner} 
	\end{equation}
	\RIIcolor{with 
		$\hat{\F}_{p}\cdot \n_{pc} = \halb \left( 
		\hat{\F}_{ac}^- \cdot \hat{\n}_{ac} + 
		\hat{\F}_{bc}^+ \cdot \hat{\n}_{bc}
		\right) = \hat{\F}_{p} \cdot \halb \left( \hat{\n}_{ac} + \hat{\n}_{bc} \right) $, since the corner flux  $\hat{\F}_{p}$ is unique $\hat{\F}_{p} = \hat{\F}_{ac}^- = \hat{\F}_{bc}^+$. In \cite{Maire2007,Maire2008,Maire2009,Maire2010,Maire2011} the corner fluxes are also interpreted as a combination of the solution of two so-called half Riemann problems, see the previous references for details. }
	
	The last equation allows to introduce a multidimensional framework  generalizing the classical two states finite volume
	setting \eqref{eq.fv.2states}.  We will 
	directly define 
	the corner flux $\hat{\F}_{p}\cdot \n_{pc}$, which is now a function of $d+1$ states
	\begin{equation}\label{eqn.fv.corner.flux}
		\hat{\F}_{p}\cdot \n_{pc} =\hat{\F}_{p}\cdot \n_{pc}(\Q_a,\Q_b,\Q_c)
	\end{equation}
	The above  multidimensional numerical flux can be characterized with the following generalization of standard consistency and stability conditions \RIIcolor{(assuming either periodic or no-flux boundary conditions)}.
	\begin{description}
		\item[Consistency] A numerical corner flux is consistent if 
		\begin{equation}\label{eq:point_flux_consistency}
			\hat{\F}_{p}\cdot \n_{pc}(\Q,\Q,\Q) =
			\F(\Q)\cdot \n_{pc}\;,
		\end{equation}
		with $\F$ the physical flux.
		\item[Conservation] Local conservation at each corner  is equivalent to
		\begin{equation}\label{eq:point_flux_conservation}
			\sum_{c\in \mathcal{C}_p} \hat{\F}_{p}\cdot  \n_{pc}  =0\;.
		\end{equation}
		\item[Entropy stability]   Corner fluxes  verify a generalization of Tadmor's entropy conservation condition~\cite{tadmor} which reads 
		\begin{equation} \label{eq.shuffle0}
			\sum_{c\in\mathcal{C}_p}   \mathbf{W}_c^t \hat\F_p \cdot\n_{pc}  =   \sum_{c\in\mathcal{C}_p}    \boldsymbol{\Psi}_c  \cdot\n_{pc},
		\end{equation}
		having denoted by $\mathbf{W}$ the vector of entropy variables,
		and  with $\boldsymbol{\Psi}$ the entropy potential vector. The derivation of this condition is provided in appendix~\ref{app.entropy}. Any corner flux with more dissipation than the one of an entropy conservative flux, 
		defines and entropy stable scheme.
	\end{description}
	In the following  we will provide two 
	examples: \RIIcolor{the first can be seen as a multidimensional  generalizations 
		of the Osher-Solomon flux and the second of the 1d upwind flux difference splitting}.
	
	\section{A new genuinely multidimensional Osher-type flux}
	\label{sec.MultiDOsher}
	Here, the nodal numerical flux tensor $\hat{\F}_{p}$ in \eqref{eqn.fv.corner}   will be obtained via a novel extension of the Osher-Solomon flux~\cite{osherandsolomon} to the genuinely multidimensional case.
	Inspired by the ideas outlined in~\cite{OsherUniversal,OsherNC}, the genuinely multidimensional Osher-type flux, which provides a full \textit{numerical flux tensor} instead of a usual numerical flux in a special normal direction, is assumed to have the following form
	\begin{equation}
		\hat{\F}_{p} = \frac{1}{d+1} \sum_{c \in \mathcal{C}_p} \F\left(\Q_{c}^n\right) - \frac{h}{d+1} \, \left( \int \limits_{T_0} \diag\left(|\A_i(\psi(\bxi)|\right)\, d\bxi \, \right) \cdot \nabla_p^c \Q_c^n, 
		\label{eqn.multid.osher} 
	\end{equation}
	where the first term with the sum contains the centered part of the flux, while the term with the integral contains the numerical dissipation. 
	In the above formula~\eqref{eqn.multid.osher}, $\mathcal{C}_p$ is the set of cells $\Omega_c$ attached to a point $p$ (for example in Figure~\ref{fig:ncp-geo} $\mathcal{C}_p = \left\{c, a, b \right\}$) corresponding to the vertexes of the simplex $T_p$. 
	Instead, the unit simplex is denoted by $T_0$ and has the canonical vertexes $0$ and $1$ in 1d, 
	$(0,0)$, $(1,0)$ and $(0,1)$ in 2D, 
	and $(0,0,0)$, $(1,0,0)$, $(0,1,0)$ and $(0,0,1)$ in 3D. 
	Moreover, $h$ is a characteristic length scale and $\nabla_p^c \Q_c^n$ is a discrete gradient operator whose precise definitions will be given later. 
	
	We now recall that the absolute value of a \RIIcolor{diagonalizable} matrix $\mathbf{M}$ is defined as 
	\begin{equation}
		|\mathbf{M}| = \mathbf{R} \, |\boldsymbol{\Lambda}| \, \mathbf{R}^{-1},
	\end{equation} 
	with $\boldsymbol{\Lambda}$ the diagonal matrix of the eigenvalues of the matrix $\mathbf{M}$, and $\mathbf{R}$ the associated matrix of right eigenvectors. Since the system~\eqref{eq.pde1} is assumed to be hyperbolic, all matrices $\A_i(\psi(\bxi))$ in~\eqref{eqn.multid.osher}
	can be diagonalized, hence $|\A_i(\psi(\bxi))|$ can be computed.  
	
	While in the 1d Riemann solvers~\cite{OsherUniversal,OsherNC} it is employed a linear 1d \textit{segment path} $\psi = \psi(\xi_1)$, 
	here we use a piecewise linear \textit{surface} (manifold) inside the simplex $T_p$ given by 
	\begin{equation}
		\bpsi(\bxi) = \sum_{k = 1}^{d+1} \varphi_k(\bxi) \, \Q_{c(k,p)}^n,
		\label{eqn.multid.path}
	\end{equation}
	with $k \in \left\{1,2,\dots,d+1\right\}$ a local vertex index in the simplex element $T_p$, 
	$c=c(k,p)$ the global element number associated to node $p$ and vertex $k$ of the simplex $T_p$, 
	and the $\varphi_k(\bxi)$ being the classical P1 Lagrange basis functions given inside the unit simplex 
	and expressed in terms of reference coordinates $\bxi$. 
	In the $d=2$ case, they read 
	\begin{equation}
		\varphi_1(\bxi) = 1 - \xi_1 - \xi_2, \qquad 
		\varphi_2(\bxi) = \xi_1, \qquad 
		\varphi_3(\bxi) = \xi_2, 
	\end{equation}
	and for $d=3$ they read
	\begin{equation}
		\varphi_1(\bxi) = 1 - \xi_1 - \xi_2 - \xi_3, \qquad 
		\varphi_2(\bxi) = \xi_1, \qquad 
		\varphi_3(\bxi) = \xi_2, \qquad 
		\varphi_4(\bxi) = \xi_3.
	\end{equation}
	In the following we will denote the physical coordinates of the simplex $T_p$ by $\X_k = \x_{c(k,p)}$, $k \in \left\{1,2,\dots,d+1\right\}$. 
	
	The discrete gradient $\nabla_p^c \Q_c^n$ (which is a constant inside the simplex $T_p$)
	in physical coordinates $\x$  can be immediately computed from the linear surface $\bpsi(\bxi)$ as 
	\begin{equation}
		\nabla_p^c \Q_c^n = \frac{\partial \bpsi(\bxi)}{\partial \bxi}  \JJ^{-1} = %
		\sum_{k = 1}^{d+1} \frac{\partial \varphi_k(\bxi)}{\partial \bxi}  \, \Q_{c(k,p)}^n \, \JJ^{-1} = \frac{1}{|\Omega_p|} \sum_{c \in \mathcal{C}_p} \n_{cp} \, \Q_c^n.
		\label{eqn.nablaQ.osher}
	\end{equation}
	Here, $\JJ$ is the Jacobian matrix associated with the mapping from the unit simplex $T_0$ to the simplex $T_p$ 
	\begin{equation}
		\x = \x(\bxi) = \sum_{k=1}^{d+1} \varphi_k(\bxi) \, \X_k,
	\end{equation}
	and reads
	\begin{equation}
		\JJ = \frac{\partial \x}{\partial \bxi} = \sum_{k=1}^{d+1} \frac{\partial \varphi_k(\bxi)}{\partial \bxi} \, \X_k = 
		\left( \X_{d+1} - \X_1, \X_{d} - \X_1, \dots \X_2 - \X_1 \right). 
	\end{equation}	
	For the sake of simplicity we define the characteristic length scale as $h=\left( \frac{1}{d} |\JJ| \right)^{1/d}$.
	The integral in~\eqref{eqn.multid.osher} is approximated at the aid of numerical quadrature formulae of suitable accuracy, see e.g.~\cite{stroud}. 
	
	Inserting~\eqref{eqn.nablaQ.osher} into~\eqref{eqn.multid.osher} yields the following final expression
	\begin{equation}
		\hat{\F}_{p} = \frac{1}{d+1} \sum_{c \in \mathcal{C}_p} \F\left(\Q_{c}^n\right) - \frac{h}{d+1} \, \left( \int \limits_{T_0} \diag\left(|\A_i(\bpsi(\bxi)|\right) \, d\bxi \, \right) \cdot \left( \frac{1}{|\Omega_p|} \sum_{c \in \mathcal{C}_p} \n_{cp} \, \Q_c^n \right),
		\label{eqn.multid.osher.final}
	\end{equation}
	which together with the definition of the manifold~\eqref{eqn.multid.path} provides a closed expression for the numerical flux tensor $\hat{\F}_{p}$ in terms of the states $\Q_c^n$ adjacent to the vertex $\x_p$. 
	This completes the description of the new genuinely multidimensional Osher-type Riemann solver.

	\paragraph{The one-dimensional case as a special case of the general multidimensional setting} 
	The 1d case already treated in~\cite{OsherNC,OsherUniversal} is immediately included as a special case of the more general framework outlined in this Section. 
	Indeed, it is enough to take as basis function $\varphi_1(\bxi) = 1 - \xi_1$ and $\varphi_2(\bxi) = \xi_1$, 
	with the simplex $T_p$ for node $p=i+\halb$ being the interval $[x_i,x_i+1]$, i.e. $\mathcal{C}_p = \left\{ i, i+1 \right\}$, 
	and the 1d unit simplex being $T_0 = [0,1]$. 
	Thus the integration manifold reduces to the segment $\bpsi(\xi_1) = (1 - \xi_1) \Q_i^n + \xi_1 \Q_{i+1}^n = \Q_i^n + \xi_1 (\Q_{i+1}^n-\Q_i^n)$. 
	Moreover, we simply have $h=\Delta x$ hence 
	\begin{equation}
		h \nabla_p^c \Q_c = \Q_{i+1}^n - \Q_{i}^n,
	\end{equation} 
	which leads to the 1d numerical flux 
	\begin{equation}
		\hat{\f}_{i+\halb} = \halb \left( \f_1\left(\Q_i^n\right) + \f_1\left(\Q_{i+1}^n\right) \right) - 
		\halb \left( \int_{T_0} |\A_1(\bpsi(\xi_1))| \, d \xi_1 \right) \left( \Q_{i+1}^n - \Q_{i}^n \right),
	\end{equation} 
	already proposed in~\cite{OsherNC,OsherUniversal}, where the integral was also approximated at the aid of Gaussian quadrature formulae of suitable accuracy.
	
	\section{Corner fluxes coming from the residual distribution framework}
	\label{sec.rd.o1}
	
\begin{figure}[!h]
		\centering\includegraphics[width=0.6\textwidth]{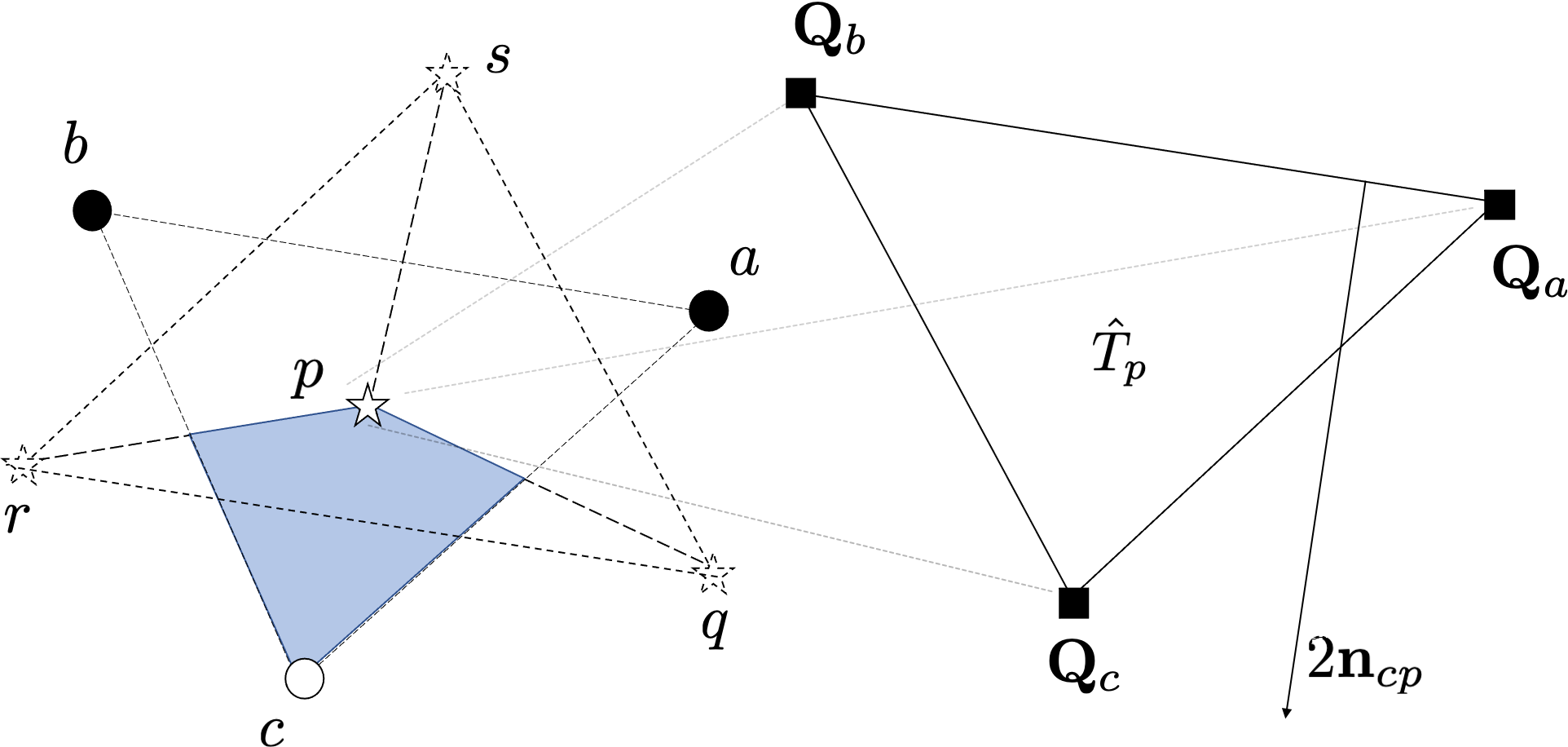}
		\caption{In this Figure we show the triangle $\hat T_p$ used to compute the residual distribution fluctuations $\boldsymbol{\phi}_p$ around each vertex $p$ of a polygonal element $c$.}
		\label{fig:tp-hat}
	\end{figure}

	The construction of the multidimensional Osher-type solver exploits \RIIcolor{two} geometrical properties of our Voronoi-like cells.
	In particular, it exploits the fact that every cell corner identifies uniquely $d+1$ cells which can be used to define 
	a unique local gradient using   classical linear interpolation on the \RIIcolor{triangle} $T_p$.
	\RIIcolor{Here,} we  take this idea one step further. Let us start again from the
	multidimensional prototype with corner fluxes already given in~\eqref{eqn.fv.corner}
	\begin{equation}\label{eq:lpc_npc1}
		\Q_c^{n+1} = \Q_c^n - \frac{\Delta t}{|\Omega_c|} \sum_{p\in\mathcal{P}_c}  \hat{\F}_{p}\cdot  \n_{pc},
	\end{equation}
	and focus on the definition of $\hat{\F}_{p}\cdot  \n_{pc}$ for a fixed $p$. 
	We will now exploit the point-wise conservation constraint on the numerical flux \eqref{eq:point_flux_conservation}. In particular, 
	we require that for each $c\in\mathcal{C}_p$ the numerical flux should be given by the internal one plus a fluctuation: 
	\begin{equation}\label{eq:point_flux_phi}
		\hat{\F}_{p}\cdot  \n_{pc} = \F_c\cdot  \n_{pc} + \boldsymbol{\phi}_{pc},
	\end{equation}
	where $\F_c = \F (\Q_c^n)$ and the $\boldsymbol{\phi}_{pc}$ are {\it nodal fluctuations} which must be defined.  
	The conservation condition \eqref{eq:point_flux_conservation} 
	requires that these fluctuations verify the constraint
	\begin{equation} 
		0= 
		\sum_{c\in\mathcal{C}_p}  \F_c\cdot  \n_{pc} +\sum_{c\in\mathcal{C}_p}   \boldsymbol{\phi}_{pc}.
	\end{equation}
	Setting 
	\begin{equation} \label{eq:point_flux_RD}	
		\boldsymbol{\phi}_p:=
		\sum_{c\in\mathcal{C}_p}  \F_c\cdot  \n_{cp}\;,
	\end{equation}
	Conservation at each node is equivalent to 
	\begin{equation} \label{eq:point_flux_RD.conservation}	
		\sum_{c\in\mathcal{C}_p}   \boldsymbol{\phi}_{pc}  =\boldsymbol{\phi}_p\;.
	\end{equation}
	Conservation at each corner is thus equivalent to  the requirement that the fluctuations $\boldsymbol{\phi}_{pc}$  should be a splitting of the residual $\boldsymbol{\phi}_p$
	defined in~\eqref{eq:point_flux_RD}.  The latter  
	can be seen as the  exact  integral of the divergence of a linear approximation of the flux, 
	over  the triangle $\hat T_p$ of Figure~\ref{fig:tp-hat},
	obtained by mirroring  the triangle $rqs$ w.r.t. the segment $[r,q]$.
	For this equivalence to be correct, the  values of the unknown at the nodes 
	of $\hat T_p$ should be  the traces of the $d+1$ discrete finite volume solution at node $p$.
	In the first order case  this trivially reduced to requiring that 
	the nodal values of   $\hat T_p$ should be   $\Q_a$, $\Q_b$, and $\Q_c$.
	Note, that in general $\hat T_p$ does not coincide exactly with $T_p$, unless the segments  $\overline{pq}$ and $\overline{ac}$  in Figure~\ref{fig:ncp-geo} intersect in correspondence of their mid-points, and similarly for  $\overline{pr}$ and $\overline{bc}$.\\

	The analogy presented above is a multidimensional generalization of the classical duality
	between finite volume methods and residual \RIIcolor{distributions} (see e.g.~\cite{abgrall2022hyperbolic} and references therein). 
	Besides conservation, all the properties discussed previously can be generalized to this setting. 
	In particular, the consistency requirement \eqref{eq:point_flux_consistency}, reduces in this case
	to the condition
	\begin{equation} \label{eq:point_fluct_consistency}	
		\phi_p(\Q,\Q,\Q) =0\;.
	\end{equation}
	In Appendix A we also show that the generalized Tadmor shuffle condition \eqref{eq:point_fluct_Tadmor}, is equivalent to the
	entropy conservation condition
	\begin{equation} \label{eq:point_fluct_Tadmor}	
		\sum\limits_{c\in\mathcal{C}_p} 
		\mathbf{W}_c^t\pmb{\phi}_{pc}=
		\sum\limits_{c\in\mathcal{C}_p} \G_c\cdot\n_{cp}\,,
	\end{equation}
	where the right hand side in the last expression is  the  integral over    $\hat T_p$  of the divergence of a linear approximation of the entropy  flux $\G(\Q)$. \\

	The problem of finding a splitting 
	\eqref{eq:point_flux_RD} is classical in the context of second order Residual Distribution~(RD) schemes~\cite{RD-ency,RD-ency2,HDR-MR}. 
	We will recall hereafter some   examples. 
	Before that let us note that in absence of a high order reconstruction we have  
	\begin{equation}\label{eq:lpc_npc_RD}
		\begin{aligned}
			\Q_c^{n+1} = & \Q_c^n - \frac{\Delta t}{|\Omega_c|} \sum_{p\in\mathcal{P}_c}  \hat{\F}_{p}\cdot  \n_{pc}\\
			=& \Q_c^n - \frac{\Delta t}{|\Omega_c|}  \F_c\cdot \sum_{p\in\mathcal{P}_c}  \n_{pc}  - \frac{\Delta t}{|\Omega_c|} \sum_{p\in\mathcal{P}_c}  \boldsymbol{\phi}_{pc} \\
			=& \Q_c^n   - \frac{\Delta t}{|\Omega_c|} \sum_{p\in\mathcal{P}_c}  \boldsymbol{\phi}_{pc}, 
		\end{aligned}
	\end{equation}
	having used the fact that $\sum_{p\in\mathcal{P}_c}  \n_{pc}=0$. This shows that \emph{at first order} the above is not only an analogy, but we have  a full equivalence between finite volumes with corner fluxes and RD schemes.
	
	\subsection{Example:  multidimensional Osher-Solomon  flux in fluctuation form}
	\label{sec.rd.o1-ex-Osher}
	
	Using~\eqref{eq:point_flux_phi}, we can readily deduce a fluctuation corresponding to the multidimensional Osher-type flux~\eqref{eqn.multid.osher}, that we have introduced in the previous section, which in 2D, and
	with the notation of Figures~\ref{fig:ncp-geo} and~\ref{fig:tp-hat}, reads
	\begin{equation} \label{eq:Osher_phi} 
		\boldsymbol{\phi}_{pc}^{\textsf{OS}} =
		\frac{1}{3}(\F_a-\F_c)\cdot\n_{pc}+ \frac{1}{3}(\F_b-\F_c)\cdot\n_{pc}- \frac{h}{3} \, \left[\left( \int \limits_{T_0} \diag\left(|\A_i(\bpsi(\bxi)|\right) \, d\bxi \, \right) \cdot \nabla_h \Q_p\right]\cdot\n_{pc}. 
	\end{equation}
	
	\subsection{Example: Rusanov method and an entropy stable \RIIcolor{splitting}}  
	\label{sec.rd.o1-ex-central}

	A Rusanov flux, or multidimensional variant  of local the Lax-Friedrich's scheme,  is obtained  with the following definition
	\begin{equation}\label{eq.rv-rd}
		\boldsymbol{\phi}_{cp}^{\textsf{Rv}}= \dfrac{1}{3}\boldsymbol{\phi}_p  + \alpha_p ( \Q_c - \bar \Q_p ),
	\end{equation}
	where $\bar \Q_{p}$ is the average 
	\begin{equation}\label{eq.rv-rda}
		\bar \Q_{p} =\dfrac{1}{d+1}\sum_{c\in \mathcal{C}_p}   \Q_c.
	\end{equation}
	Usually the constant $\alpha_p$ is chosen as the mesh size times the largest spectral radius of the system matrix projected 
	in the normal directions. In other words 
	\begin{equation}\label{eq.rv-rdb} 
		\alpha_p = h\times \max( \rho(K_{\n_{cp}}(\Q_c)), \rho(K_{\n_{ap} }(\Q_{a})) , \rho(K_{\n_{bp} }(\Q_{b}))   ) .
	\end{equation}
	The Rusanov splitting  has many interesting   properties,  including
	a discrete maximum principle when using (convex combinations of) the explicit Euler time integration, and  
	energy stability for linear symmetric systems (see e.g.~\cite{HDR-MR,RD-ency}). It is a very dissipative method, but its form is useful
	to construct schemes with desirable properties, 
	or design corrections for existing schemes.

	As an example, we can use the above splitting to obtain entropy conservative or stable schemes.   Inspired by the previous definition, and by recent work by R. Abgrall~\cite{abgrall2018general}, we can easily provide an extension of the previous
	definition  which is  strictly entropy conservative/stable.   
	In particular, instead of  using \eqref{eq.rv-rdb} to evaluated $\alpha_p$, we now compute it by explicitly imposing the entropy conservation condition \eqref{eq:point_fluct_Tadmor}:
	\begin{equation} \label{eq.EC1}
		\alpha_p  = \dfrac{  \sum_{c\in\mathcal{C}_p}    \mathbf{G}  \cdot\n_{cp} -   \overline{\mathbf{W}}_p^t \boldsymbol{\phi}_p  }{(\sum_{c\in\mathcal{C}_p} \mathbf{W}_c^t ( \mathbf{Q}_c - \overline{\mathbf{Q}}_p))}
	\end{equation}
	where $\overline{\mathbf{W}}_p = \sum_{c\in\mathcal{C}_p}\mathbf{W}_c/(d+1)$, and where   the denominator can be easily shown to be always positive definite for all convex entropies. 
	
	An entropy stable flux can be obtained by adding   to the above  any 
	positive contribution, or by adding to $\boldsymbol{\phi}_{pc}$
	an extra dissipation operator as e.g. the one  of the multidimensional Osher-Solomon method. 
	
	\subsection{Example:  a multidimensional upwind flux using the N scheme}
	\label{sec.rd.o1-ex-N}
	
	In the scalar case, the N scheme is Roe's multidimensional upwind optimum first order positive coefficient
	method on triangles~\cite{Roe:87,Roe:90,Roe92}.  Its extension to systems of conservation laws can be performed  using
	an exact linearization of the flux Jacobians, as in Roe's one dimensional method, and  either 
	seeking a decomposition of the system in scalar waves, or via 
	a compact matrix formulation~\cite{deconinck93,vanderweide}.
	For cases in which    an exact   mean value linearization of the flux Jacobians  is not available or impractical, 
	the authors of~\cite{crd02,rcd05} have introduced 
	a more general formulation compatible with    the use of arbitrary linearized states  for the system  Jacobians. 
	Using the notation introduced in~\eqref{eq.pde3}, the  N scheme splitting is  defined by 
	\begin{equation}\label{eq.N-rd}
		\boldsymbol{\phi}_{cp}^{\text{N}}=\K_{\n_{cp}}^+(\Q_c - \tilde \Q_p )\,.
	\end{equation}
	As in~\cite{crd02,rcd05} the tilde state $\tilde \Q_p$ is obtained by imposing the consistency condition~\eqref{eq:point_flux_RD}:
	\begin{equation}\label{eq.N-rd1}
		\tilde \Q_p =  \Q^+_{p}  - \mathbf{N}^{-1}\boldsymbol{\phi}_p, 
	\end{equation}
	where  $\mathbf{N}$  denotes   the matrix
	\begin{equation}\label{eq.N-rd2}
		\mathbf{N} =  \sum_{c\in\mathcal{C}_p}\K_{\n_{cp}}^+ ,
	\end{equation}
	and where
	\begin{equation}\label{eq.N-rd3}
		\Q^+_{p}:= \mathbf{N}^{-1} \sum_{c\in\mathcal{C}_p}\K_{\n_{cp}}^+  \Q_c.
	\end{equation}
	The multidimensional upwind N scheme is the optimum first order scheme, where  ``optimum''  means it has  
	the 
	least amount of numerical dissipation among first order positive coefficient discretizations for scalar advection~\cite{Roe:87,Roe:90,Roe92}.  
	The N scheme verifies 
	a discrete maximum principle for scalar conservation laws, when using (convex combinations of) the explicit Euler time integration.
	For linear symmetric systems, it can also be shown to be energy stable (see e.g.~\cite{HDR-MR,barth02,RD-ency} for all the previous statements).
	\RIIcolor{For nonlinear problems, only asymptotic entropy inequalities have been established (see e.g.~\cite{barth02}), which do not preclude violations of this condition.} 
	Indeed, when considering high Mach number flows, while the low dissipation of the N scheme has shown high potential to resolve complex flows with contact waves and boundary layers on coarse meshes~\cite{CiCP-13-479}, ad-hoc fixes for unwanted features such as the carbuncle phenomenon and expansion shocks may be necessary~\cite{GARICANOMENA201643,rd-efix}.
	In this work, we also introduce a correction  inspired by the  Rusanov splitting:
	\begin{equation}\label{eq.N-carbuncle}
		\boldsymbol{\phi}_{cp}^{\text{N-C}}
		=\boldsymbol{\phi}_{cp}^{\text{N}} + \mu_{\text{C}}(\Q) 
		(\Q_c -\overline{\Q}_p)
	\end{equation}
	where $\mu_{\text{C}}(\Q) $  is a small  extra artificial viscosity   only 
	active in strong shocks. 
	This modification has in particular no effect on contact waves.
	The precise definition of $\mu_{\text{C}}$, as well as the shock detection
	follows almost exactly~\cite{RODIONOV2017308}. 
	
	\paragraph{Roe's   one-dimensional upwind flux difference splitter as a special case of the  multidimensional setting} In one space dimension, looking at a cell $i$, the computational triangle $\hat{T}_p$ reduces to 
	a segment at the interface $i\pm 1/2$. The corner normals reduce to $\pm 1$. If we focus on the interface $i+1/2$,  we can easily see that  the $\K_{\n_{cp}}$ matrices reduce to $\pm A_1(\Q_{i+1/2}))$, and in particular for cell $i$ the matrix $\K_{\n_{cp}}^+$ is replaced by $[-A_1(\Q_{i+1/2}))] ^+= -A_1^-(\Q_{i+1/2})$. One easily shows that the nodal fluctuation becomes (see~\cite{abgrall2022hyperbolic} for example) 
	\begin{equation}
		\phi_{i+1/2} =  \f_1\left(\Q_{i+1}^n\right)-\f_1\left(\Q_{i}^n\right)\\
	\end{equation}
	while  the matrix  $\mathbf{N}$ in \eqref{eq.N-rd2} reduces to
	\begin{equation}
		\mathbf{N}_{i+1/2}  = -A_1^-(\Q_{i+1/2}) + A_1^+(\Q_{i+1/2}) = |A_1(\Q_{i+1/2})|.
	\end{equation}
	All things assembled, and after few manipulations we obtain the numerical fluxes
	\begin{equation}
		\begin{aligned}
			\hat{\f}_{i+\halb} = &
			\f_1\left(\Q_i^n\right)
			- \dfrac{A_1^-(\Q_{i+1/2})|A_1(\Q_{i+1/2})|^{-1}}{2} (\f_1\left(\Q_{i+1}^n\right)-\f_1\left(\Q_{i}^n\right))
			\\
			=& \halb \left( \f_1\left(\Q_i^n\right) + \f_1\left(\Q_{i+1}^n\right) \right) - \dfrac{A_1(\Q_{i+1/2})|A_1(\Q_{i+1/2})|^{-1}}{2} (\f_1\left(\Q_{i+1}^n\right)-\f_1\left(\Q_{i}^n\right)).
		\end{aligned}
	\end{equation} 
	When $A_1(\Q_{i+1/2})$ is evaluated averaging the  Roe parameter, the above is   precisely  Roe's flux~\cite{roe1d}.

	\section{Higher order extension using the ADER approach with CWENO reconstruction}
	\label{sec.highorder}
	
	\subsection{Spatial CWENO reconstruction}
	\label{ssec.cweno}
	
	In the context of better than second order accurate  finite volume schemes a  reconstruction procedure is necessary to obtain better than piecewise constant approximations of the data within each cell. For this purpose one needs to compute a high order accurate but at the same time non-oscillatory polynomial representation $\w_h^n(\x,t^n)$ of the solution $\Q(\x,t^n)$ for each polygon $\Omega_c$, starting from the values of the cell averages $\mathbf{Q}_a^n$ in $\Omega_c$ and its neighbors contained in a so-called reconstruction stencil $\mathcal{S}_c^s$. For a polynomial of degree $M \geq 0$ in two space dimensions the number of degrees of freedom per conserved variable is $\mathcal{M} = (M+1)(M+2)/2$. Moreover, as already mentioned in~\cite{barthlsq}, on general unstructured meshes it is necessary to use more than $\mathcal{M}$ cells to get a stable numerical scheme. Furthermore, one reconstruction stencil alone is not enough in order to obtain an essentially non-oscillatory scheme, but ENO/WENO reconstruction techniques~\cite{shu_efficient_weno,harten_eno} relying on multiple stencils, as those introduced in~\cite{abgrall_eno,HuShuVortex1999,DumbserKaeser07,ZhangShu3D,MixedWENO2D,MixedWENO3D,ADERCWENO} for unstructured meshes, are needed in order to circumvent the Godunov theorem which states that there are no better than first order accurate linear monotone schemes~\cite{godunov}. In this paper we follow in particular the successful approach of CWENO schemes introduced by Lewy, Puppo, Semplice and Russo \textit{et al.} in a series of papers in~\cite{CWENO1,CWENO2,CWENO3,cravero2018cweno,ADERCWENO}, and which was for the first time applied on general unstructured polygonal meshes in~\cite{gaburro2020high,gaburro2021unified}. Here, we closely follow the work outlined in~\cite{ADERCWENO,gaburro2020high}, but for the sake of completeness and to make the paper self-contained for the reader, we report here the entire algorithm.
	
	The reconstruction starts from the computation of a so-called \textit{central polynomial} $\Popt$ of degree $M$. 
	In order to define $\Popt$ in a robust manner, following~\cite{barthlsq,DumbserKaeser07,ADERCWENO}, we consider a stencil
	$\mathcal{S}_c^0$ which contains a total number of $n_e= f \cdot \mathcal{M}$ control volume and which includes cell $\Omega_c$ and its neighbors
	\begin{equation}
		\mathcal{S}_c^0 = \bigcup \limits_{k=1}^{n_e} \Omega_{a(k)},
		\label{stencil}
	\end{equation}
	with the safety factor $f \ge 1.5$. 
	Stencil $S_c^0$ is constructed by starting from the central polygon $\Omega_c$ and then adding recursively of neighbors and neighbors of elements that have been already included in the stencil, until the desired total number of $n_e$ elements is reached.  
	The polynomial $\Popt(\x,t^n)$ is then defined by imposing that its average on each cell $\Omega_{a} \in \mathcal{S}_c^0$ matches the known cell average  
	${\mathbf{Q}}^{n}_{a}$. 
	Since $n_e > \mathcal{M}$, \RIIcolor{this obviously leads} to an overdetermined linear system, which is solved using a constrained least-squares technique introduced in~\cite{DumbserKaeser06b}.
	Here, the polynomial $\Popt$ has exactly the cell average ${\mathbf{Q}}^{n}_{c}$ on the polygon $\Omega_c$ and matches all other cell averages ${\mathbf{Q}}^{n}_{a}$ for $a \neq c$ only in a least-square sense. The polynomial $\Popt$ is expressed in terms of spatial reconstruction basis functions $\psi_\ell(\x,t^n)$ of degree $M$,  
	\begin{equation}
		\label{eqn.recpolydef} 
		\Popt(\x,t^n) = \sum \limits_{\ell=0}^{\mathcal{M}-1} \psi_\ell(\x,t^n) \hat{\mathbf{p}}^{n}_{\ell,i}.  
	\end{equation}
	Since there is no simple mapping to universal reference elements for general polygonal meshes, as in~\cite{gaburro2020high,gaburro2025high} we use the following reconstruction basis functions that employ rescaled and shifted Taylor monomials
	\begin{equation} 
		\label{eq.psi_spatial}
		\psi_\ell(\x,t^n) |_{P_i^n} = \frac{(x - x_c)^{p_\ell}}{p_\ell! \, h_c^{p_\ell}} \, \frac{(y - y_c)^{q_\ell}}{q_\ell! \, h_c^{q_\ell}}, \qquad 
		\ell = 0, \dots, \mathcal{M}-1, \quad \ 0 \leq p_\ell + q_\ell \leq M,
	\end{equation} 
	with $h_c$ the radius of the circumcircle and  $\xbin$ the barycenter of the polygon $\Omega_c$. The  expansion coefficients $\hat{\mathbf{p}}^{n}_{\ell,i}$ in~\eqref{eqn.recpolydef} therefore represent the rescaled derivatives  of the Taylor expansion about $\xbin$.
	
	According to what was written above we now require integral conservation in the central stencil $\mathcal{S}_c^0$ as follows: 
	\begin{equation} 
		\frac{1}{|\Omega_a|} \int \limits_{\Omega_a} \Popt(\x,t^n) d\x = \Q_a^n, \qquad 
		\forall \Omega_a \in \mathcal{S}^0_c.
		\label{CWENO:Popt} 
	\end{equation} 
	and the integrals appearing in~\eqref{CWENO:Popt} are conveniently computed in each polygon by simply summing up the contributions over all sub-triangles of a polygon and on the sub-triangles we employ suitable conical products of the one-dimensional Gauss-Jacobi formula, see~\cite{stroud} for details. 
	
\begin{figure}[!b]
		\centering
		\includegraphics[width=0.33\textwidth]{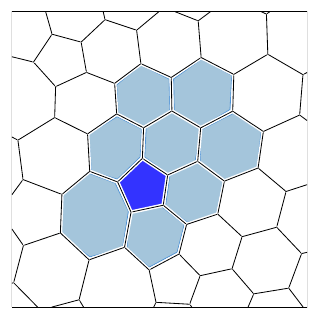}%
		\includegraphics[width=0.33\textwidth]{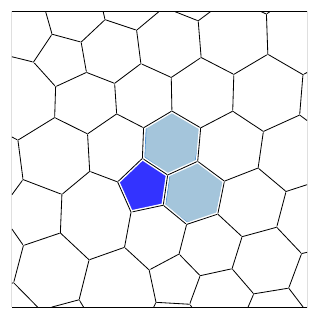}%
		\includegraphics[width=0.33\textwidth]{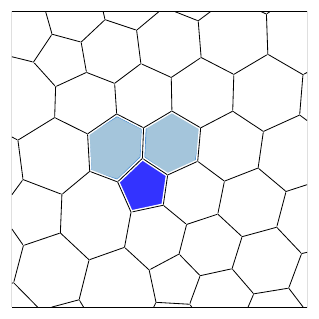}\\
		\includegraphics[width=0.33\textwidth]{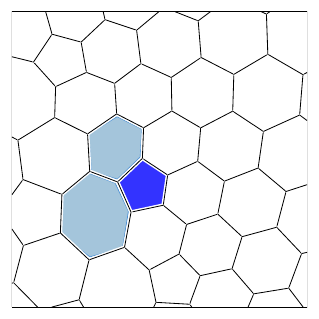}%
		\includegraphics[width=0.33\textwidth]{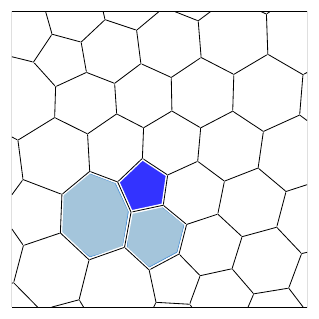}%
		\includegraphics[width=0.33\textwidth]{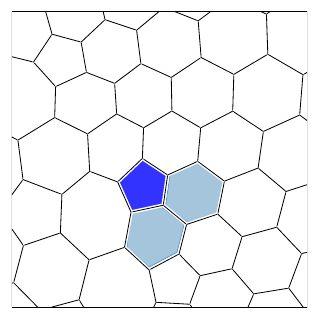}%
		\caption{Illustration of the different reconstruction stencils used for the CWENO reconstruction of order three ($M=2$) with a safety factor of $f=1.5$ for a pentagonal central element $\Omega_c$ in blue. Top-left: central stencil needed for the reconstruction of $\Popt$ (in light blue). In the other panels we report the 5 sectorial stencils containing the element itself and two consecutive neighbors needed to reconstruct piecewise linear sectorial polynomials.}
		\label{fig.cweno_stencils}
	\end{figure}

	In order to obtain a nonlinear reconstruction procedure, which is needed to avoid spurious oscillations in the vicinity of shock waves, the CWENO method  uses also other polynomials of lower degree. In particular, in this work we employ piecewise linear reconstructions for the lower order reconstruction polynomials.  
	Given a polygon $\Omega_c$ with $N_c$ Voronoi neighbors contained in the set $\mathcal{V}(\Omega_c)$, we construct $N_c$ interpolating polynomials of degree $M^s=1$ referred to as sectorial polynomials. In order to do this we consider $N_c$ stencils $S_c^s$ with $s\in[1,N_c]$, each of them containing three cells, namely the central polygon $\Omega_c$ and two consecutive neighbors in the set
	$\mathcal{V}(\Omega_c)$.  
	An example of such reconstruction stencils is reported for the case $M=2$ in Figure~\ref{fig.cweno_stencils}. 
	
	For each stencil $\mathcal{S}_c^s$ we compute a linear polynomial $\mathbf{P}_s(\x,t^n) \in \mathbb{P}_1$ by solving the reconstruction systems 
	\begin{equation} 
		\frac{1}{|\Omega_a|} \int \limits_{\Omega_a} \mathbf{P}_s(\x,t^n) d\x = \Q_a^n, \qquad 
		\forall \Omega_a \in \mathcal{S}^s_c 
		\label{CWENO:Ps}
	\end{equation} 
	on all sectorial stencils $\mathcal{S}^s_c$. 
	These reconstruction equations are \textit{not} overdetermined and thus do not require the constrained least squares procedure. 
	Following the general framework introduced in~\cite{cravero2018cweno}, we select a set of positive coefficients $\lambda_0,\ldots,\lambda_{N_p}$ such that
	\begin{equation}
		\sum_{s=0}^{N_c}\lambda_s=1
		\label{eqn.sumCWENO}
	\end{equation}
	and we define a new auxiliary polynomial 
	\begin{equation}
		\label{CWENO:P0}
		\mathbf{P}_0(\x,t^n) = \frac{1}{\lambda_0}\left( \left( \sum_{s=0}^{N_c} \lambda_s \right) \Popt(\x,t^n) - \sum_{s=1}^{N_c} \lambda_s \mathbf{P}_s(\x,t^n) \right) \in\mathbb{P}_M,
	\end{equation}
	so that the linear combination of the polynomials $\mathbf{P}_0,\ldots,\mathbf{P}_{N_c}$ with the coefficients $\lambda_0,\ldots,\lambda_{N_c}$ is equal to $\Popt$ and conservation is ensured. Specifically, we consider the linear weights used in~\cite{ADERCWENO}, namely $\lambda_0=10^5$ for $\mathcal{S}_i^0$ and $\lambda_s=1$ for the sectorial stencils. These weights are then automatically renormalized in order to sum to unity, according to the requirement~\eqref{eqn.sumCWENO}, see also Equation~\eqref{CWENO:P0} which contains the automatic rescaling.
	All polynomials $\mathbf{P}_s$ with $s\in[0,N_c]$ are then non-linearly combined with each other as in classical WENO schemes, see e.g.~\cite{shu_efficient_weno,HuShuVortex1999,DumbserKaeser06b}. 
	We therefore obtain the final nonlinear CWENO reconstruction polynomial $\w_{h}(\x,t^n)$ in $\Omega_c$ as
	\begin{equation}
		\label{CWENO:Prec}
		\w_{h}(\x,t^n) = \sum_{s=0}^{N_c} \omega_s \mathbf{P}_s(\x,t^n), \qquad \x \in \Omega_c, 
	\end{equation}
	where the normalized {\em nonlinear weights} $\boldsymbol{\omega}_s$ read
	\begin{equation}
		\label{eqn.weights}
		\boldsymbol{\omega}_s = \frac{\tilde{\boldsymbol{\omega}}_s}{ {\sum \limits_{m=0}^{N_c} \tilde{\boldsymbol{\omega}}_m} }, 
		\qquad \textnormal{ with } \qquad 
		\tilde{\boldsymbol{\omega}}_s = \frac{\lambda_s}{\left(\boldsymbol{\sigma}_s + \epsilon \right)^r}. 
	\end{equation} 
	In the equation above the non-normalized weights $\tilde{\boldsymbol{\omega}}_s$ depend on the linear weights $\lambda_s$ and the oscillation indicators $\boldsymbol{\sigma}_s$. 
	The parameters $\epsilon=10^{-14}$ and $r=4$ are chosen according to~\cite{DumbserKaeser06b,ADERCWENO,gaburro2020high}. 
	As in~\cite{gaburro2020high}, the oscillation indicators $\boldsymbol{\sigma}_s$ appearing in~\eqref{eqn.weights} are simply given by the quadratic form 
	\begin{equation}
		\boldsymbol{\sigma}_s = \sum_{l} \left(\hat{\mathbf{p}}^{n,s}_{l,i}\right)^2.
		\label{eqn.OI}
	\end{equation}  
	Note that in smooth regions of the domain the central stencil with the optimal order polynomial is preferred due to the large linear weight $\lambda_0$, while in the vicinity of \RIIcolor{discontinuities} the reconstruction switches to the lower order sectorial polynomials to reduce oscillations. 
	
	\subsection{High order element-local space-time ADER predictor}
	\label{ssec.ader}
	
	Unlike the original ADER schemes of Toro and Titarev~\cite{toro3,toro2005ader}
	here we employ the element-local space-time DG predictor as introduced in~\cite{DumbserEnauxToro,dumbser2008unified}. This local space-time predictor solution is valid \textit{locally} inside each $\Omega_c$ for the current time interval $[t^n,t^{n+1}]$ and is given by high order piecewise polynomials 
	$\q_h(\x,t)$ of degree $M$ in space-time which read 
	\begin{equation} 
		\q_h = \q_h(\x, t) = \sum_{\ell=0}^{\mathcal{Q}-1} \theta_\ell (\x, t) \hat{\q}_\ell^n, \qquad (\x,t) \in \Omega_c \times [t^n,t^{n+1}], \qquad \mathcal{Q} = \mathcal{M} \cdot (M+1),  
		\label{eqn.qh}
	\end{equation} 
	with $\theta_\ell(\x, t)$ the \textit{modal space--time} basis already introduced in~\cite{gaburro2020high,GaburroRicchiutoPrimitive}:
	\begin{equation}
		\label{eq.theta}
		\theta_\ell(x,y,t)|_{C_i^n} = \frac{(x - x_c)^{p_\ell}}{{p_\ell}! \, h_c^{p_\ell}} \, \frac{(y - y_c)^{q_\ell}}{{q_\ell}! \, h_c^{q_\ell}}
		\, \frac{(t - t^n)^{q_\ell}}{{q_\ell}! \, h_c^{q_\ell}}, 
		\qquad \ell = 0, \dots, \mathcal{Q}-1, 
		\quad 0 \leq p_\ell + q_\ell + r_\ell \leq M. \\
	\end{equation}  
	To obtain a weak formulation of the governing PDE system in space-time we multiply the system of conservation laws~\eqref{eq.pde1} by a test function 
	$\theta_k$,  integrate over the space-time control volume $\mathcal{C}_c = \Omega_c \times [t^n,t^{n+1}]$ and substitute the discrete solution $\q_h(\x,t)$ for $\mathbf{Q}$. This yields 
	\begin{equation}
		\int \limits_{\mathcal{C}_c} \theta_k(\x,t) \de{\q_h}{t} \, d\x dt  + 
		\int \limits_{\mathcal{C}_c} \theta_k(\x,t) \nabla \cdot \F(\q_h) \, d\x dt  = 0.
		\label{eqn.PDEweak1}
	\end{equation}
	In order to introduce the initial condition in this weak formulation of the local Cauchy problem, given by the spatial reconstruction polynomial obtained via the CWENO procedure outlined in the previous section, we integrate the first term by parts in time and then use upwinding in time (causality principle), thus obtaining
	\begin{equation}
		\int \limits_{\Omega_c} \theta_k(\x,t^{n+1}) \q_h(\x,t^{n+1}) \, d\x 
		-\int \limits_{\mathcal{C}_c} \frac{\partial}{\partial t} \theta_k(\x,t) \q_h \, d\x dt  + 
		\int \limits_{\mathcal{C}_c} \theta_k(\x,t) \nabla \cdot \F(\q_h) \, d\x dt  = \int \limits_{\Omega_c} \theta_k(\x,t^n) \w_h(\x,t^n) \, d\x.
		\label{eqn.PDEweak2}
	\end{equation}
	The element-local space-time predictor $\q_h^n$ can then be easily computed via a simple fixed-point iteration, the convergence of which has been shown in~\cite{frontiers2020}. For details, the reader is referred to~\cite{dumbser2008unified,jackson2017eigenvalues,frontiers2020}.  
	
	\subsection{Final high order ADER-CWENO scheme based on multidimensional Riemann solvers}
	
	For better than second order accurate schemes in space and time we need higher order accurate numerical quadrature formulas to integrate the fluxes along the element edges. As done in~\cite{BalsaraMultiDRS,ADERdivB,LagrangeMHD} we use the Keplerian barrel rule\footnote{This quadrature formula is widely known under the name Simpson rule, but it had already been found more than 100 years earlier by Johannes Kepler in 1615 in his work \textit{Nova stereometria doliorium vinariorum, in primis Austriaci }, see~\cite{kepler}, hence we deliberately prefer to refer to the original reference.}, which is a Gauss-Lobatto quadrature formula with 3 quadrature points, two being located at the vertexes of an edge and one in the edge barycenter. As such, the high order version of the ADER-CWENO scheme based on vertex fluxes reads
	\begin{equation}
		\Q_c^{n+1} = \Q_c^n 
		- \frac{1}{3} \frac{\Delta t}{|\Omega_c|} \sum_{p \in \mathcal{P}_c}  \n_{pc} \cdot  \hat{\F}_{p}
		- \frac{2}{3} \frac{\Delta t}{|\Omega_c|} \sum_{\Omega_a \in \mathcal{N}_c} |\partial \Omega_{ac}| \, \hat{\f}_{ac},
		\label{eqn.finafv.highorder}
	\end{equation}
	with $\hat{\F}_{p}$ the time-averaged \textit{vertex fluxes} computed in the vertex-extrapolated states of the high order space-time predictor solution $\q_h$ and $\hat{\f}_{ac}$ the time-averaged \textit{edge fluxes} computed at the aid of a conventional 1d Riemann solver. Throughout this paper we employ the 1d Riemann solvers
	to which the multidimensional fluxes reduce to:
	the 1d Osher solver~\cite{OsherNC,OsherUniversal} for
	the flux of Section~\ref{sec.MultiDOsher}, and the Roe method for the
	multidimensional upwind N scheme solver.
	Both choices are consistent with the aim of this paper of using \textit{complete} Riemann solvers anywhere. 
	More precisely, to compute the time average flux $\hat{\f}_{ac}$ 
	with the appropriate order of accuracy, we evaluate it as 
	\begin{equation}
		\hat{\f}_{ac} = \sum_j \omega_j \, \hat{\f}_{ac}\left( \q_h(\x_{ac}^-,t^n+\tau_j \Delta t), \q_h(\x_{ac}^+,t^n+\tau_j \Delta t), \hat{\n}_{ac} \right), 
		\label{eqn.f1d.ho}
	\end{equation}
	where $\tau_j$ and $\omega_j$ are the points and weights of a classical Gauss-Legendre quadrature formula of suitable order in 1d. 
	Two quadrature points would be enough both for third and fourth order of accuracy, however we employ $M+1$ points to respect the underlying original structure of our ADER solver. 
	Furthermore, $\q_h^- = \q_h(\x_{ac}^-,t^n+\tau_j \Delta t)$ and $\q_h(\x_{ac}^+,t^n+\tau_j \Delta t)$ are the boundary-extrapolated states that enter the 1d Riemann solver and $\x_{ac}^-$ and $\x_{ac}^+$ is the edge barycenter seen from the left element $\Omega_c$ and from the right element $\Omega_a$, respectively. Likewise, the vertex flux must be computed in the vertex-extrapolated states. In case of the new multidimensional Osher flux we obtain
	\begin{eqnarray}
		\hat{\F}_{p} &=& \sum_j \omega_j \, \frac{1}{d+1} \sum_{c \in \mathcal{C}_p} \F\left(\q_h(\x_{cp},t^n+\tau_j \Delta t)\right) - \nonumber \\ 
		&& \sum_j \omega_j \, \frac{h}{d+1} \, \left( \int \limits_{T_0} \diag\left(|\A_i(\bpsi(\bxi)|\right) \, d\bxi \, \right) \cdot \left( \frac{1}{|\Omega_p|} \sum_{c \in \mathcal{C}_p} \n_{cp} \, \q_h(\x_{cp},t^n+\tau_j \Delta t) \right). \qquad 
		\label{eqn.multid.osher.ho} 
	\end{eqnarray}
	Here, $\x_{cp}$ is the vertex coordinate $\x_p$ taken from inside control volume $\Omega_c$,
	i.e. the predictor $\q_h(\x_{cp},t^n+\tau_j \Delta t) $ is also taken from within 
	$\Omega_c$ and is computed in the location $\x_p$ of the vertex. 
	
	The other node fluxes presented in this paper can be computed in a similar manner, simply replacing $\Q_c^n$ with the vertex-extrapolated value $\q_h(\x_{cp},t^n+\tau_j \Delta t)$ and using Gaussian quadrature in time to compute the time-averaged flux.

	\section{Numerical results}
	\label{sec.tests}

	In order to numerically verify the capabilities of our multidimensional Riemann solvers and the applicability of our strategy for their high order extension in space and time, 
	we test them on a series of benchmarks typical of the compressible Euler equations of gas-dynamics. 
	
	The Euler equations can be cast in the form~\eqref{eq.pde1} in the two-dimensional framework by choosing   
	\begin{equation}
		\label{eulerTerms}
		\Q = \left( \begin{array}{c} \rho   \\ \rho u  \\ \rho v \\ \rho E \end{array} \right)\!, \quad
		\mathbf{F} = \left( \begin{array}{ccc}  \rho u       & \rho v        \\ 
			\rho u^2 + p & \rho u v          \\
			\rho u v     & \rho v^2 + p      \\ 
			u(\rho E + p) & v(\rho E + p)   
		\end{array} \right). 
	\end{equation}
	Here, the vector of conserved variables $\Q$ involves the fluid density $\rho$, the momentum density vector $\rho \v=(\rho u, \rho v)$ and the total energy density $\rho E$. 
	The fluid pressure $p$ is related to the conserved quantities $\Q$ using the equation of state for an ideal 
	gas    
	\begin{equation}
		\label{eqn.eos} 
		p = (\gamma-1)\left(\rho E - \frac{1}{2} \rho \mathbf{v}^2 \right)\!, 
	\end{equation}
	where $\gamma$ ($\gamma=1.4$ in our tests) is the ratio of specific heats so that the speed of sound takes the form $c=\sqrt{\frac{\gamma p}{\rho}}$.
	
	Let us notice that for all the benchmarks we have deliberately chosen to work with coarse meshes in order to better appreciate the influence of the Riemann solvers on the final result and to make the results quickly replicable.

	\subsection{Order of convergence of our high order ADER-CWENO schemes based on multidimensional Riemann solvers}
	\begin{table*}[!t]
		\caption{Isentropic Shu-type vortex. Numerical convergence results for our ADER-CWENO scheme based on the multidimensional Osher and the N scheme as Riemann solvers and employed from the first order ($M=0$) to the fourth order ($M=3$) of accuracy.
			The error norms refer to the variable $\rho$ at time $t=1.0$ in the $L_2$ norm.} \vspace{2pt}
		\label{tab.orderOfconvergenceDG_shu}
		\centering
		\begin{tabular}{cc|cc|cc} 
			\hline 
			&&  \multicolumn{2}{c|}{MultiD Osher} & \multicolumn{2}{c}{N scheme} \\
			\hline           
			& $h$ & $\epsilon_{L_2}(\rho)$ & $\mathcal{O}(\rho)$ &  $\epsilon_{L_2}(\rho)$ & $\mathcal{O}(\rho)$ \\
			\hline   
			\multirow{5}{6em}{$M=0\rightarrow\mathcal{O} 1$} 
			& 4.10E-2 & 6.62E-2  & -    & 2.91E-2  & -   \\ 
			& 2.47E-2 & 4.11E-2  & 0.94 & 1.74E-2  & 1.02 \\
			& 1.69E-2 & 2.85E-2  & 0.96 & 1.18E-2  & 1.01 \\
			& 1.24E-2 & 2.12E-2  & 0.97 & 8.68E-3  & 1.01 \\		
			\hline        
			\multirow{5}{6em}{$M=1\rightarrow \mathcal{O}2$}         
			& 6.11E-2 & 2.87E-3  & -    & 2.26E-3  &   - \\		
			& 4.10E-2 & 1.10E-3  & 2.40 & 9.86E-4  & 2.07 \\		
			& 3.07E-2 & 5.64E-4  & 2.33 & 5.37E-4  & 2.12\\ 
			& 2.47E-2 & 3.55E-4  & 2.10 & 3.47E-4  & 1.99 \\		
			\hline        
			\multirow{5}{6em}{$M=2\rightarrow \mathcal{O}3$}         
			& 6.11E-2 & 2.76E-3  & -    & 2.45E-3  & -    \\		        
			& 4.10E-2 & 7.30E-4  & 3.31 & 6.19E-4  & 3.42 \\		
			& 3.51E-2 & 4.86E-4  & 2.70 & 4.07E-4  & 2.77 \\		
			& 3.07E-2 & 3.23E-4  & 3.02 & 2.72E-4  & 3.00 \\ 		
			\hline 
			\multirow{5}{6em}{$M=3\rightarrow \mathcal{O}4$}         
			& 1.19E-1 & 2.17E-2  & -    & 2.29E-2  & -    \\ 		
			& 8.08E-2 & 4.50E-3  & 4.10 & 4.65E-3  & 4.10 \\		
			& 6.11E-2 & 1.54E-3  & 3.83 & 1.54E-3  & 3.96  \\		        
			& 4.10E-2 & 2.67E-4  & 4.36 & 3.12E-4  & 3.97 \\	      
			\hline 
		\end{tabular}   
	\end{table*}
	
	To check the order of accuracy of our multidimensional Riemann solvers embedded in the strategy presented in Section~\ref{sec.highorder}, 
	we employ a smooth isentropic vortex test case inspired to the one 
	proposed in the example 3.3. of~\cite{hu1999weighted}, 
	which is a stationary equilibrium of the Euler system of equations. 

	\RIIcolor{The computational domain is the square $\Omega=[0;10]\times[0;10]$ where we have imposed transmissive boundary conditions (using wall boundary conditions leads to essentially the same results, since the velocity at the boundaries tends to zero).}
	The initial condition is given by some perturbations $\delta$ that are superimposed onto a homogeneous background field $\mathbf{V}_0=(\rho,u,v,p)=(1,0,0,1)$
	(note that w.r.t.~\cite{hu1999weighted} here we do not translate the vortex). 
	The perturbations for density and pressure are
	\begin{equation}
		\label{rhopressDelta}
		\delta \rho = (1+\delta T)^{\frac{1}{\gamma-1}}-1, \quad \delta p = (1+\delta T)^{\frac{\gamma}{\gamma-1}}-1, 
	\end{equation}
	with the temperature fluctuation $\delta T = -\frac{(\gamma-1)\epsilon^2}{8\gamma\pi^2}e^{1-r^2}$ and the 
	vortex strength is $\epsilon=5$.
	The velocity field is affected by the following perturbations
	\begin{equation}
		\label{ShuVortDelta}
		\left(\begin{array}{c} \delta u \\ \delta v  \end{array}\right) = \frac{\epsilon}{2\pi}e^{\frac{1-r^2}{2}} \left(\begin{array}{c} 
			-(y-5) \\ \phantom{-}(x-5)  \end{array}\right).
	\end{equation}
	
	When the simulation reaches the final time $t=1$, we compute the discrete $L^2$ error norms by evaluating the difference between the reconstructed numerical solution $\w_h$ and the exact solution $\Q_{\text{ref}}$ as
	\[
	\label{eq.error_norm}
	\epsilon_{L_2} = || \w_h -\Q_{\text{ref}} ||_{L^2}
	= \sum_{c=1}^{N_\omega} \left ( \int_{\Omega_c} | \w_h(\x,t_f) - \Q_{\text{ref}}(\x, t_f)| \, dx dy \right)^{\frac{1}{2}},
	\]
	using Gaussian quadrature rules of appropriate order.
	The obtained orders of convergence, which are achieved as expected, are reported in Table~\ref{tab.orderOfconvergenceDG_shu} for both the multidimensional Osher and the N scheme Riemann solvers employed in our CWENO-ADER scheme of order from one to four.

	\subsection{Contact and shear waves}
	In this section we concentrate on the study of isolated contact waves and shear waves.
	
	\paragraph{Steady contact wave} 
	\begin{figure}[!b]
		\centering
		\includegraphics[width=0.333\linewidth]{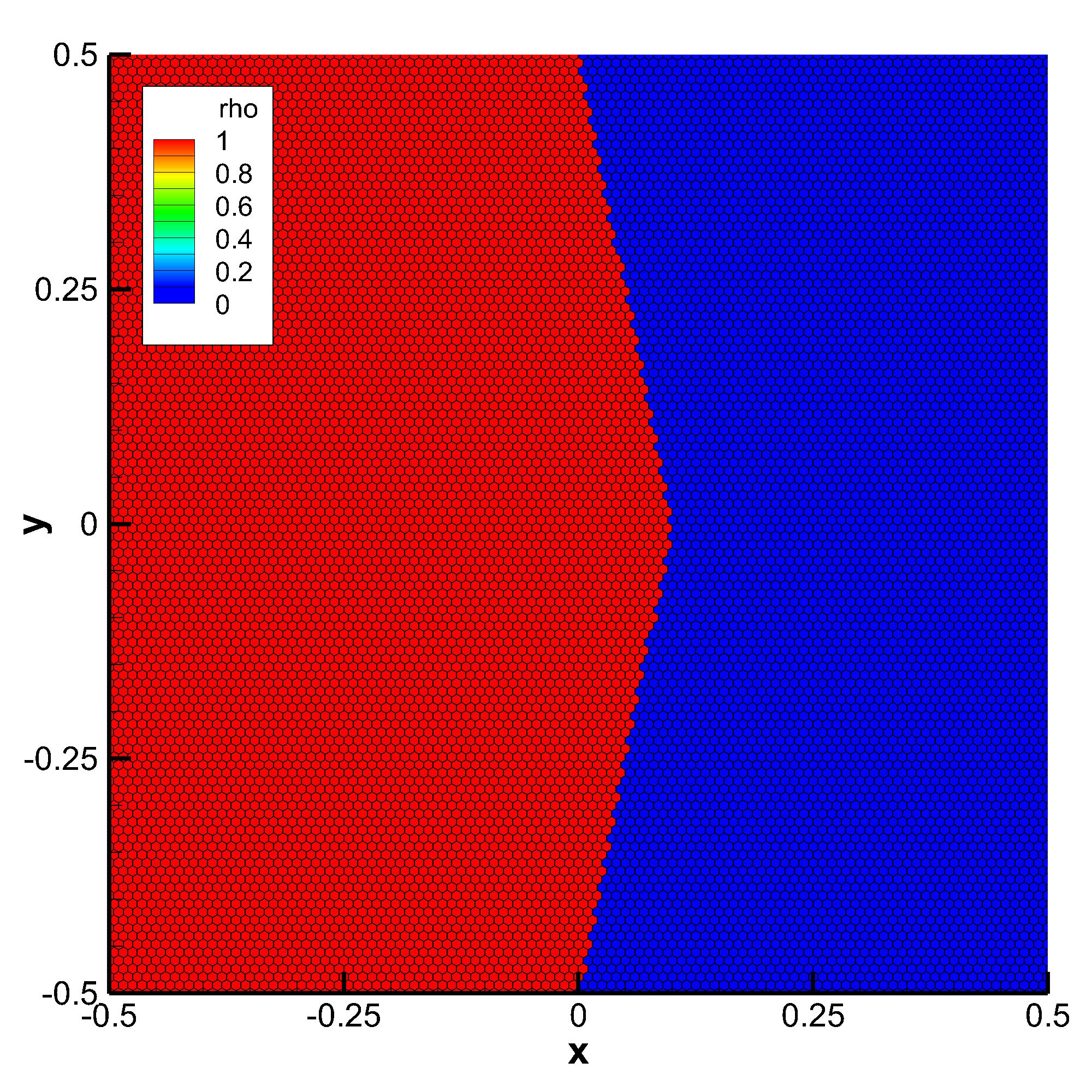}%
		\includegraphics[width=0.333\linewidth]{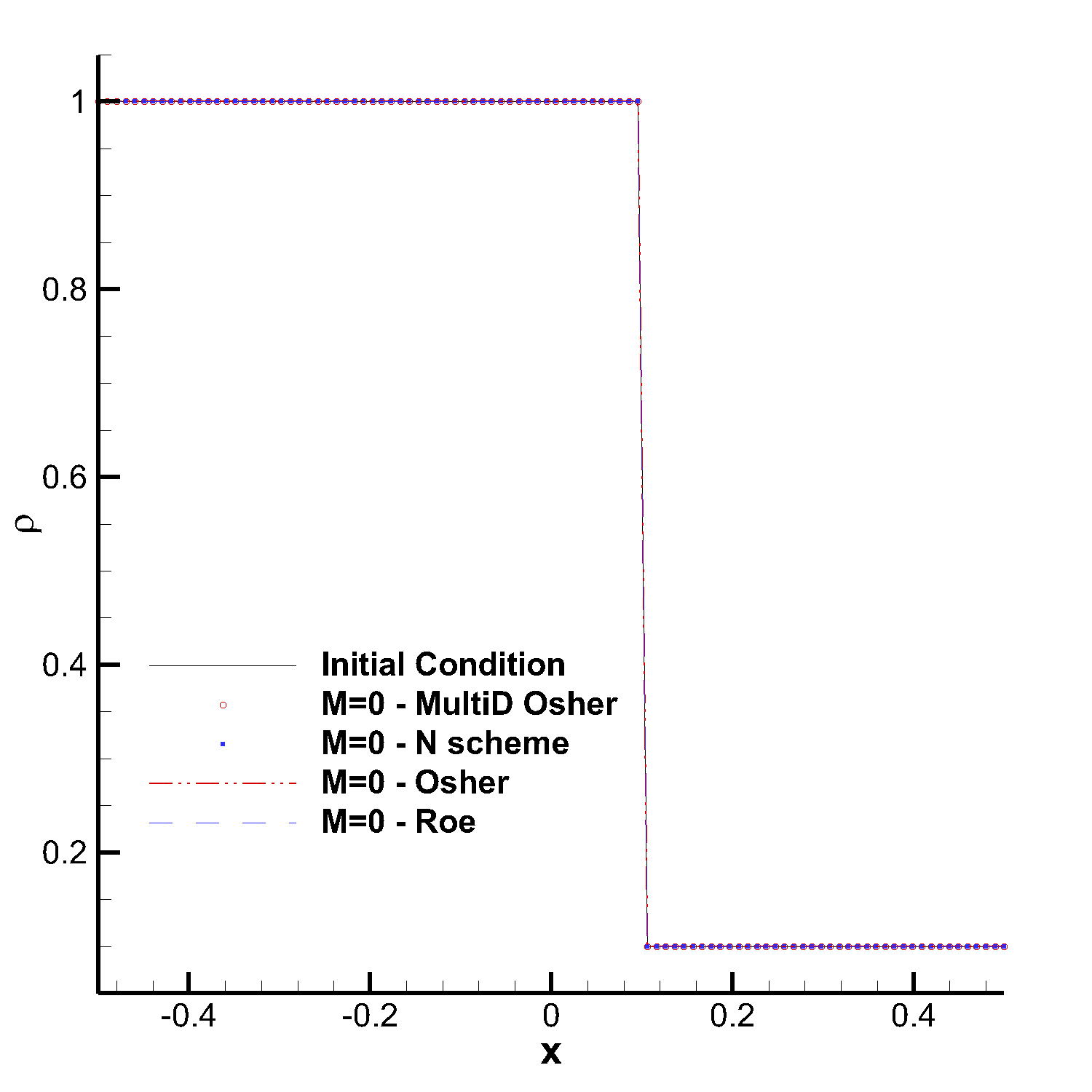}%
		\includegraphics[width=0.333\linewidth]{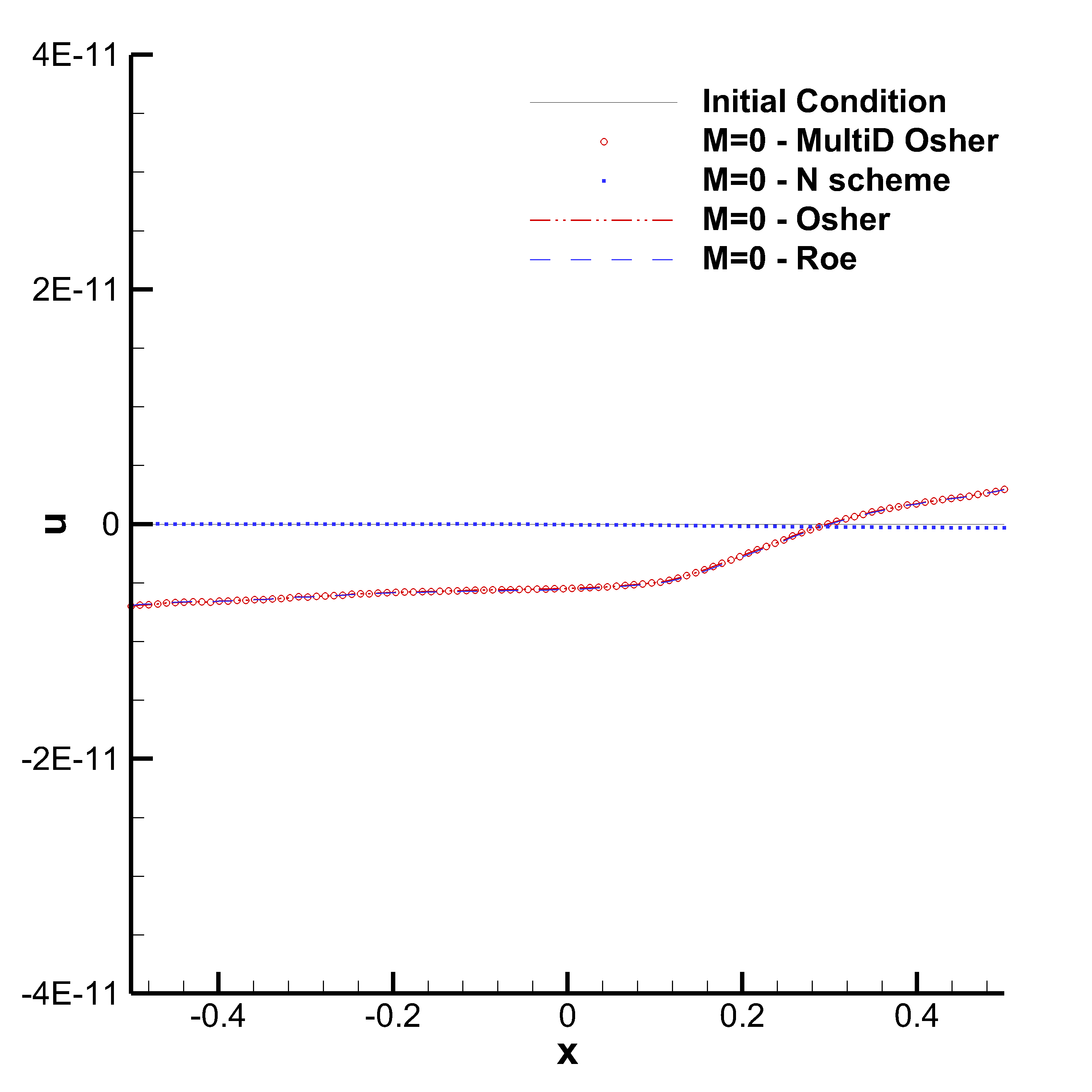}%
		\vspace{-10pt}
		\caption{Steady contact wave. On the left we show the discontinuous initial density profile, which is not aligned with the mesh. Then we show the 1d cut along $y=0$ of the final density \RIIcolor{$\rho$} (middle) and the \RIIcolor{final $x-$component of the velocity} \RIIcolor{$u$} (right) obtained with our novel multidimensional solvers and some classical 1d solvers. Note that the not aligned steady contact wave is maintained with machine precision.}
		\label{fig:contact}
	\end{figure}
	
	We first notice that both our novel multidimensional Riemann solvers as well as the classical complete Riemann solvers that we have used for comparison reasons (Osher and Roe solvers) preserve steady contact waves up to \RIIcolor{machine precision} even when the position of the initial discontinuity is not aligned with the mesh. 
	As numerical evidence of this property, we show in Figure~\ref{fig:contact} the discontinuous initial density profile of our test case
	\begin{equation}
		(\rho, u, v, p)(\x) = 
		\begin{cases} 
			(1.0, 0, 0, 1.0 )    &  \ \text{ if } y > 0 \ \&  \ y < -5x+0.5,  \\	
			(0.1, 0, 0, 1.0    ) &  \ \text{ if } y > 0 \ \&  \ y \ge -5x+0.5,\\
			(1.0, 0, 0, 1.0 )    &  \ \text{ if } y < 0 \ \& \ y \ge 5x-0.5,  \\	        
			(0.1, 0, 0, 1.0    ) &  \ \text{ if } y < 0 \ \& \ y < 5x-0.5,   \\		
		\end{cases}    
	\end{equation}
	simulated on the square $[-0.5;0.5]\times[-0.5;0.5]$ with periodic boundary conditions on top and bottom boundaries, 
	and the $1d$ cut along $y=0$ of the density and the total velocity ($V=\sqrt{u^2 + v^2}$) obtained at time t=1.0 with our first order numerical schemes, which coincide with the initial condition up to machine precision.

	\paragraph{Isolated shear waves}

\begin{figure}[!b]
		\centering
		\includegraphics[width=0.5\linewidth]{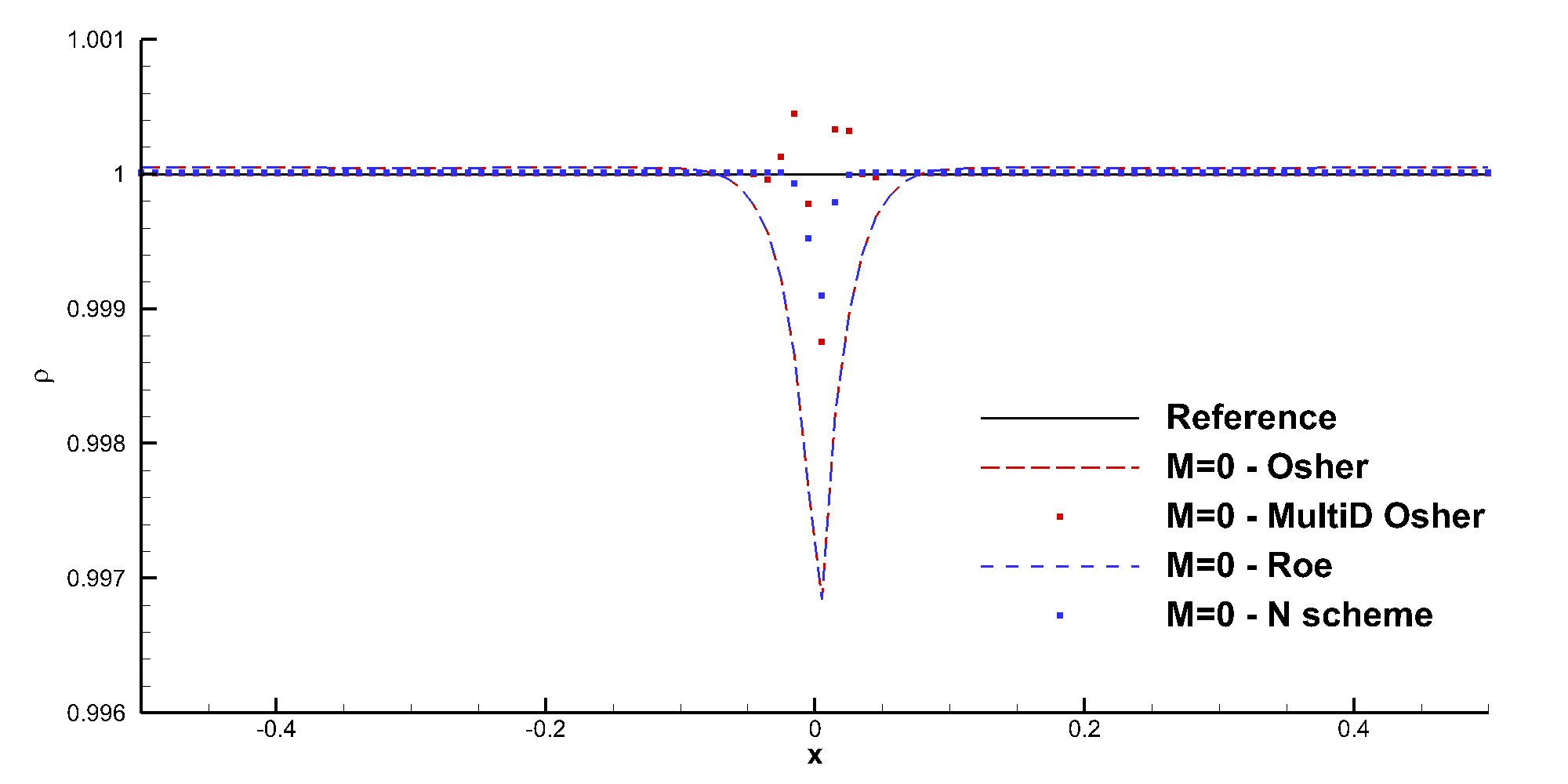}%
		\includegraphics[width=0.5\linewidth]{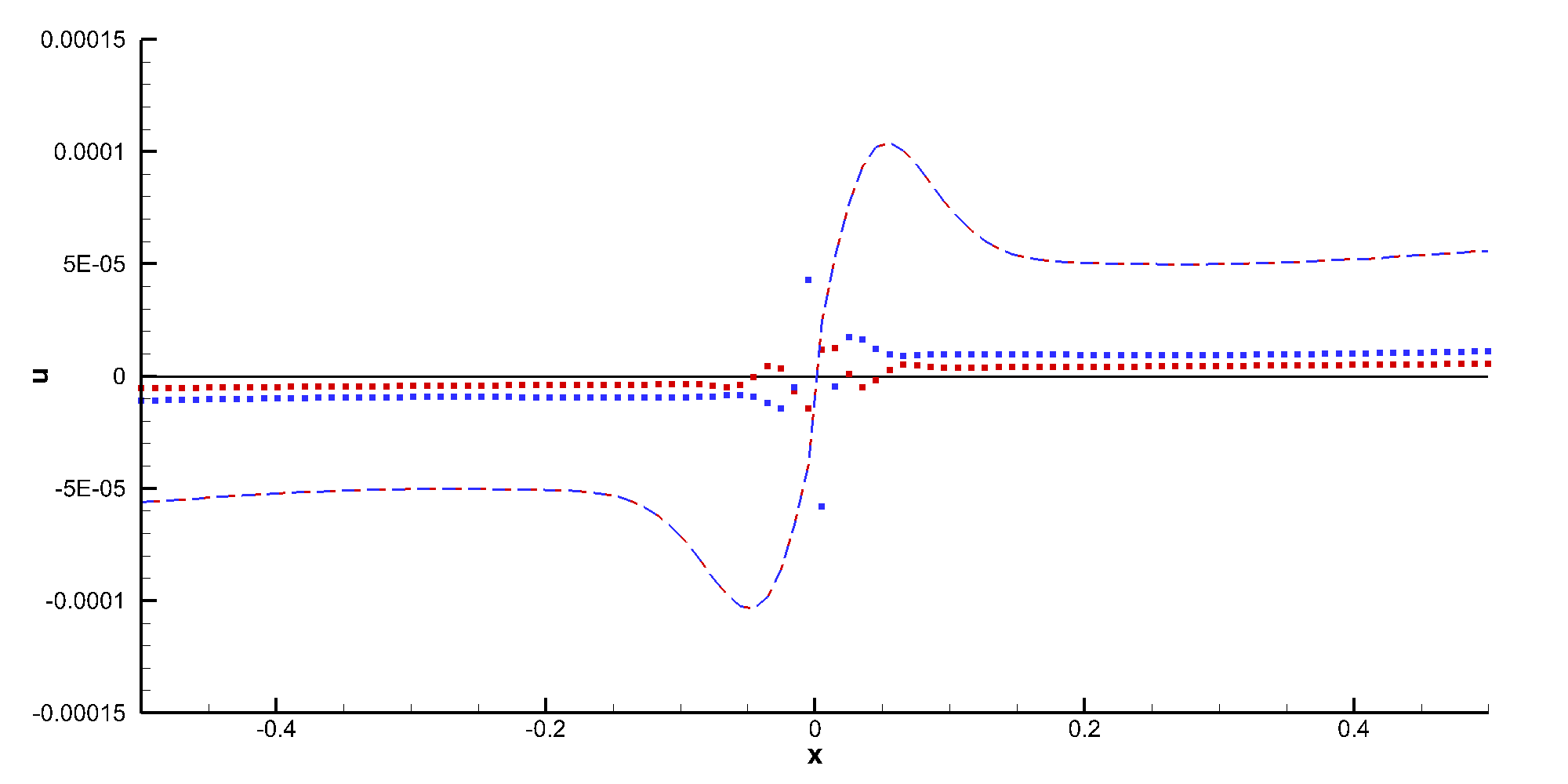}\\
		\includegraphics[width=0.5\linewidth]{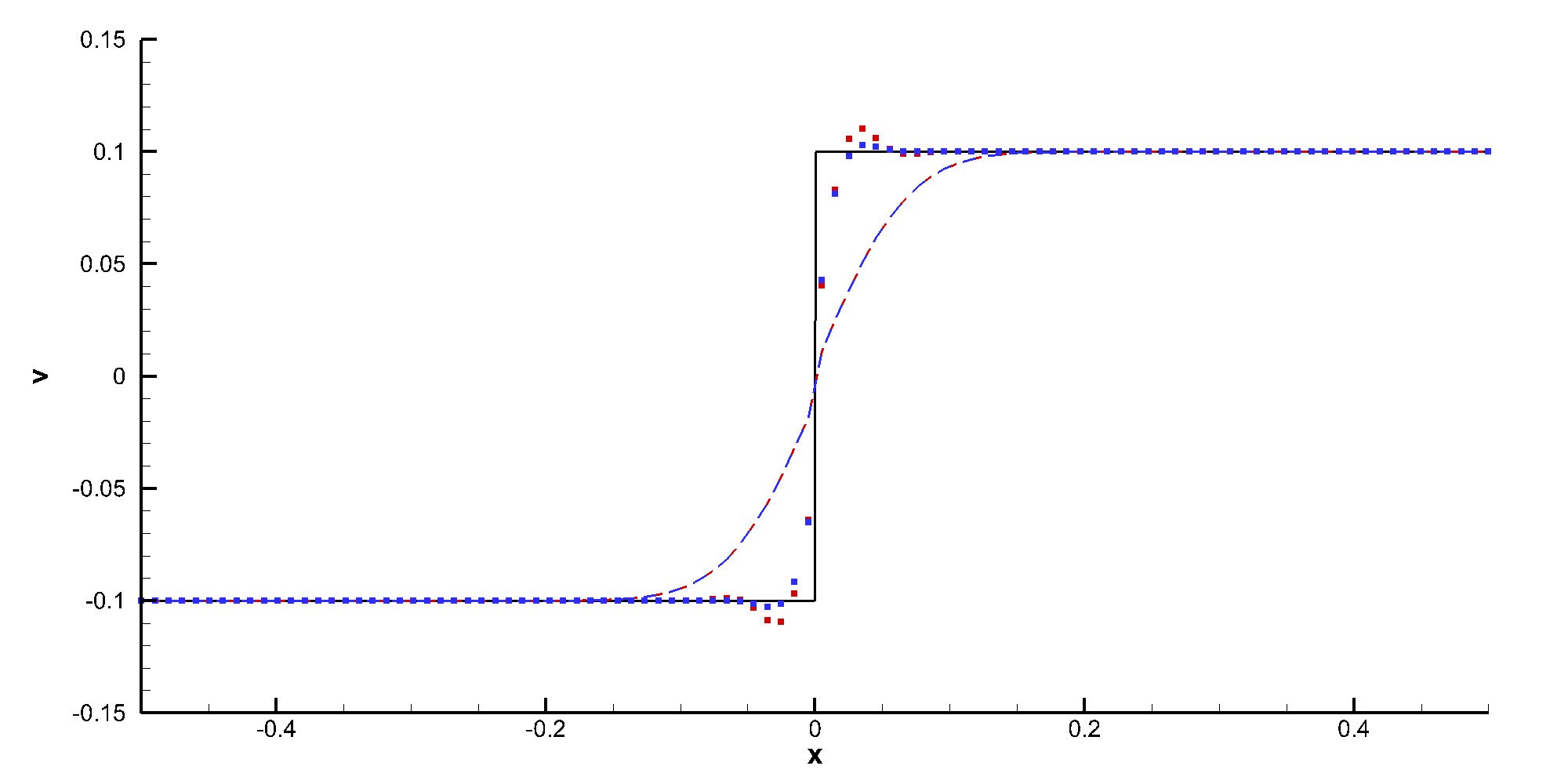}%
		\includegraphics[width=0.5\linewidth]{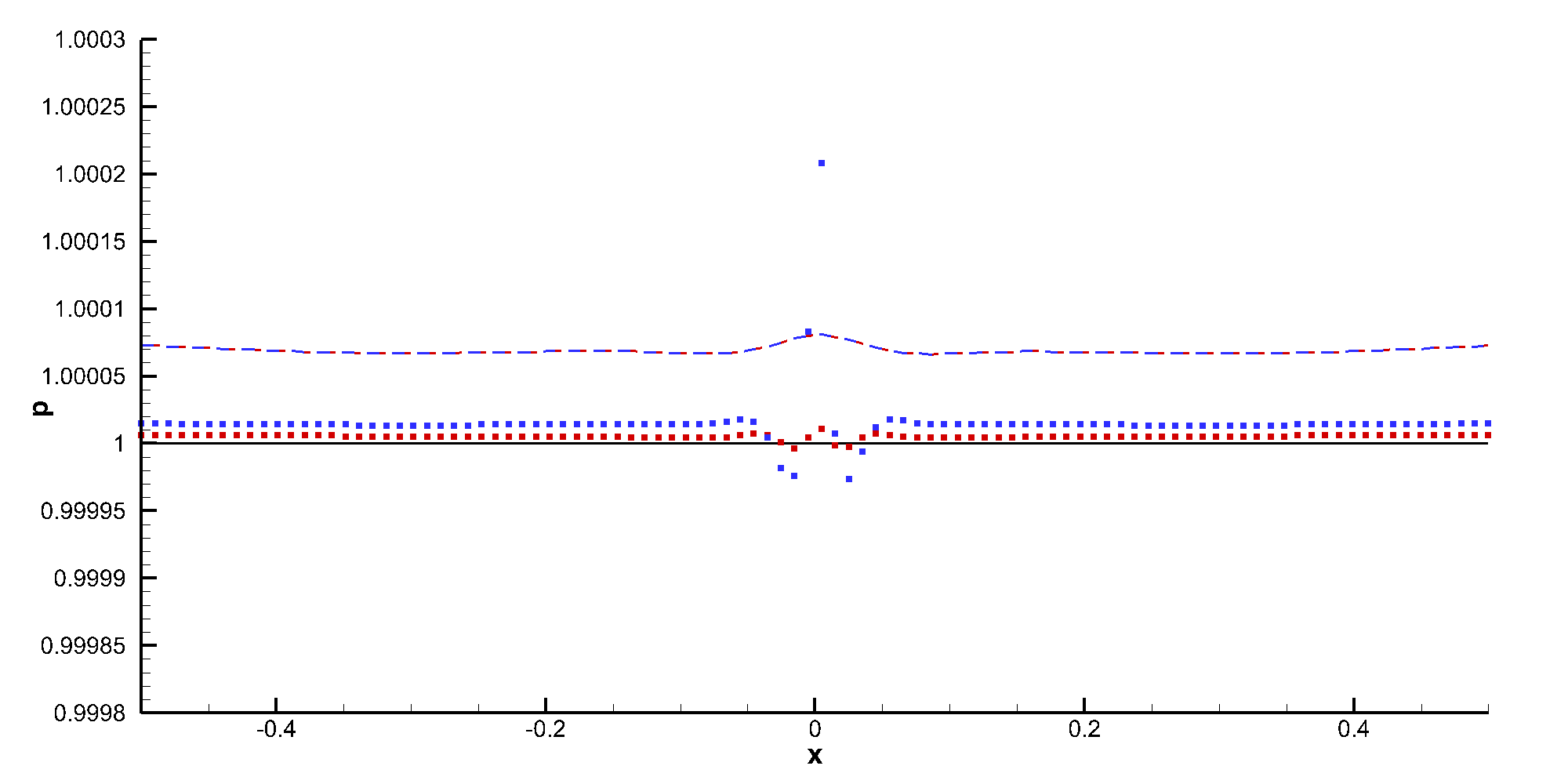}
		\vspace{-20pt}
		\caption{Shear wave not aligned with the \RIIcolor{tessellation}. 
			Here, we show the numerical results obtained with our first order schemes and our novel multidimensional Rieman solvers (squares) compared with those obtained with standard complete solvers (dashed lines). Our novel solvers significantly reduce the numerical dissipation.
		}
		\label{fig:shear}
	\end{figure}

	We then consider an isolated shear wave characterized by the following initial conditions 
	\begin{equation}
		(\rho, u, v, p)(\x) = 
		\begin{cases} 
			(1.0, 0, 0.1, 1.0 )    &  \ \text{ if } x < 0,  \\	
			(1.0, 0, -0.1, 1.0    ) &  \ \text{ if } x \ge 0, 
		\end{cases}    
	\end{equation}
	simulated on a square of dimension $[-0.5;0.5]\times[-0.5;0.5]$ with periodic boundary on top and bottom boundaries.
	We remark that, despite the mathematical discontinuity is aligned with $x=0$, our unstructured polygonal mesh, 
	is not constructed by imposing any constraints on the position of the initial generator points $c$, 
	thus it does not maintain this alignment property. 
	Hence, the results in Figure~\ref{fig:shear} refer to a shear wave not aligned with the mesh and demonstrate that our novel multidimensional solvers significantly reduce the numerical dissipation.
	
\begin{figure}[!b]
		\centering
		\includegraphics[width=0.5\linewidth]{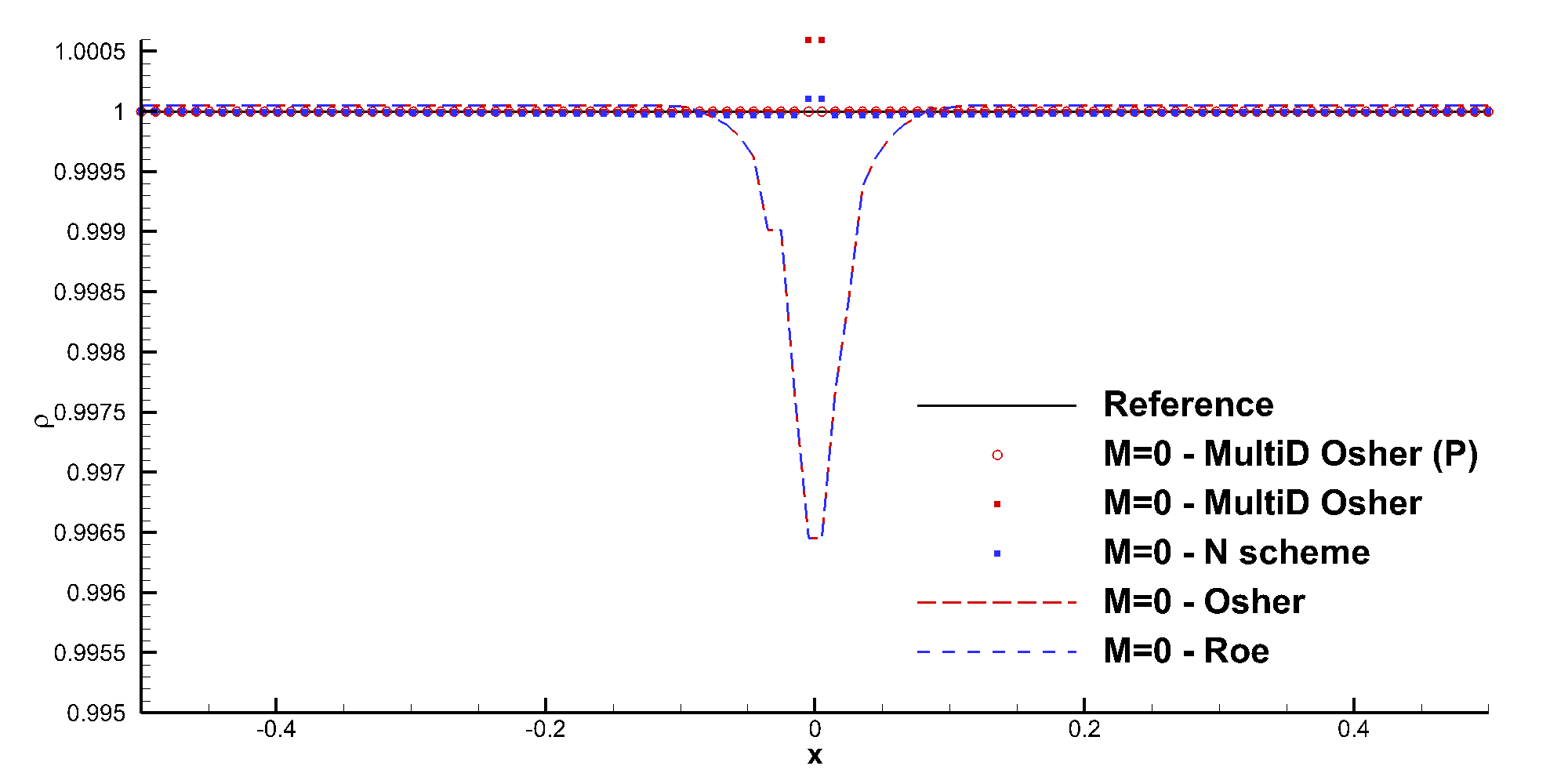}%
		\includegraphics[width=0.5\linewidth]{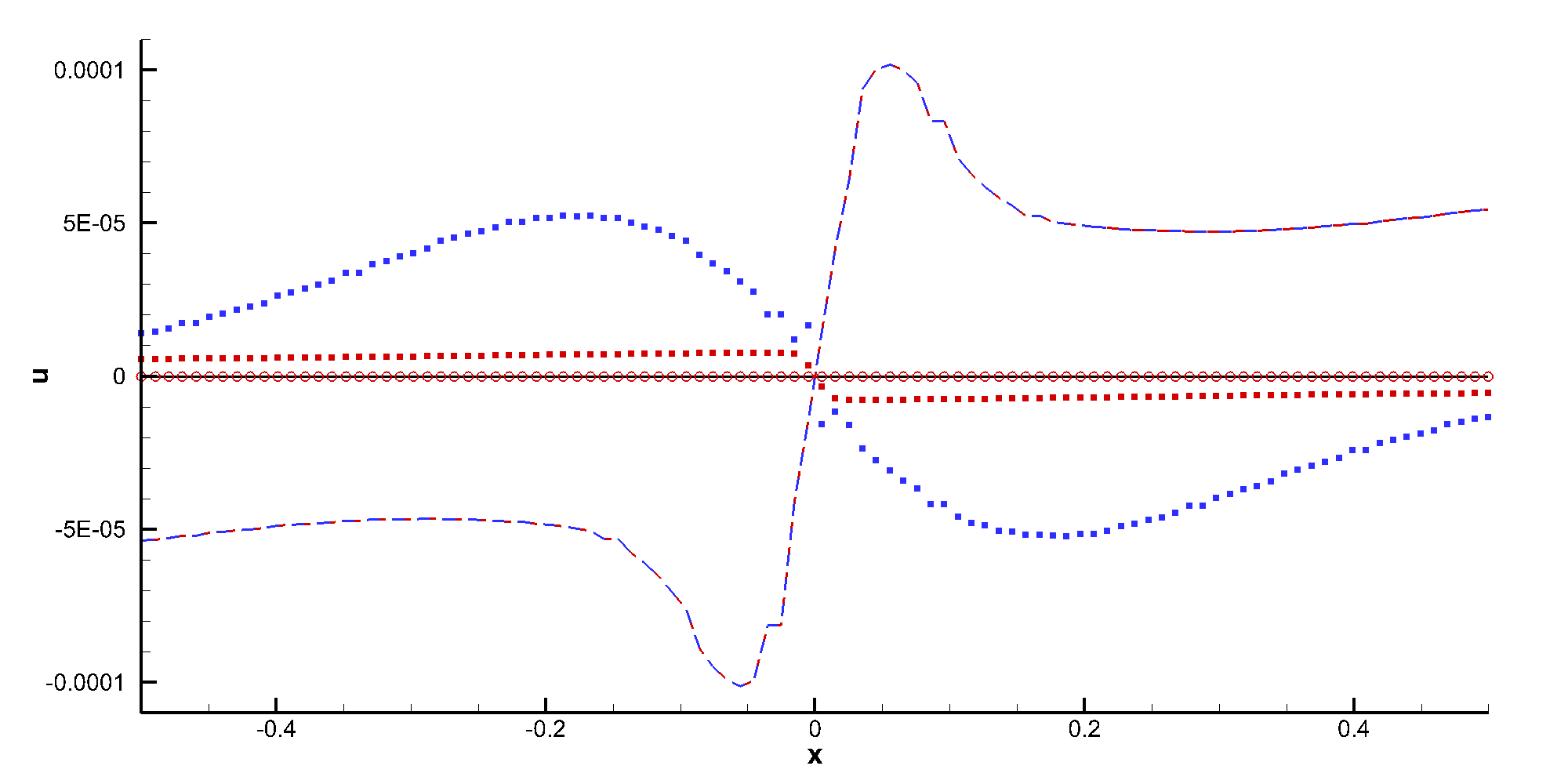}\\
		\includegraphics[width=0.5\linewidth]{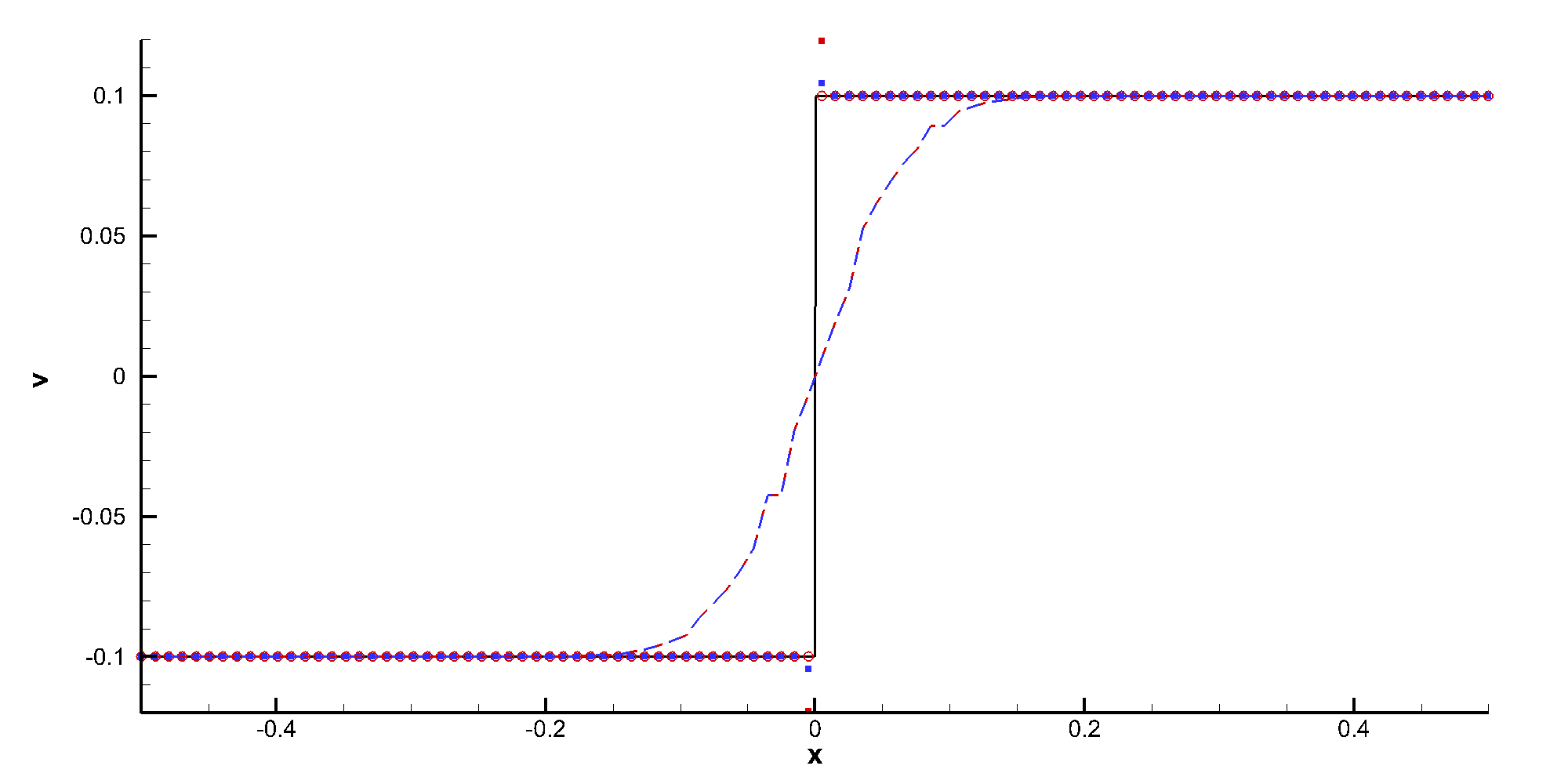}%
		\includegraphics[width=0.5\linewidth]{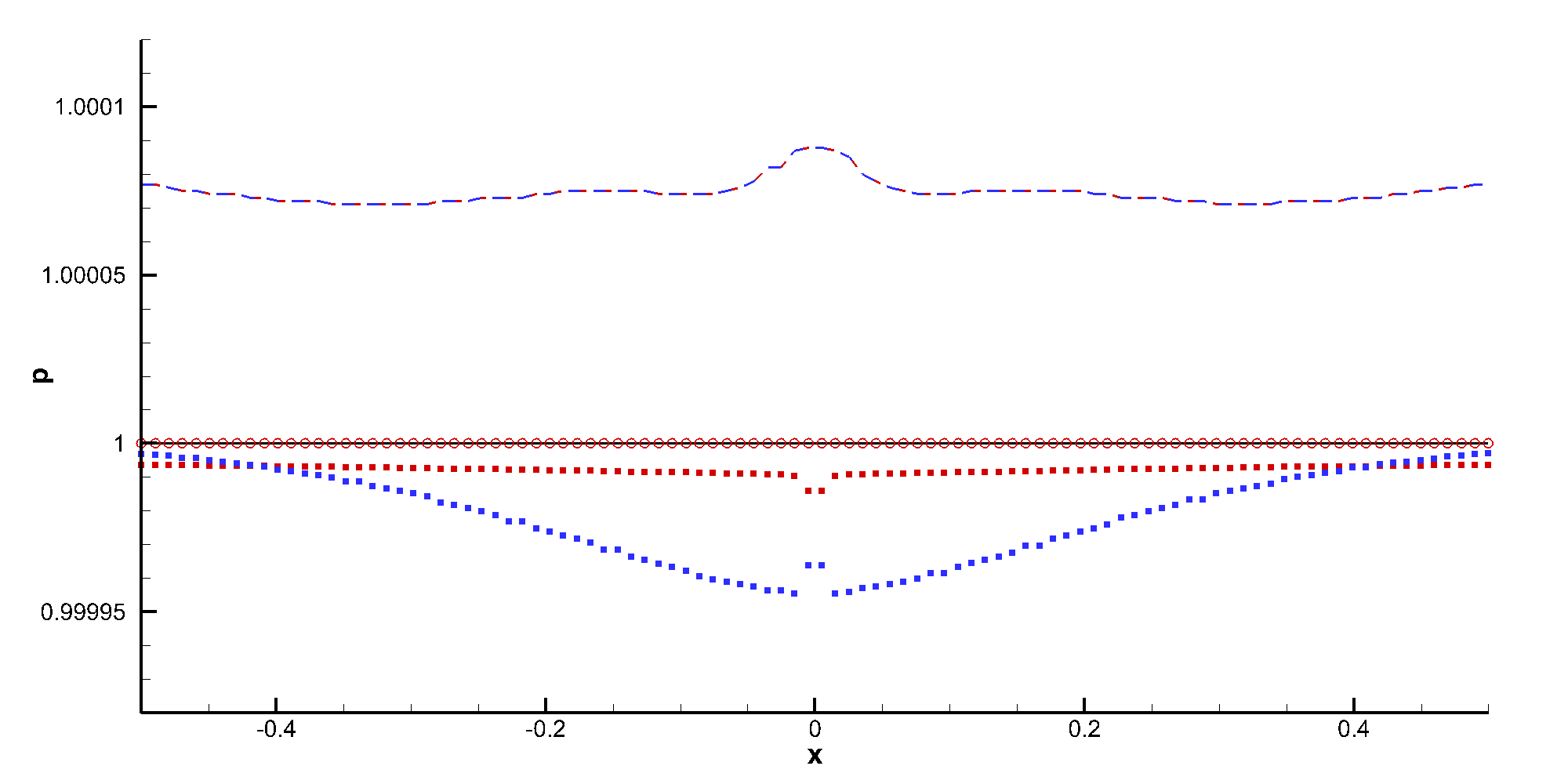}
		\vspace{-20pt}
		\caption{Shear wave simulated on a \RIIcolor{tessellation} where the generator points $c$ are aligned with the discontinuity. One can notice that our multidimensional Osher solver in primitive variables (red circles) exhibits zero numerical dissipation for the simulation of shear waves which are aligned with the mesh.}
		\label{fig:shear-aligned}
	\end{figure}

	In addition, we would like to underline that, if the generator points $c$ are chosen to be aligned with the discontinuity, 
	then our multidimensional Osher Riemann solvers,
	with the dissipation part of~\ref{eqn.multid.osher.final} computed using primitive variables (circles), 
	is able to simulate the isolated shear waves with machine precision, as can be noticed in Figure~\ref{fig:shear-aligned}. 
	
	
	\subsection{Circular Sod explosion problem}
	\label{ssec.sod}
	
	\begin{figure}[!b]
		\centering
		\includegraphics[width=0.5\linewidth,trim=1 1 1 1,clip]{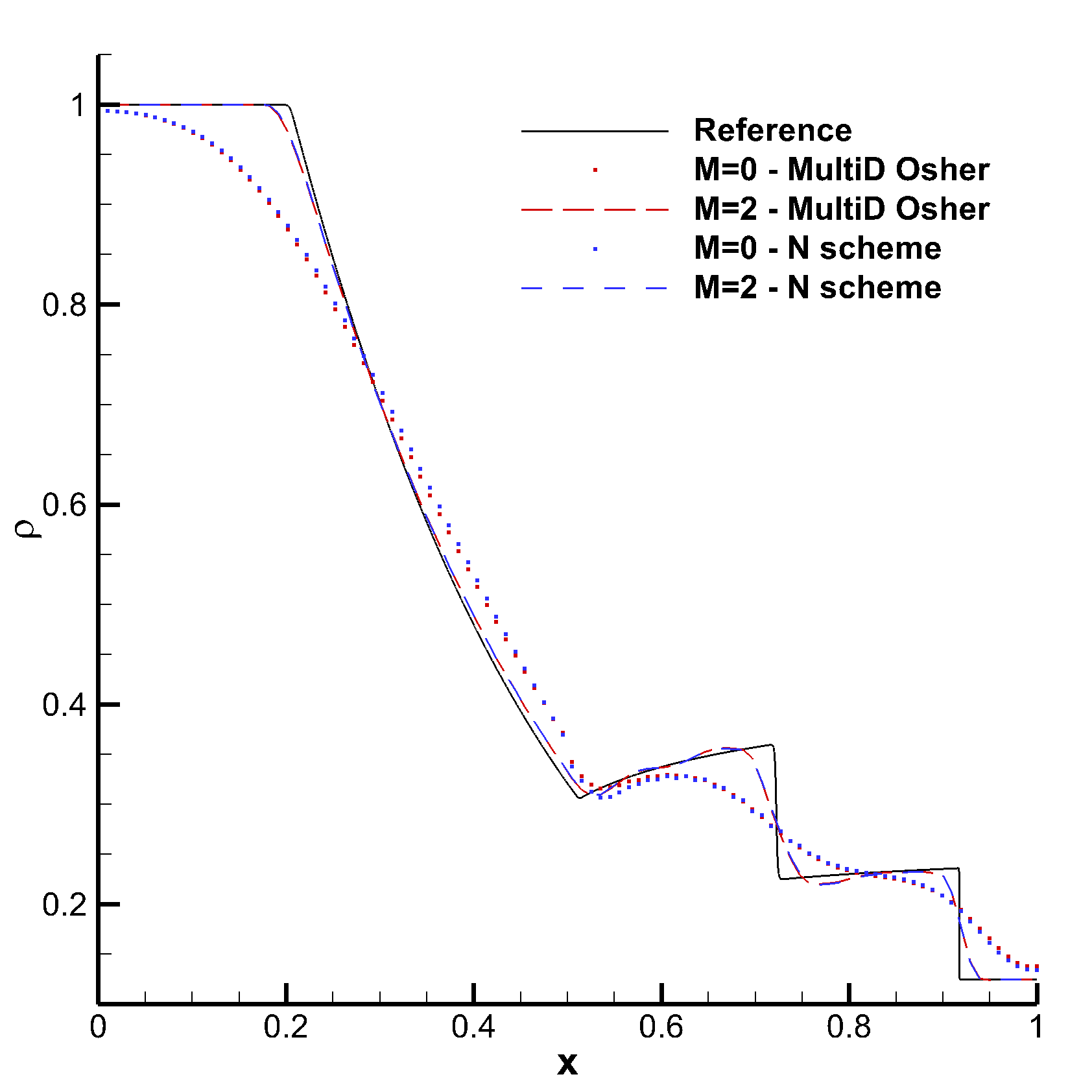}%
		\includegraphics[width=0.5\linewidth,trim=1 1 1 1,clip]{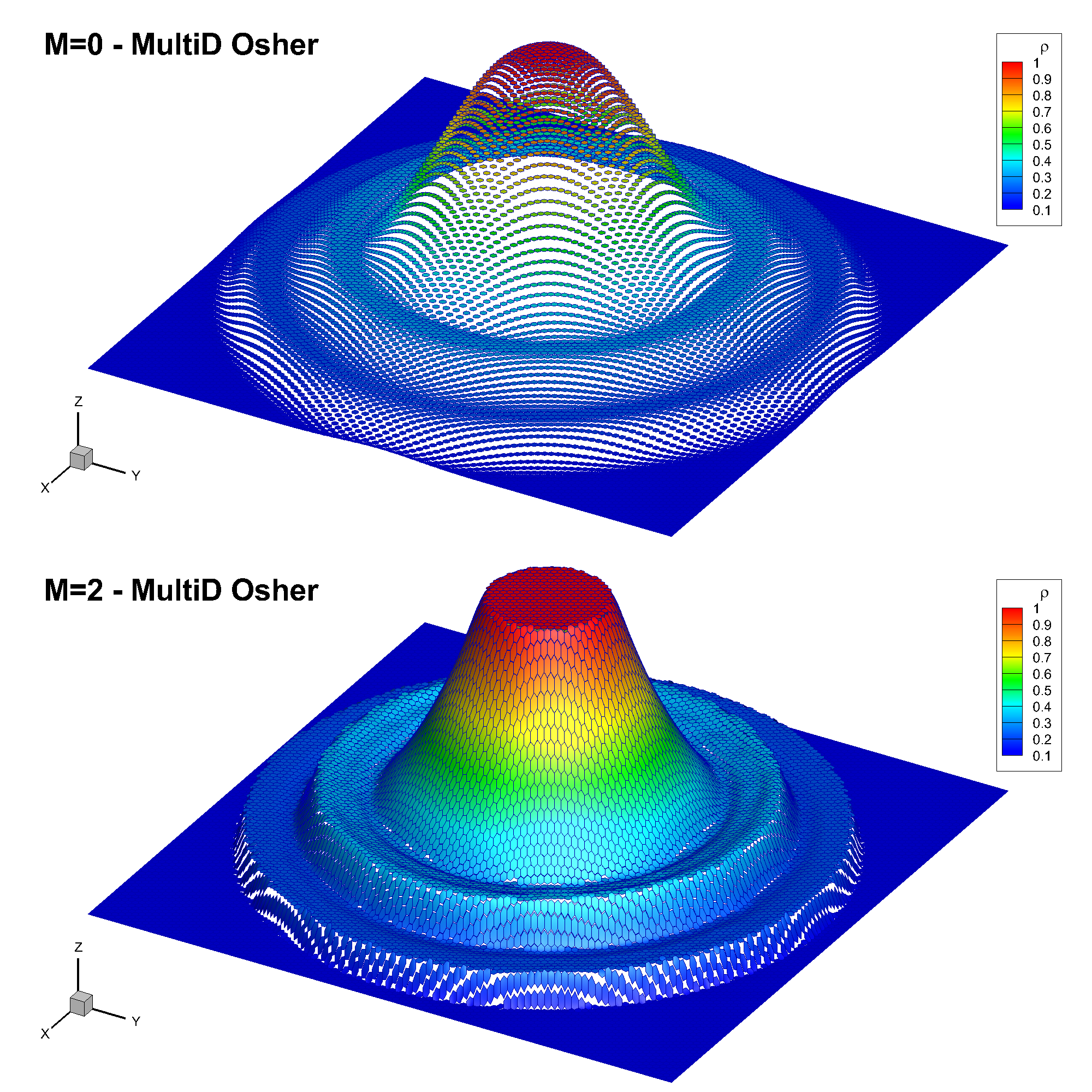}
		\caption{Circular Sod explosion. 
			We show the numerical results obtained with the first order (M=0) and the third order (M=2) CWENO-ADER schemes on the same mesh and our multidimensional Osher scheme. (The results obtained with the N scheme are really similar so we omit them).}
		\label{fig.sod}
	\end{figure}
	
	This circular explosion problem is a multidimensional extension of the classical Sod test case simulated on a square of dimension $[-1;1]\times[-1;1]$ with the following initial conditions composed of two  different states, separated by a discontinuity at radius $r_d=0.5$ 
	\begin{equation}
		(\rho, u, v, p)(\x) = 
		\begin{cases} 
			(1.0,  0, 0, 1 )        &  \ \text{ if } \  r \le r_d, \\	
			(0.125, 0, 0,  0.1    ) &  \ \text{ if } \  r > r_d. \\		
		\end{cases}
	\end{equation}
	
	Here, we numerically solve this problem with our novel multidimensional solvers and our first and third order finite volume schemes on a mesh with resolution $h \simeq 1/50$. 
	The final time is chosen to be $t_f=0.25$, so that the shock wave does not cross 
	the external boundary of the domain, where wall boundary conditions are imposed. 
	We have obtained a reference solution thanks to the rotational symmetry of the problem 
	which reduces to a one-dimensional problem with geometric source terms,
	which we have solved by using a classical second order TVD scheme on a very fine mesh. 
	
	The obtained results are reported in Figure~\ref{fig.sod} and show the increased resolution power of our high order approach. 
	
	\subsubsection{Lax shock tube}
	\label{ssec.lax}
	
\begin{figure}[!b]
		\centering
		\includegraphics[width=0.5\linewidth,trim=1 1 1 1,clip]{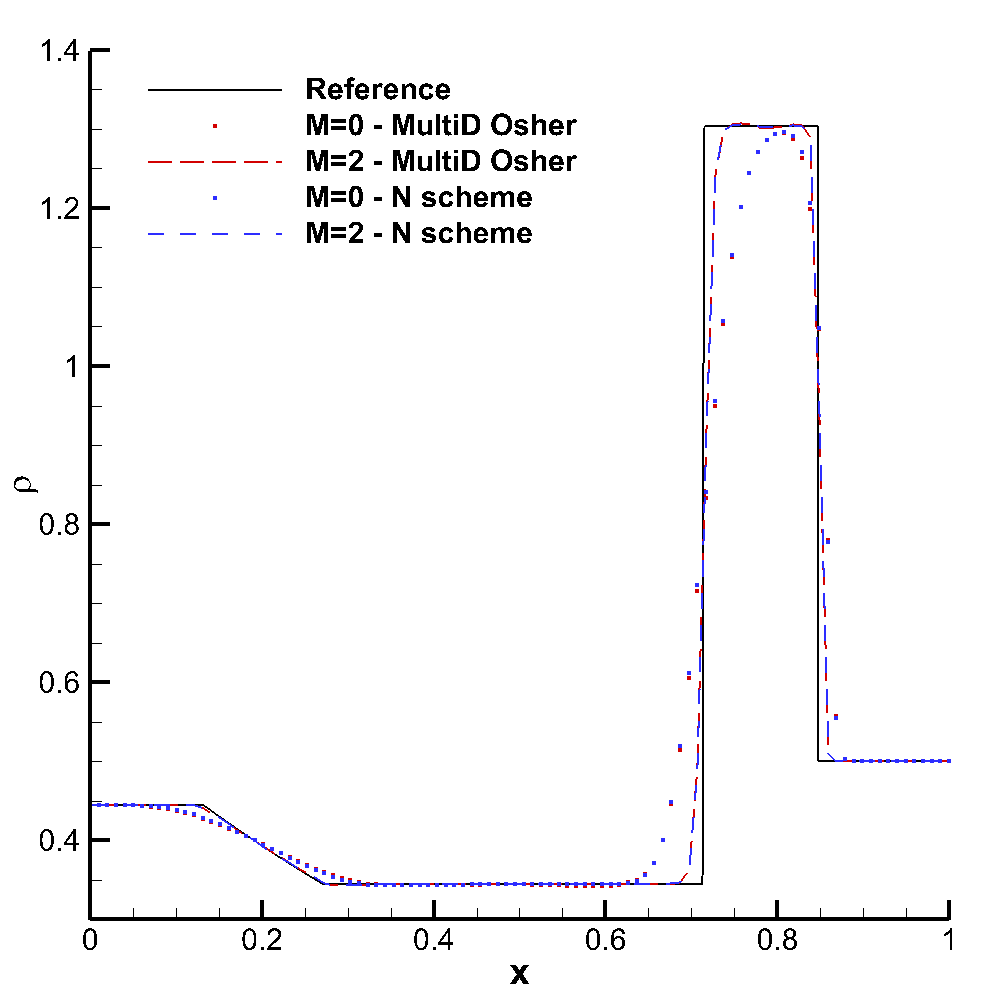}%
		\includegraphics[width=0.5\linewidth,trim=1 1 1 1,clip]{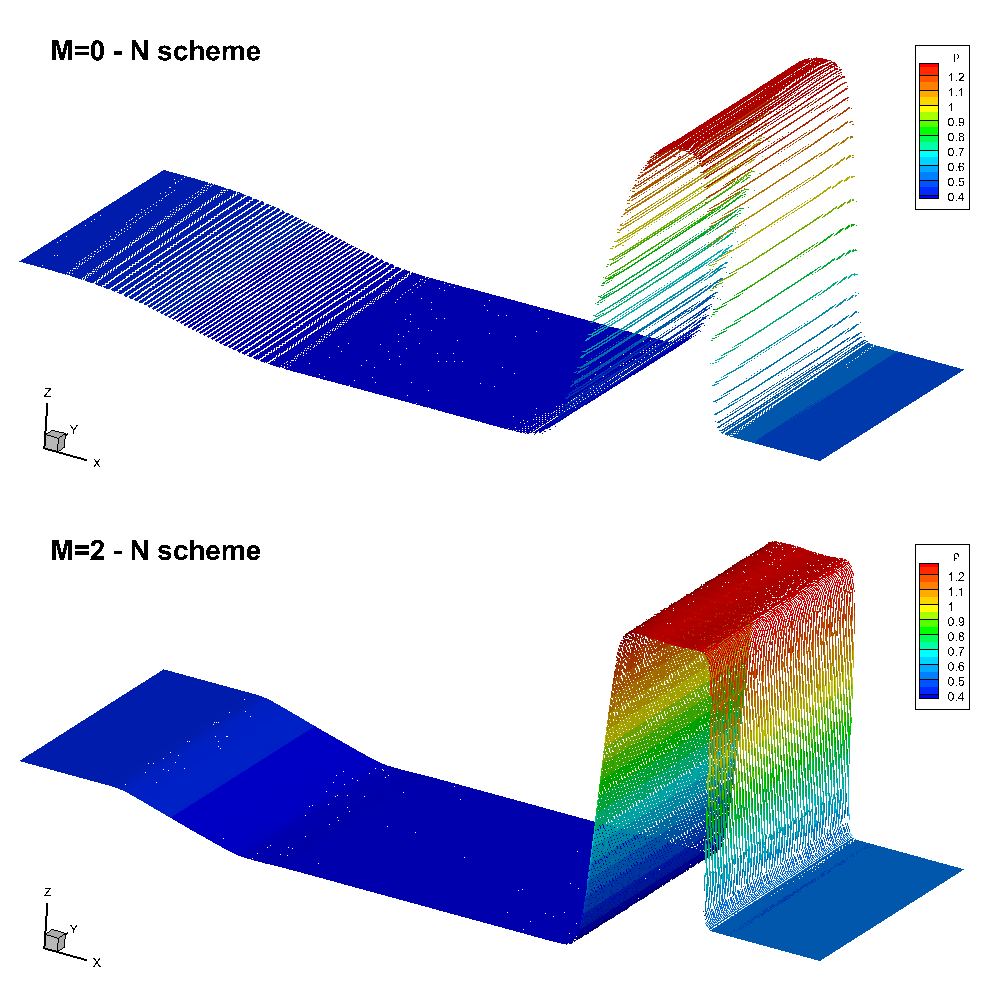}%
		\caption{Lax shock tube. 
			We show the numerical results obtained with the first order (M=0) and the third order (M=2) CWENO-ADER schemes on the same mesh and our N scheme Riemann solver. (The results obtained with the multidimensional Osher solver are really similar so we omit them).
		}
		\label{fig.lax}
	\end{figure}

	\RIIcolor{We continue with our set of Riemann problems by solving the Lax shock tube, originally presented in~\cite{lax2}, which, in addition to a shock wave, also exhibits a contact wave and a rarefaction fan.}
	The initial conditions for this problem are 
	\begin{equation}
		(\rho, u, v, p)(\x) = 
		\begin{cases} 
			(0.5  , 0, 0,  0.571    )   &  \ \text{ if } \  x > 0.5, \\
			(0.445,  0.698, 0, 3.528 )  &   \ \text{ if } \  x \le 0.5, \\	
		\end{cases}
	\end{equation}
	and it is solved on the domain $\Omega = [0,1]\times[0,0.1]$ discretized with a coarse mesh with resolution $h \simeq 1/50$. The final simulation time is $t_f =0.14$.
	We report the numerical results obtained with our multidimensional Riemann solvers and our first and third order finite volume schemes in Figure~\ref{fig.lax}.
	
	\subsection{Two dimensional Riemann Problem}
	\label{ssec.2drp}
	
	In this Section, we consider some of the truly multidimensional Riemann problems presented and studied in~\cite{laxliu98}.
	
	\paragraph{Configuration 12}
	
	\begin{figure}[!b]
		\centering
		\includegraphics[width=0.333\linewidth,trim=1 1 1 1,clip]{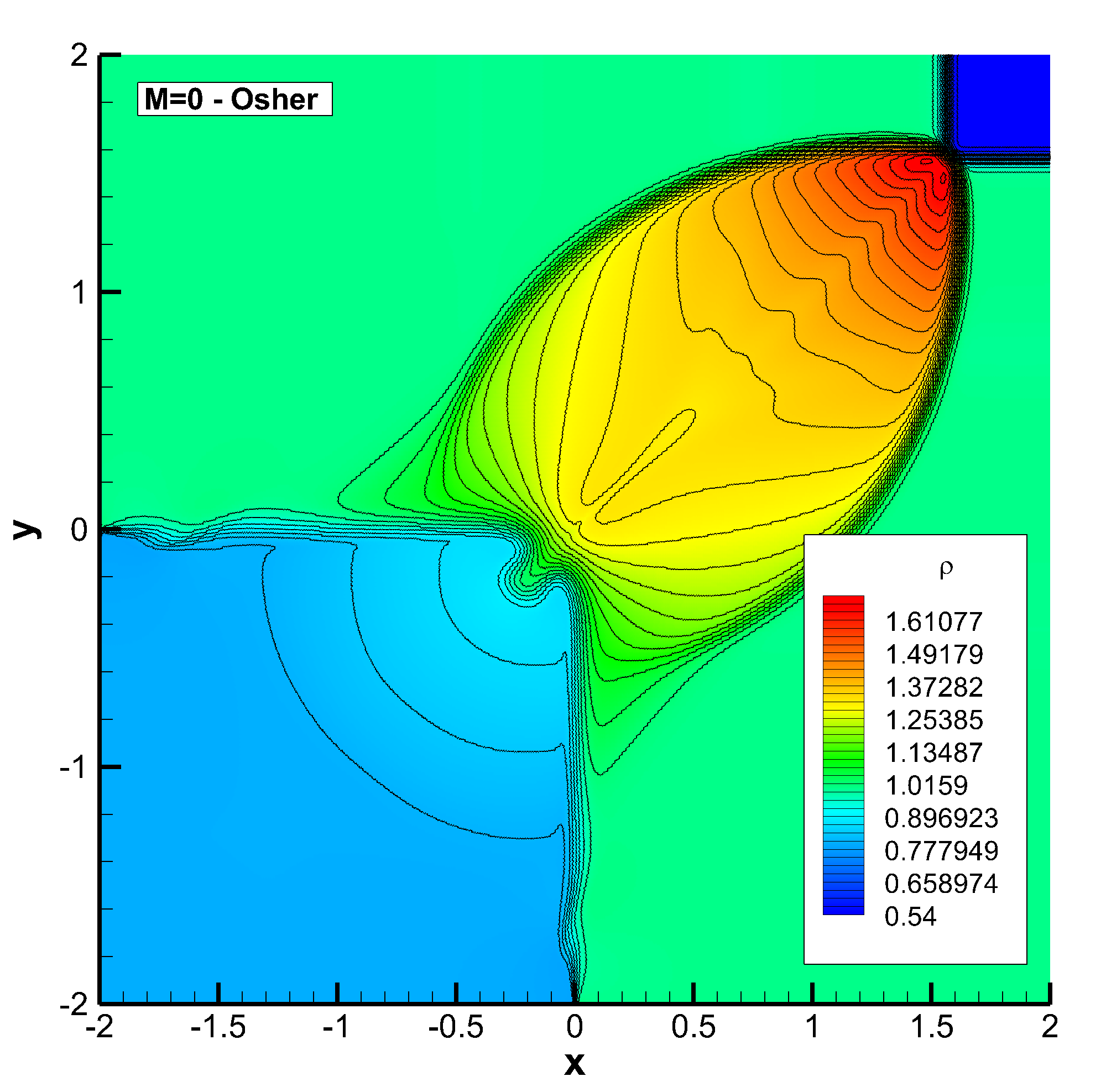}%
		\includegraphics[width=0.333\linewidth,trim=1 1 1 1,clip]{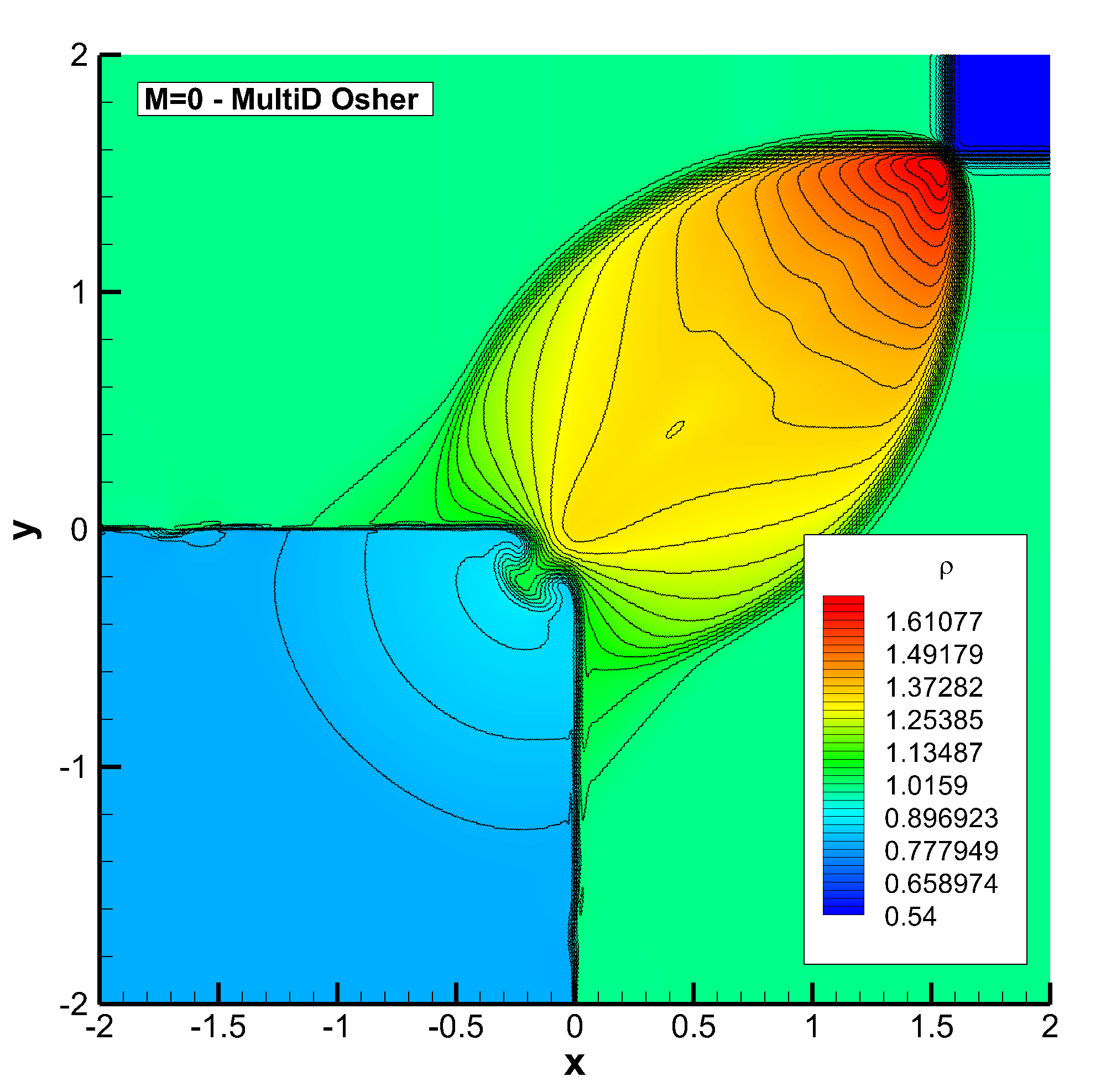}%
		\includegraphics[width=0.333\linewidth,trim=1 1 1 1,clip]{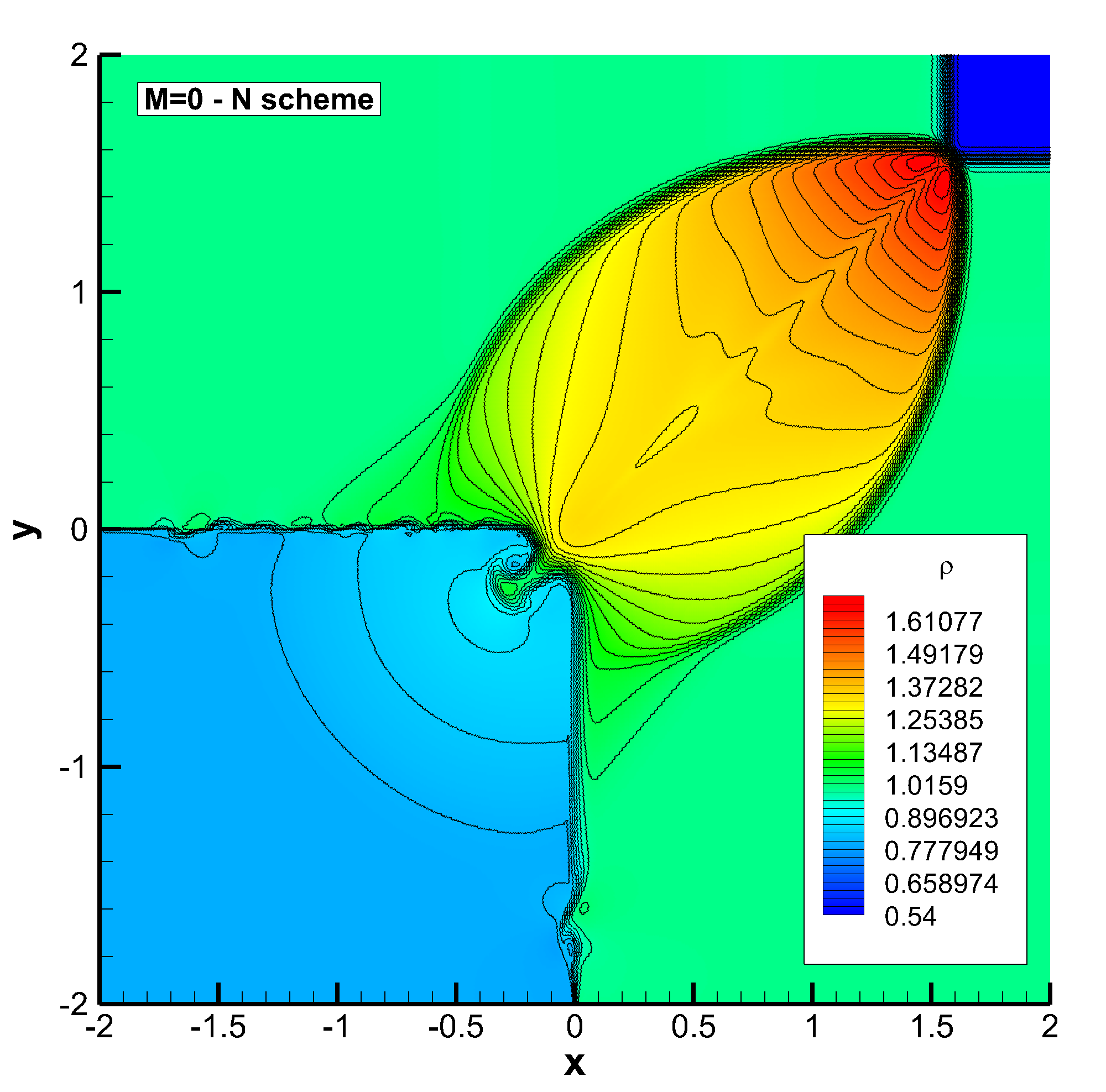}\\    
		\includegraphics[width=0.333\linewidth,trim=1 1 1 1,clip]{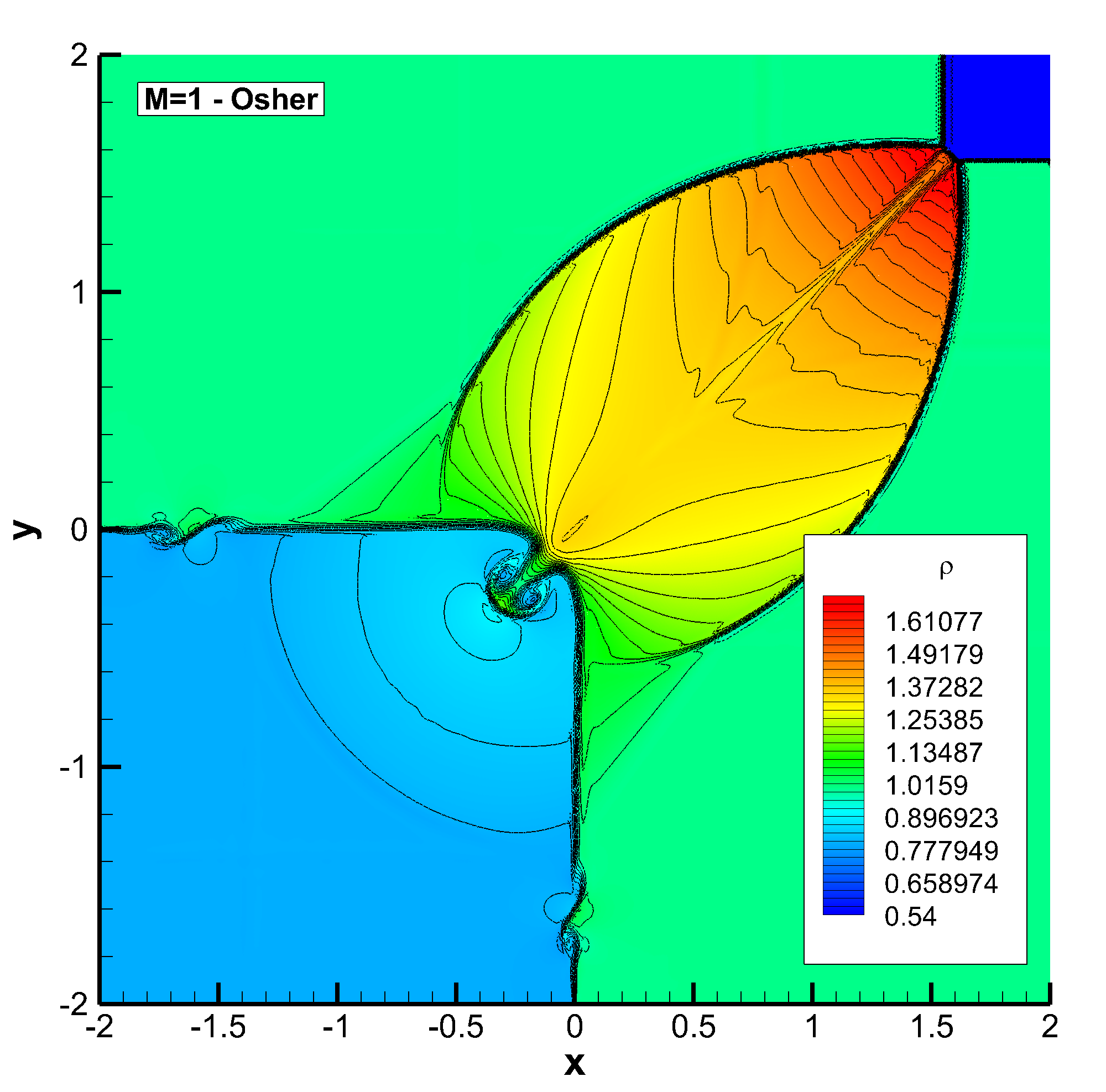}%
		\includegraphics[width=0.333\linewidth,trim=1 1 1 1,clip]{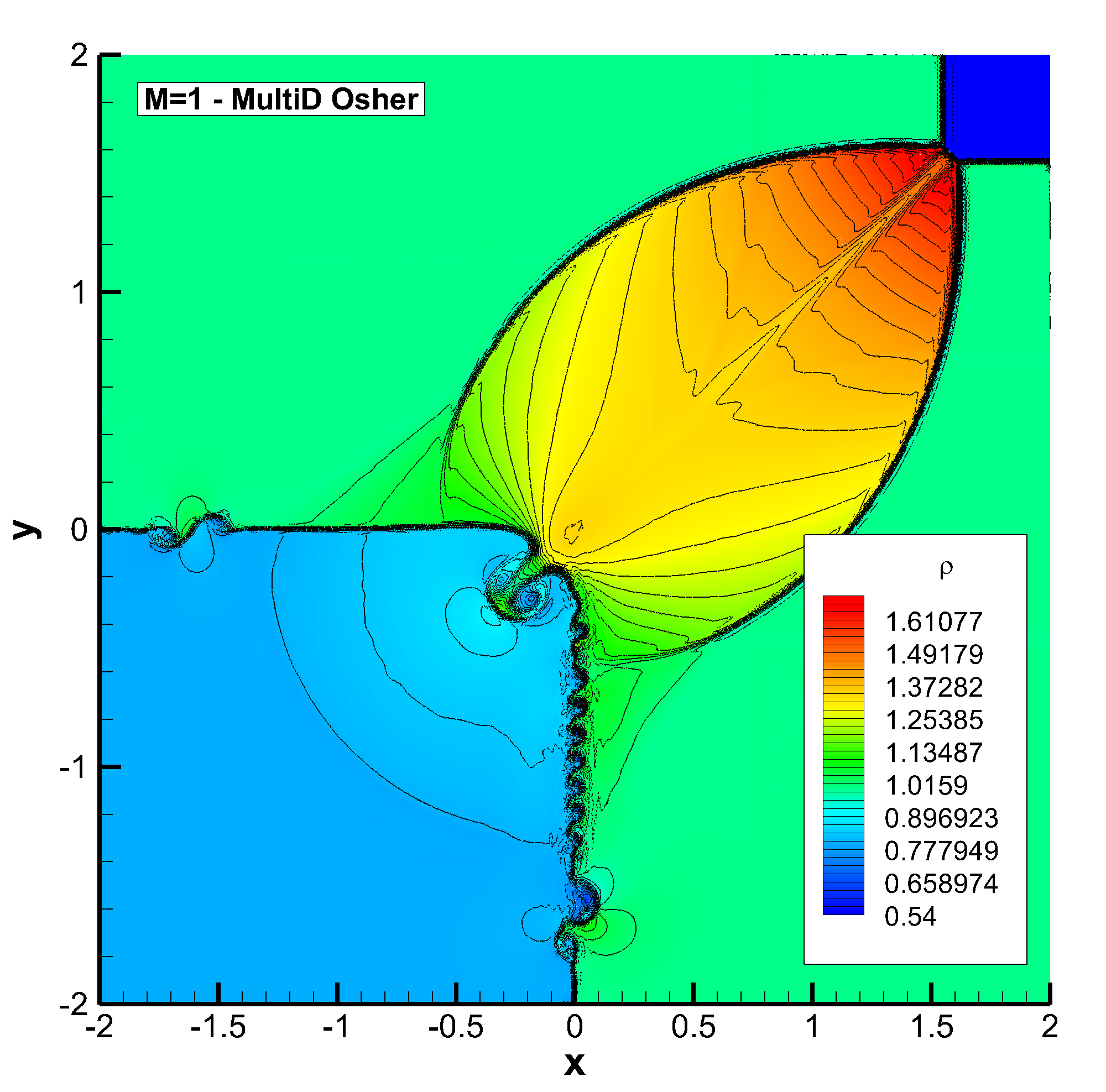}%
		\includegraphics[width=0.333\linewidth,trim=1 1 1 1,clip]{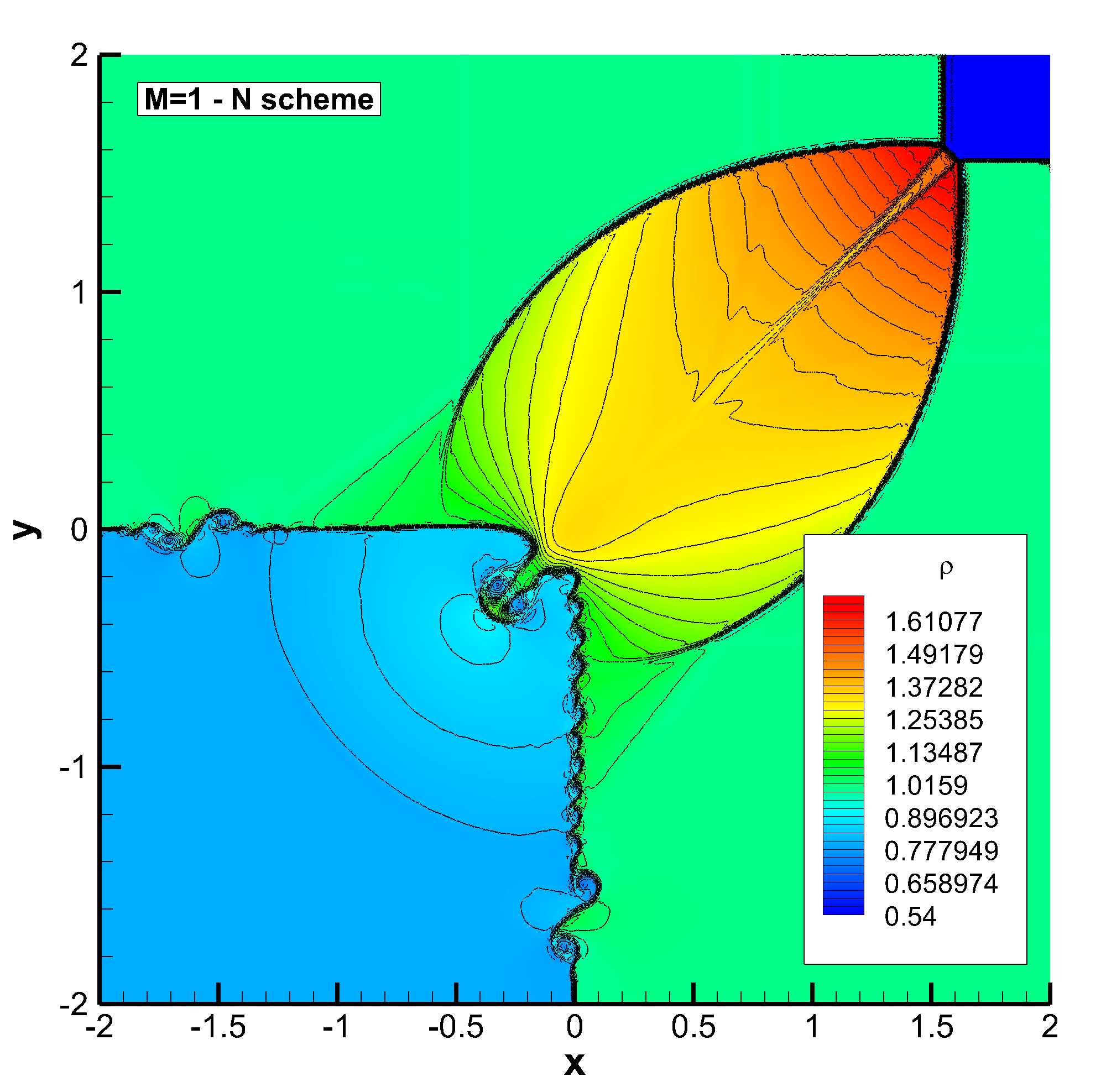}\\[10pt]
		\includegraphics[width=0.333\linewidth,trim=1 1 1 1,clip]{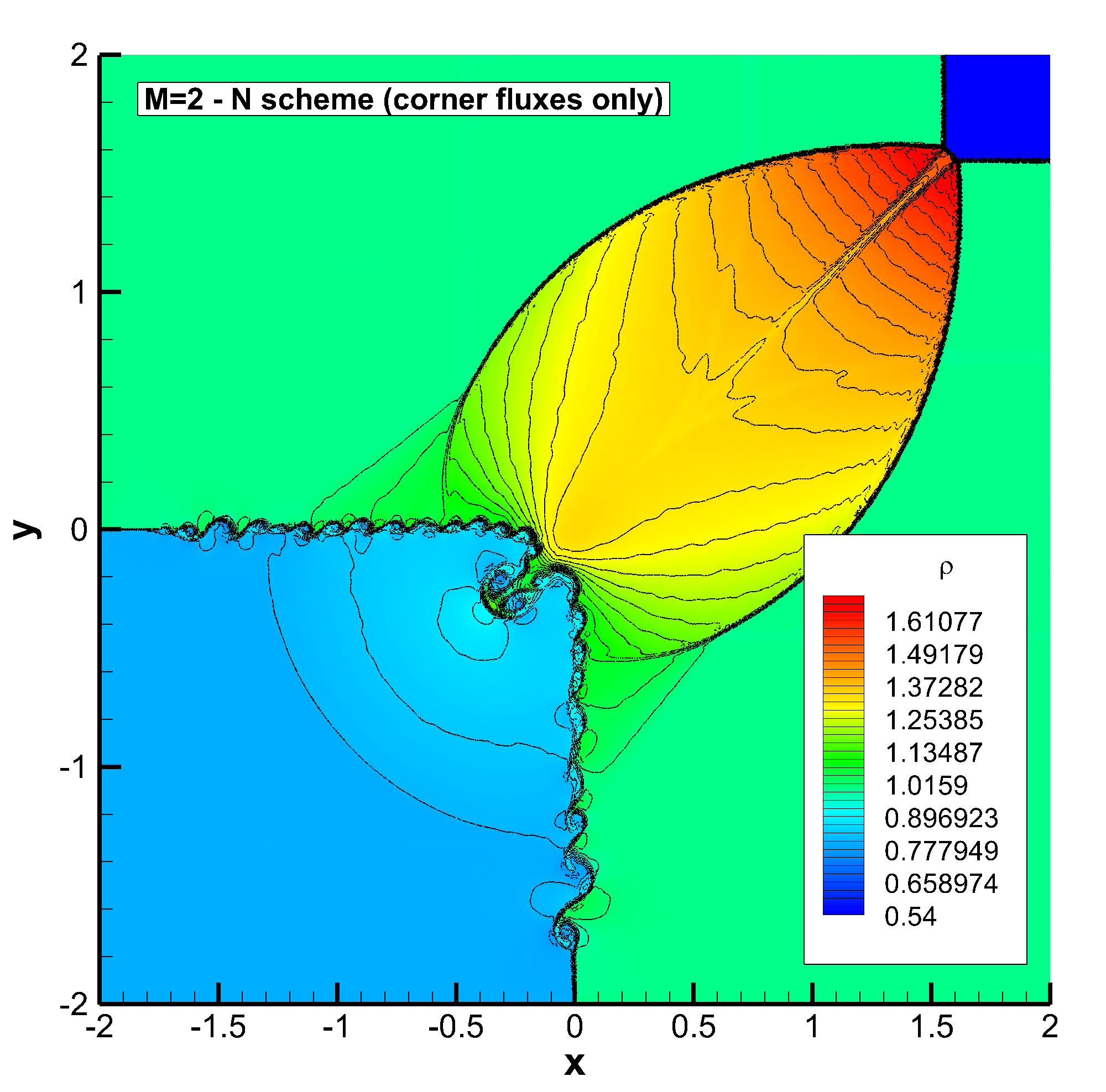}%
		\includegraphics[width=0.333\linewidth,trim=1 1 1 1,clip]{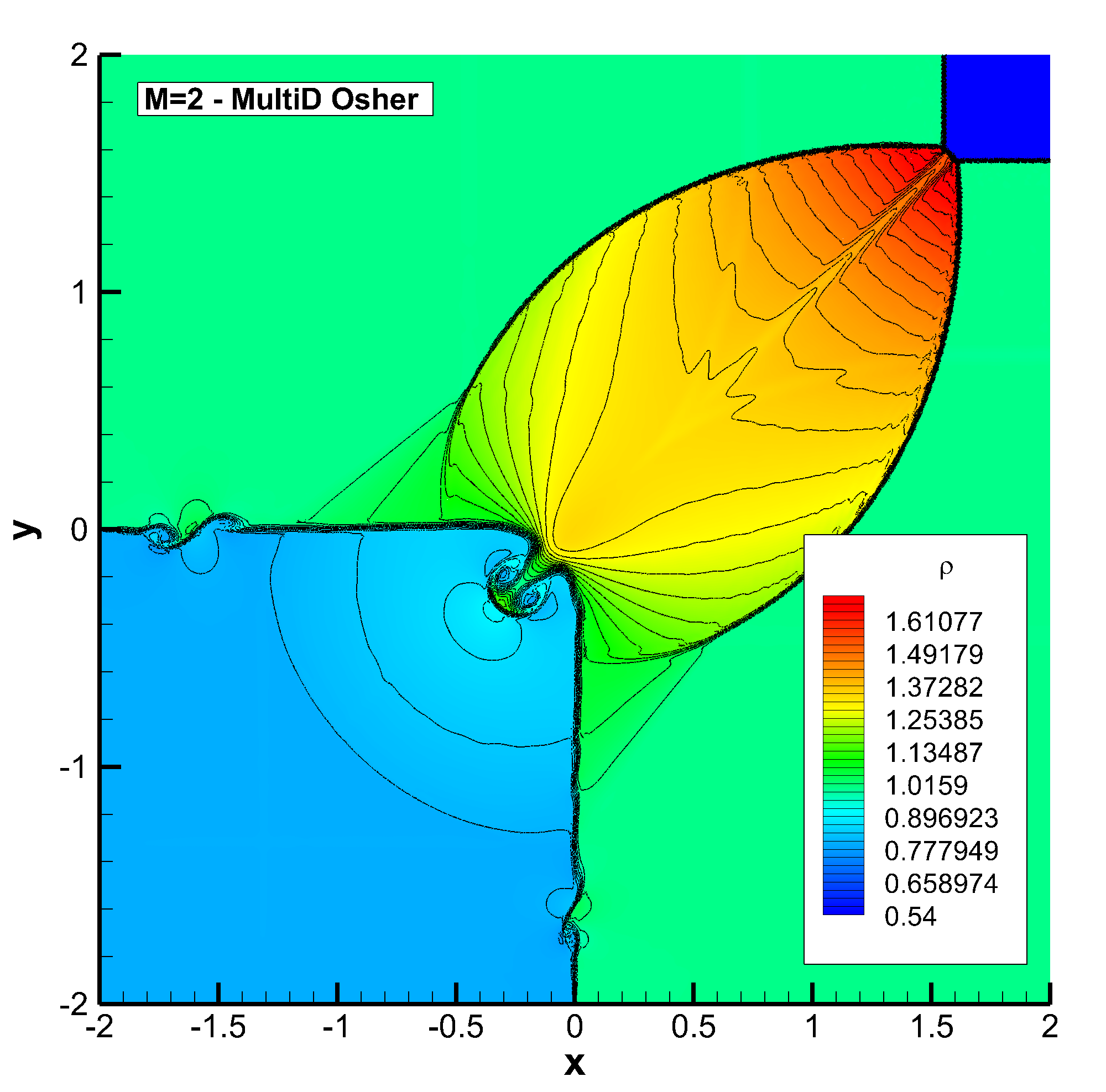}%
		\includegraphics[width=0.333\linewidth,trim=1 1 1 1,clip]{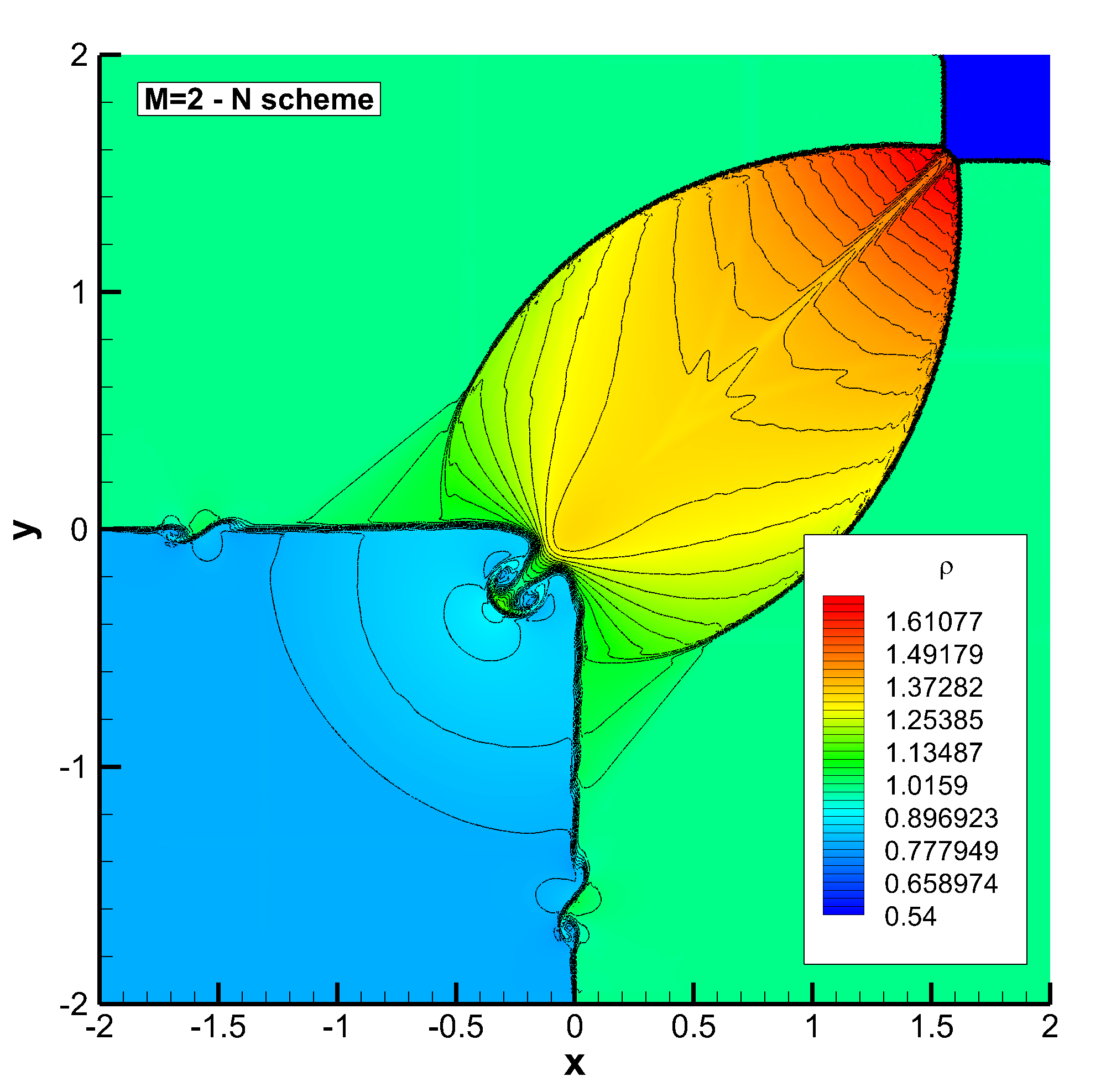}%
		\caption{Two dimensional Riemann Problem (configuration 12 of~\cite{laxliu98}). 
			Here we report the numerical results obtained with our first order \RIIcolor{(M=0)}, second order \RIIcolor{(M=1)}, and third order \RIIcolor{(M=2)} schemes \RIIcolor{over a mesh with $h\simeq1/110$.}
			In particular, we show the density profile with 40 contours line in the interval $[0.54, 1.7]$.
			We refer to the test case description for detailed comments on the obtained results. 
		}
		\label{fig.2drp_12}
	\end{figure}
	\begin{figure}[!b]
		\centering
		\includegraphics[width=0.5\linewidth,trim=1 1 1 1,clip]{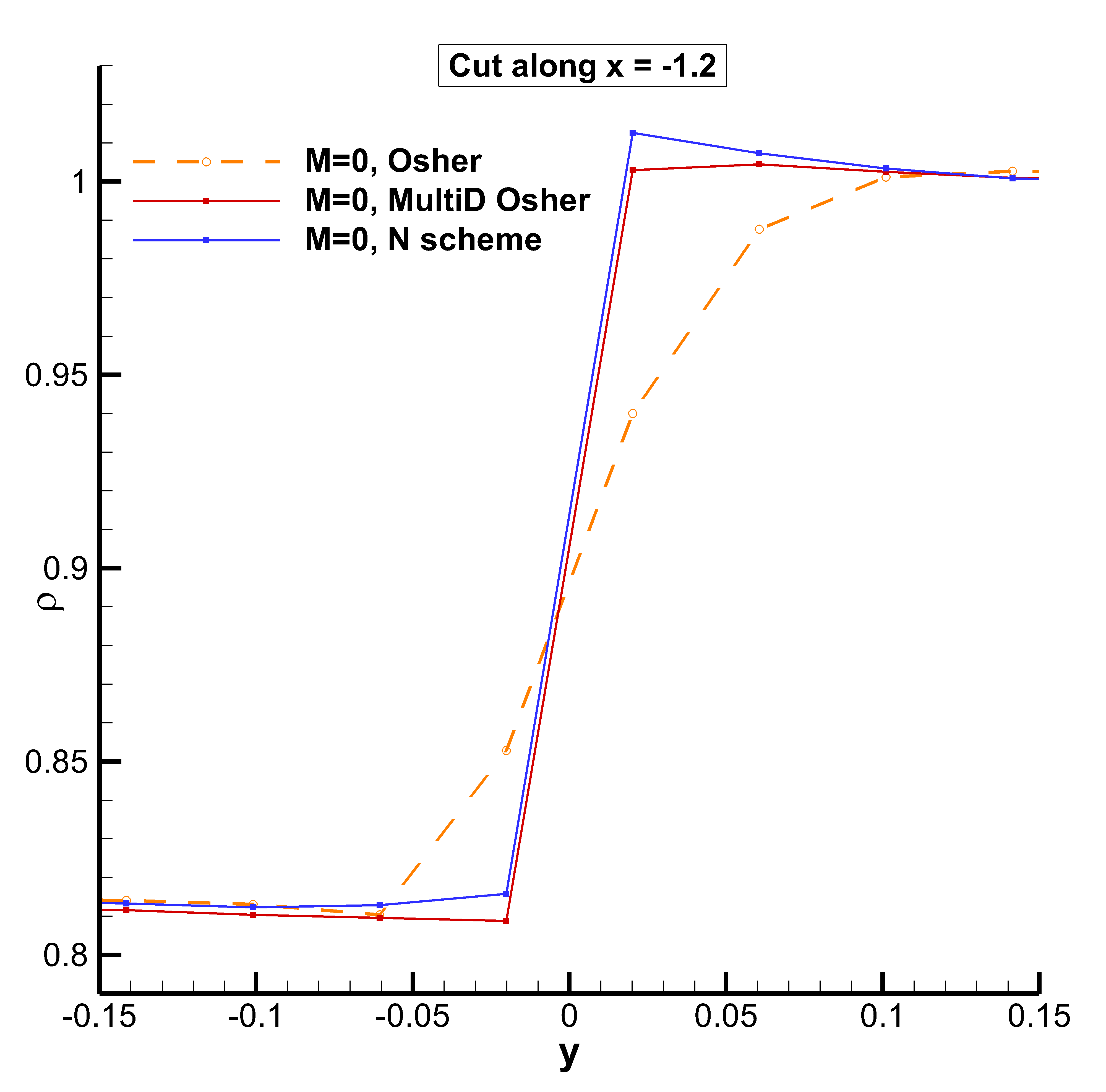}%
		\includegraphics[width=0.5\linewidth,trim=1 1 1 1,clip]{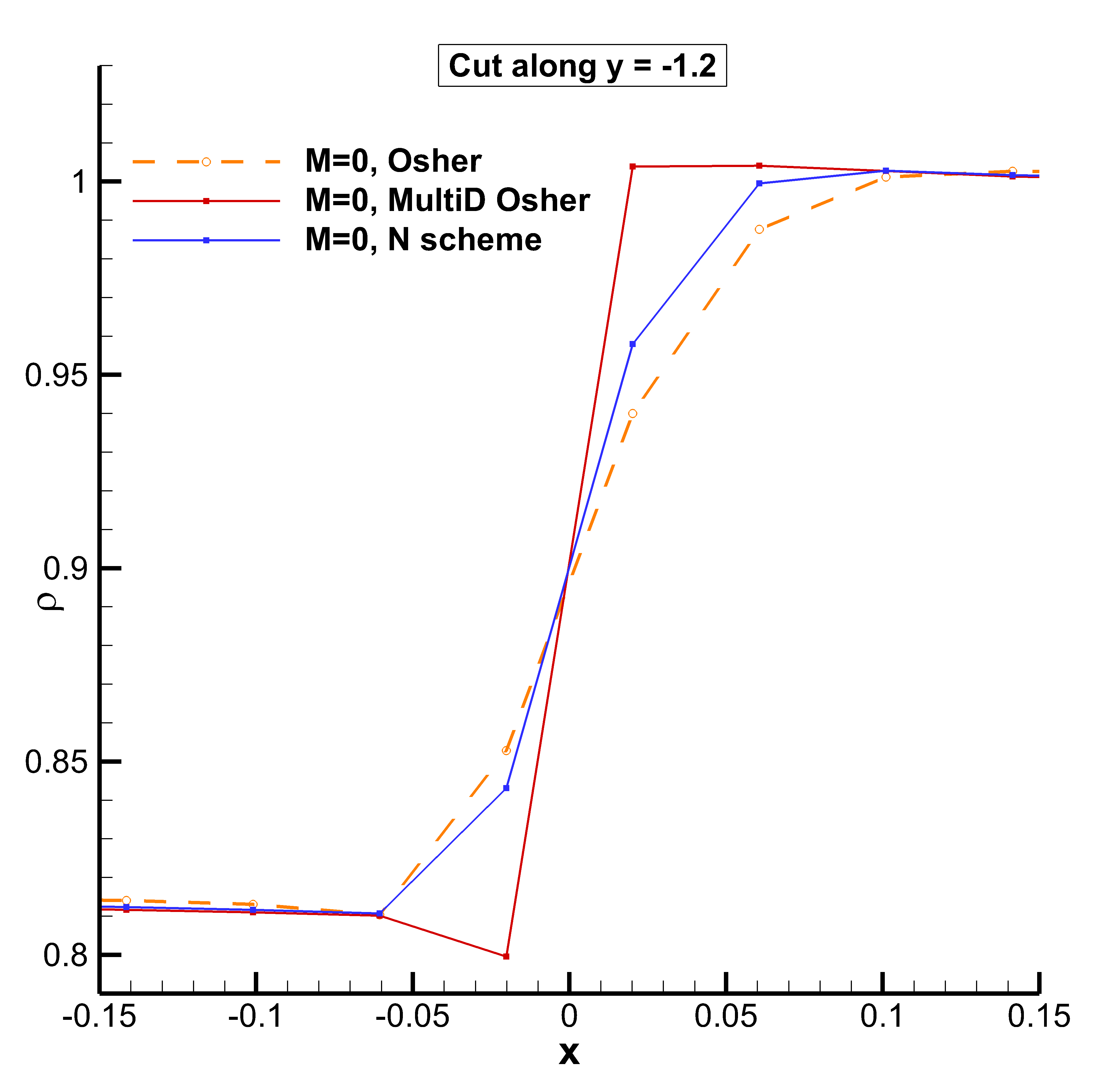}%
		\caption{Cut along $x=-1.2$ (left) and $y=-1.2$ (right) for the density values of the solution of the two dimensional Riemann Problem (configuration 12 of~\cite{laxliu98}). One can easily notice that the discontinuity between the two steady states is captured in a sharper manner by our multidimensional solvers.}
		\label{fig.2drp_12_cut}
	\end{figure}

	We start with the configuration number 12 which involves both steady contact discontinuities, and 
	strong moving shocks whose interaction produces two triple points from which
	parallel contact lines emerge. From the initial singularity a strong fluid jet emerges whose resolution is enhanced by high order and less dissipative methods as our multidimensional solvers. 
	
	Here, the computational domain is the square  $\Omega = [-2,2]\times[-2,2]$ discretized with a mesh with resolution $h \simeq 1/110$, and the discontinuous initial conditions read as follows
	\begin{equation}
		(\rho, u, v, p)(\x) = 
		\begin{cases} 
			(1.0,    0.7276, 0.0,    1.0 )    &  \ \text{ if } x < 0 \ \&  \ y > 0,  \\	
			(0.8,    0.0,    0.0,    1.0    ) &  \ \text{ if } x < 0 \ \&  \ y < 0,\\
			(0.5313, 0.0,    0.0,    0.4 )    &  \ \text{ if } x > 0 \ \& \ y > 0,  \\	        
			(1.0,    0.0,    0.7276, 1.0    ) &  \ \text{ if } x > 0 \ \& \ y < 0.   \\		
		\end{cases}    
	\end{equation}
	In Figure~\ref{fig.2drp_12}, we report the numerical results obtained at time $t_f=0.25$ with a standard 1d Osher solver compared with those achieved with our multidimensional solvers of first, second, and third order of accuracy. Moreover, in Figure~\ref{fig.2drp_12_cut} we show the cut of the density profile along $x=-1.2$ and $y=-1.2$.
	We can clearly \RIIcolor{observe} that both our multidimensional solvers exhibit increased resolution compared to the 1d solver, showing many more \RIIcolor{vortices} along the traveling shock discontinuities.
	We would also like to remark that the third order multidimensional schemes (in the last two panel of the Figure), as presented in Section~\ref{sec.highorder}, which also involve an edge contribution computed with a 1d solver, have less resolution power w.r.t. the second order schemes, which only use corner fluxes. 
	Moreover, a variant of the third order scheme (left-bottom panel), employing a third order space-time polynomial  reconstruction but a second order order flux computation, based only on corner fluxes, 
	effectively provides numerical results more accurate w.r.t the second order schemes.
	However, due to the lack of sufficient accuracy of the numerical quadrature, this variant is only second order accurate. This feature will be object of future investigations.

	\paragraph{Configuration 3}

	We now consider the configuration number 3 of the same reference. The initial solution involves again four states corresponding to    two oblique shocks and two normal shocks converging in one point.
	The interaction leads to the appearance of a strong normal shock connecting the two oblique shocks, and of two lambda shocks stemming from the normal ones.
	The lambda shocks lead to the appearance of two strong moving contact discontinuities, and of a strong jet of fluid 
	whose resolution is enhanced by high order and less dissipative methods as our multidimensional solvers.
	
	Here, the computational domain is the square  $\Omega = [0,1.2]\times[0,1.2]$ discretized with meshes of resolution $h \simeq 1/200$ or $h \simeq 1/400$, and the discontinuous initial conditions read as follows
	\begin{equation}
		(\rho, u, v, p)(\x) = 
		\begin{cases} 
			(0.5323, 1.206,  0.0,    0.3  )    &  \ \text{ if } x < 1 \ \&  \ y > 1,  \\	
			(0.138,  1.206,  1.206, 0.029 ) &  \ \text{ if } x < 1 \ \&  \ y < 1,\\
			(1.5,    0.0,    0.0,    1.5 )    &  \ \text{ if } x > 1 \ \& \ y > 1,  \\	        
			(0.5323, 0.0,    1.206, 0.3    ) &  \ \text{ if } x > 1 \ \& \ y < 1.   \\		
		\end{cases}    
	\end{equation}
	The final computational time is $t_f=1.0$.
	
\begin{figure}[!b]
		\centering
		\includegraphics[width=0.333\linewidth,trim=1 1 1 1,clip]{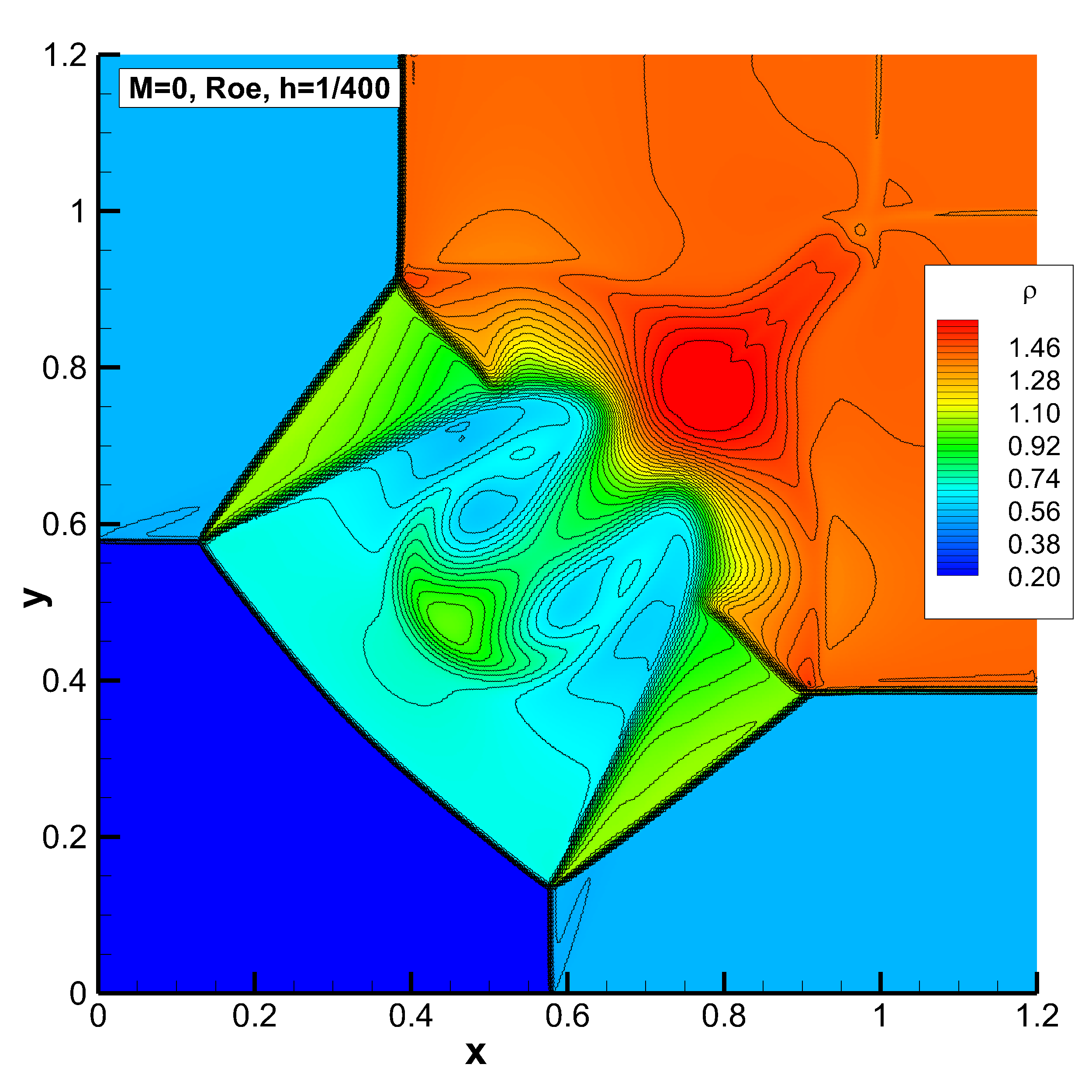}%
		\includegraphics[width=0.333\linewidth,trim=1 1 1 1,clip]{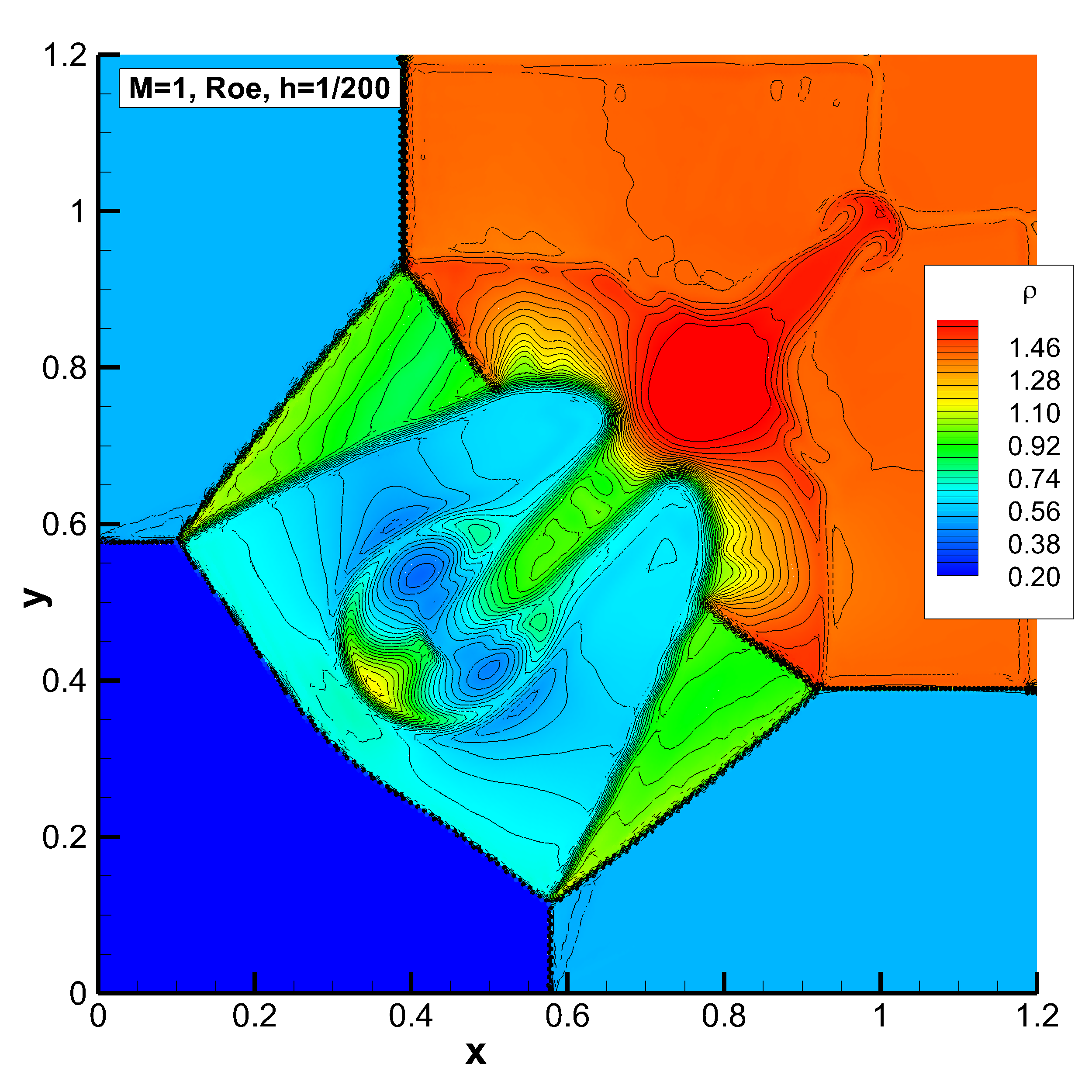}%
		\includegraphics[width=0.333\linewidth,trim=1 1 1 1,clip]{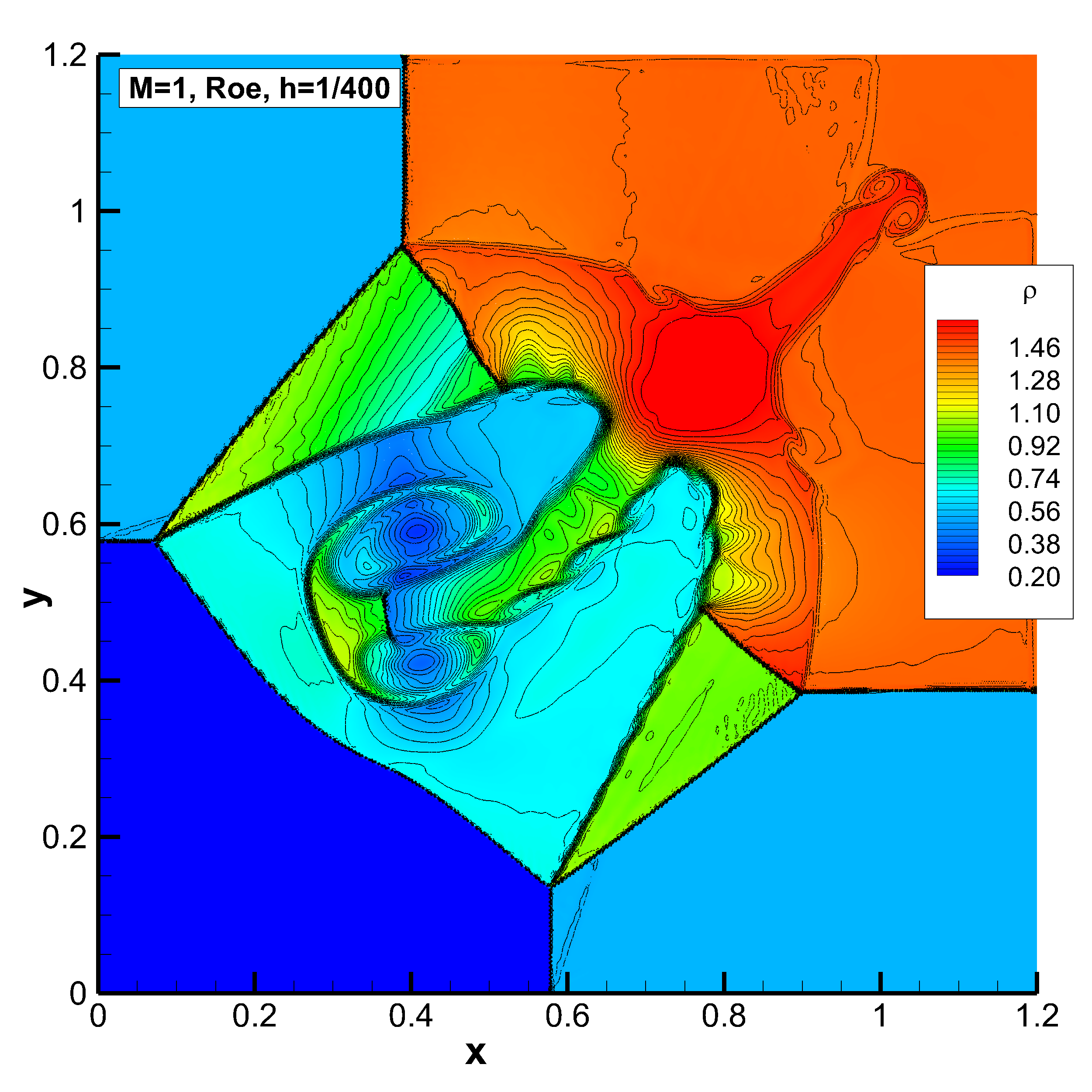}\\    
		\includegraphics[width=0.333\linewidth,trim=1 1 1 1,clip]{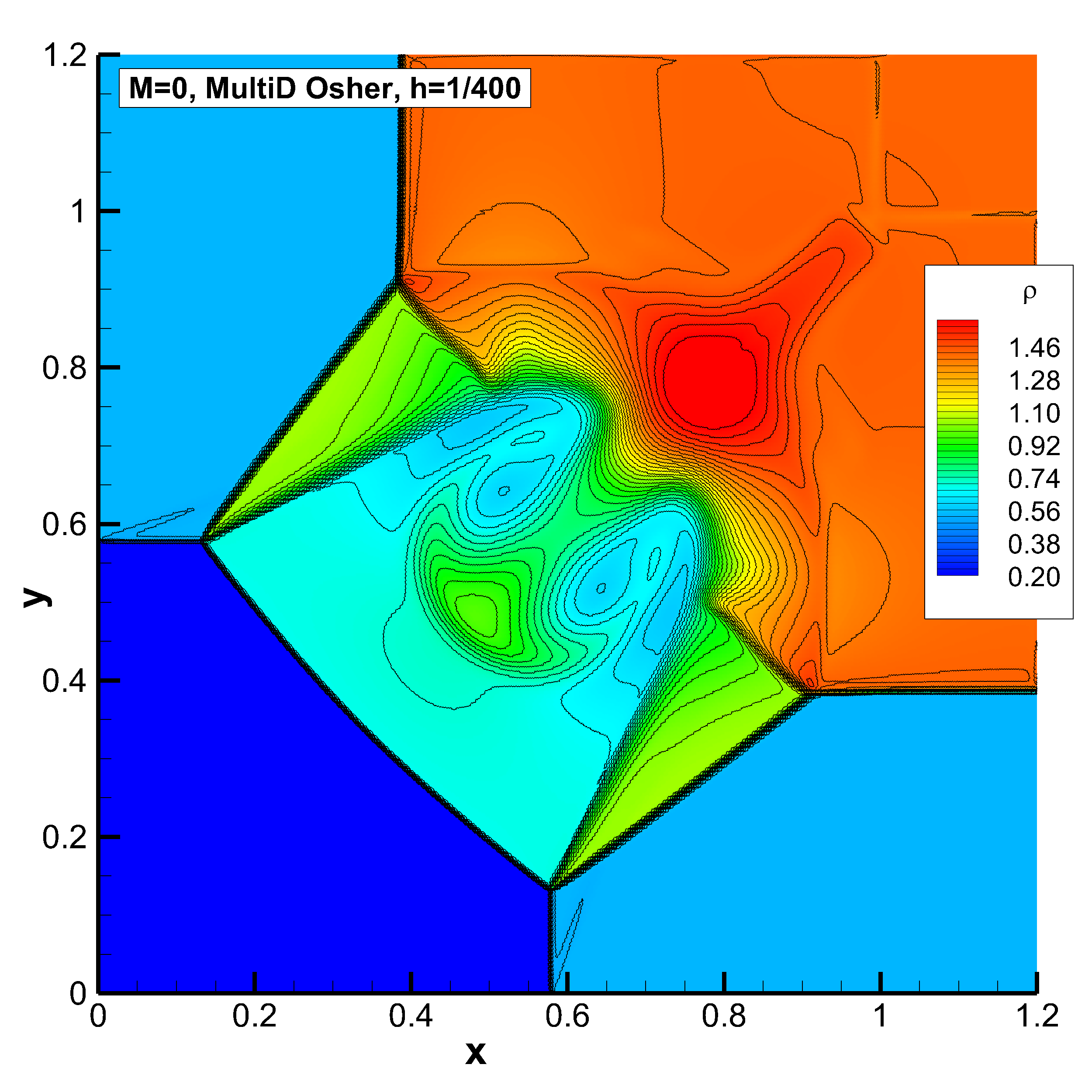}%
		\includegraphics[width=0.333\linewidth,trim=1 1 1 1,clip]{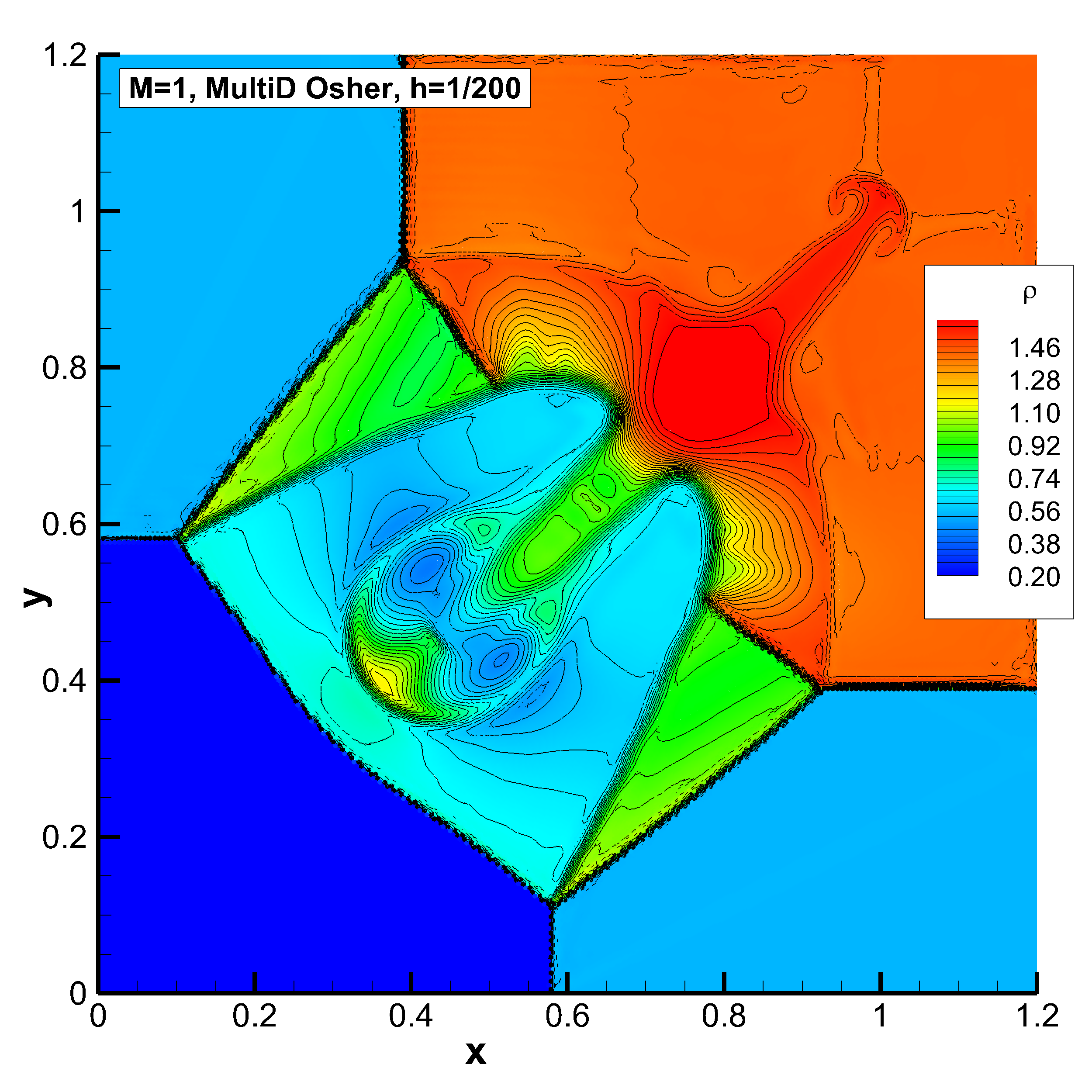}%
		\includegraphics[width=0.333\linewidth,trim=1 1 1 1,clip]{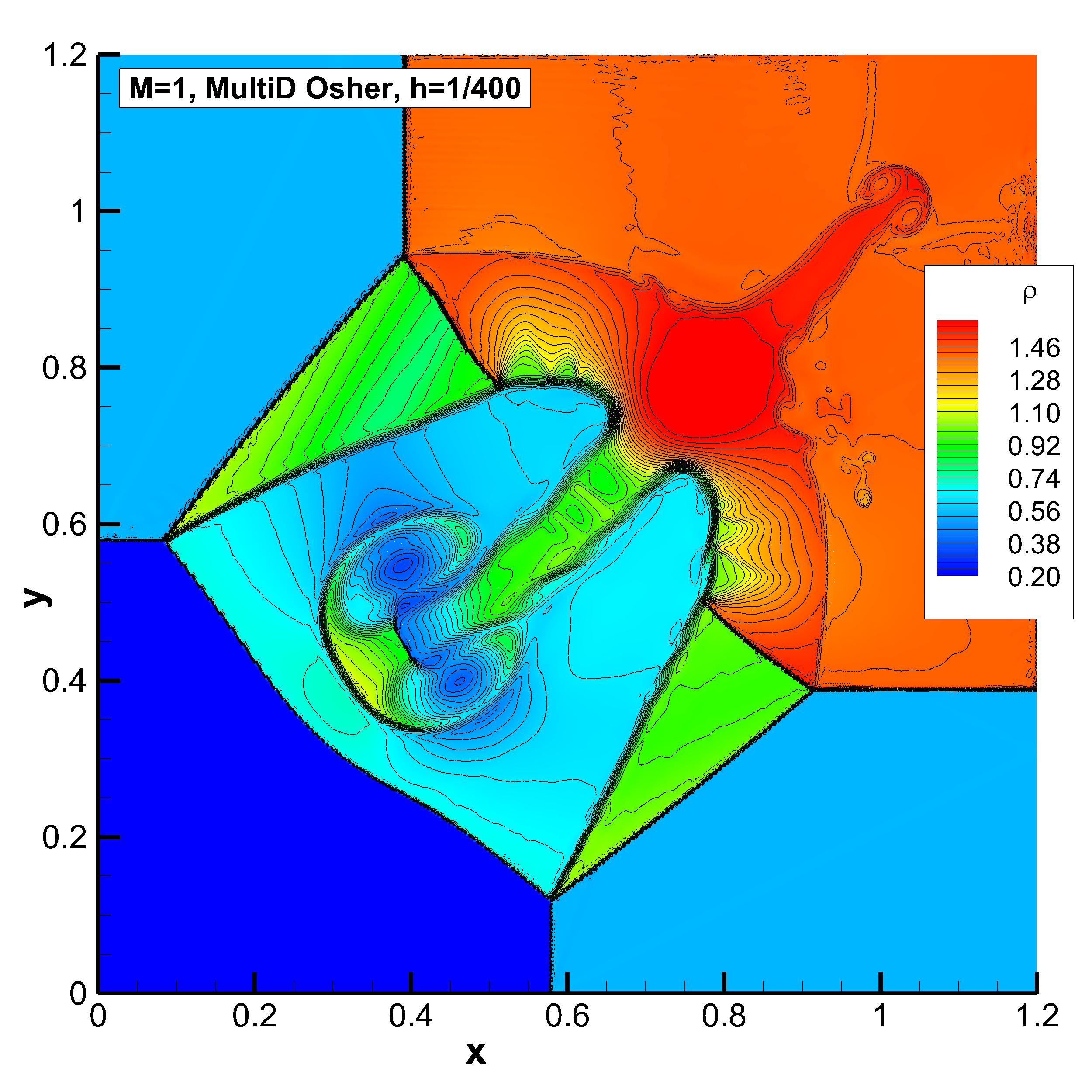}\\    
		\includegraphics[width=0.333\linewidth,trim=1 1 1 1,clip]{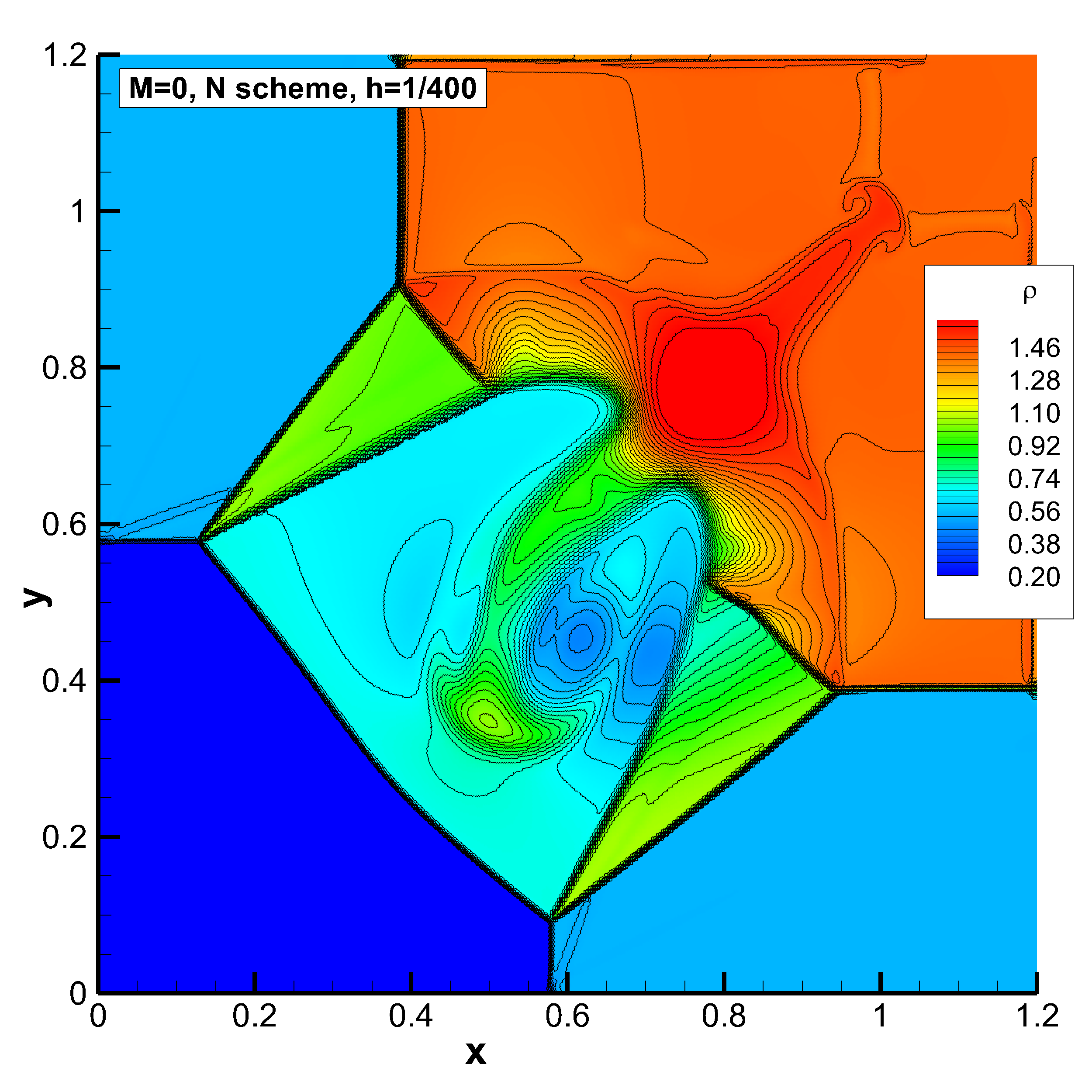}%
		\includegraphics[width=0.333\linewidth,trim=1 1 1 1,clip]{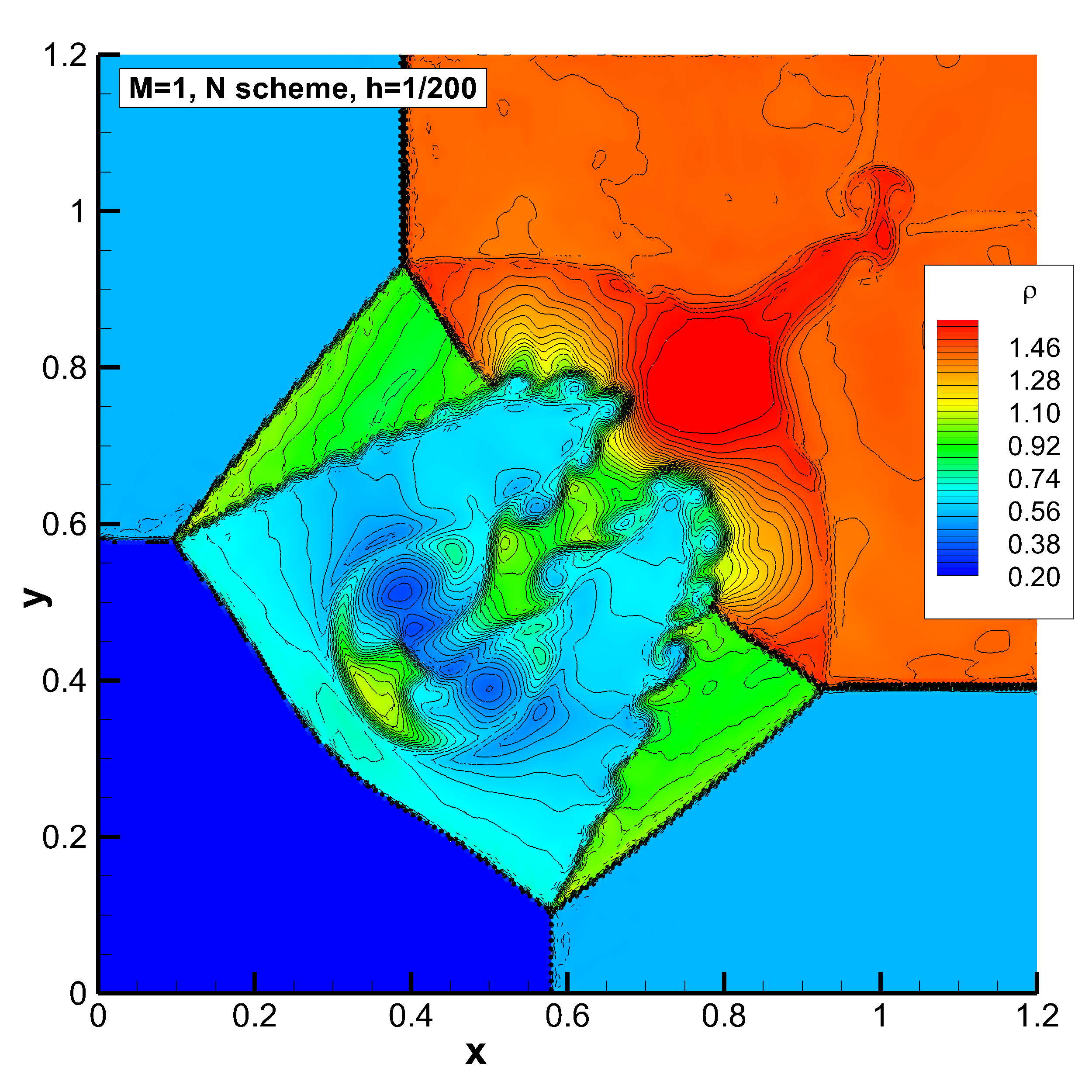}%
		\includegraphics[width=0.333\linewidth,trim=1 1 1 1,clip]{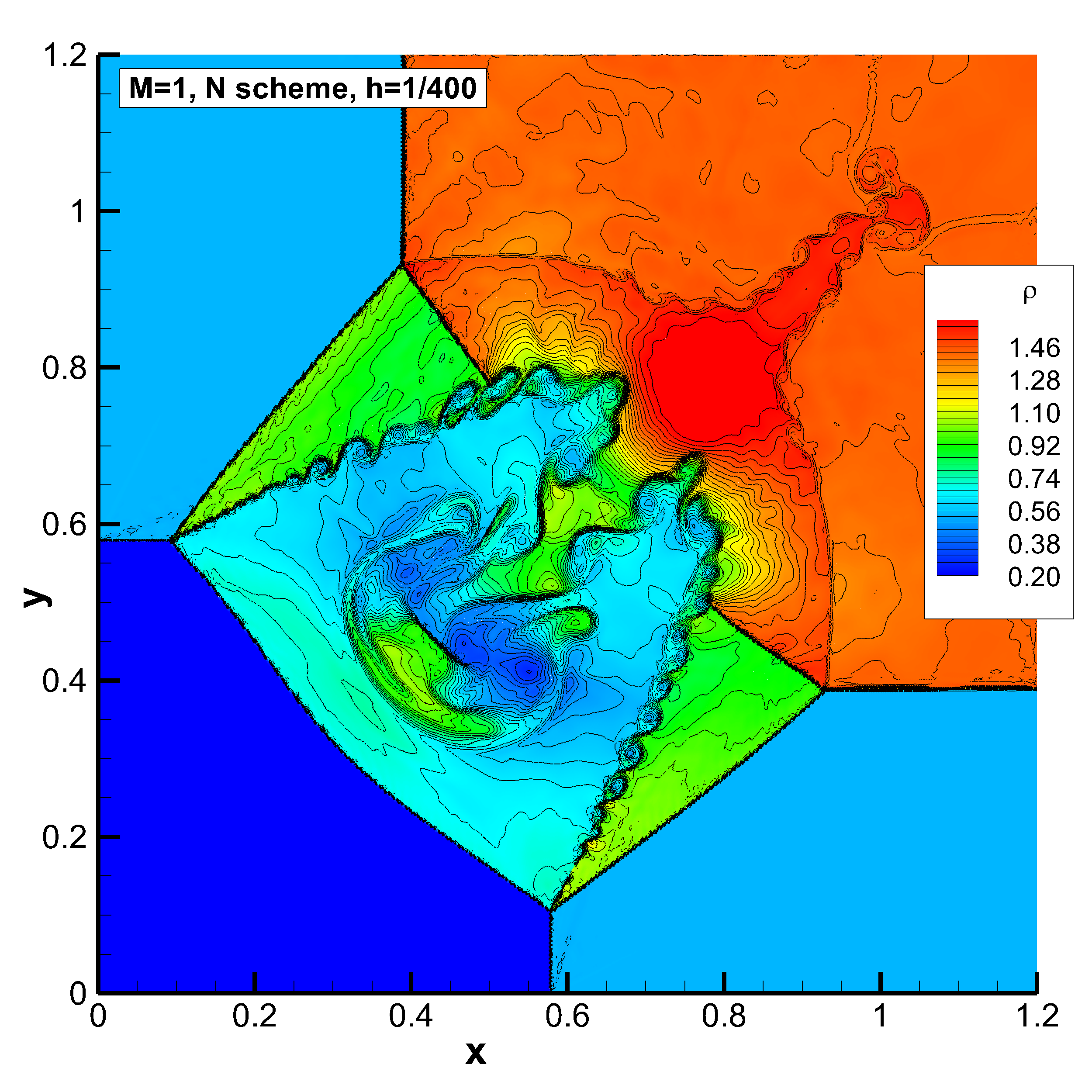}\\         
		\caption{Two dimensional Riemann problem (configuration 3 of~\cite{laxliu98}). 
			Here we report the numerical results obtained with our first order \RIIcolor{(M=0)} and second order  \RIIcolor{(M=1)} schemes over a coarse mesh with $\RIIcolor{h\simeq1/200}$ and a finer one with $h\simeq1/400$.
			In particular, we show the density profile with 40 contours line in the interval $[0.2, 1.6]$.
			We can observe  that the    
			N scheme, even if the carbuncle fix is active,     
			clearly captures the shear instabilities, since the carbuncle fix is active only across shock waves.
			\RIcolor{Also, thanks to the carbuncle fix and the truly multidimensional behavior of this Riemann solver, not dominated by excessive numerical dissipation, the N scheme shows less mesh imprinting with respect to the other solvers}. 
			Finally, we recall that the asymmetry of the results is due to the use of fully  unstructured meshes.
		}
		\label{fig.2drp_3}
	\end{figure}

	We report the results obtained with our first and second order finite volume schemes in Figure~\ref{fig.2drp_3} so to compare the different resolutions achieved with the classical Roe edge solver and with our multidimensional corner fluxes. 
	
	\subsection{Hypersonic flow past a blunt body} 
	
	\begin{figure}[!b]
		\centering
		\includegraphics[width=0.245\linewidth,trim=1 1 1 1,clip]{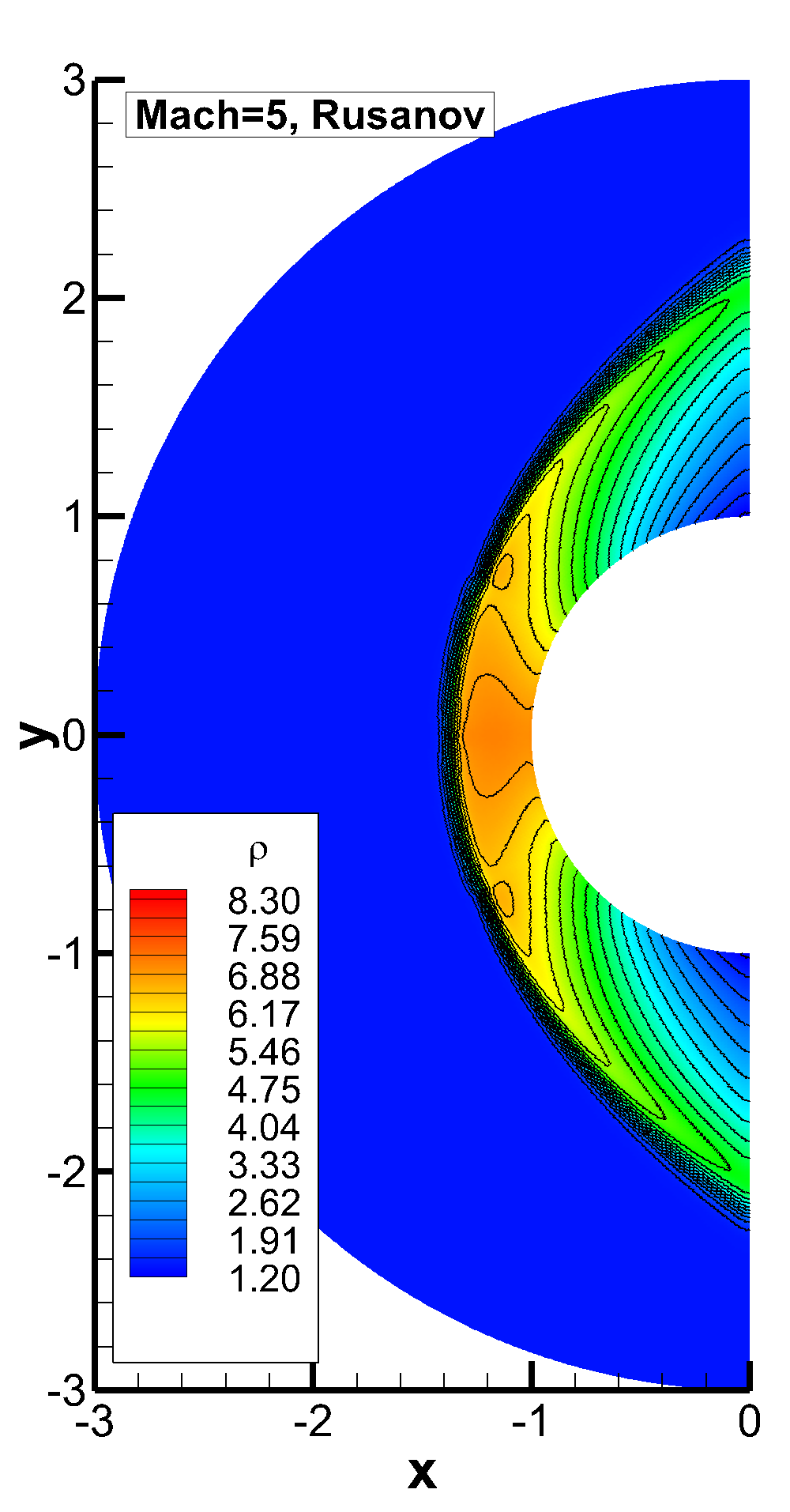}%
		\includegraphics[width=0.245\linewidth,trim=1 1 1 1,clip]{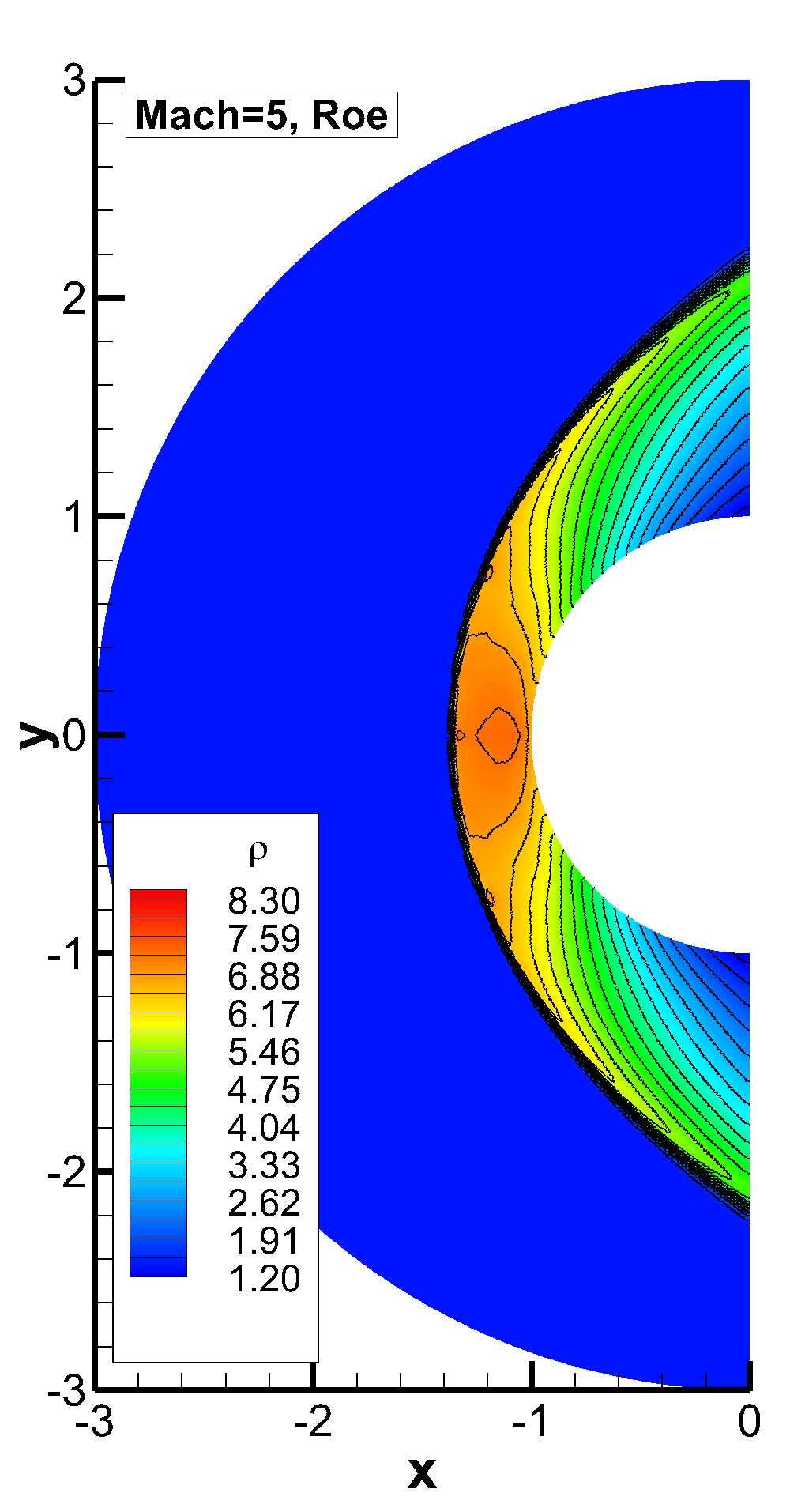}%
		\includegraphics[width=0.245\linewidth,trim=1 1 1 1,clip]{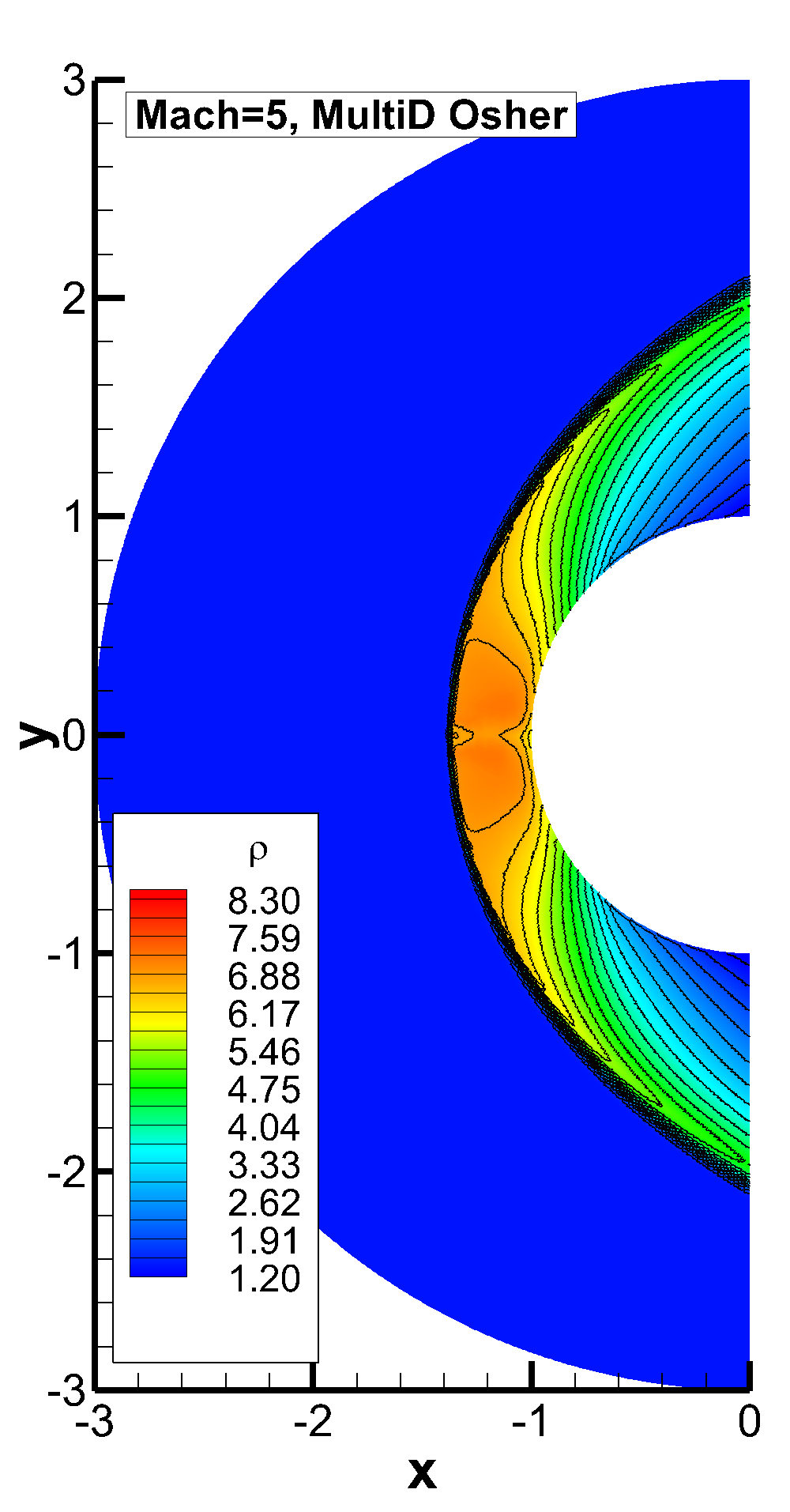}%
		\includegraphics[width=0.245\linewidth,trim=1 1 1 1,clip]{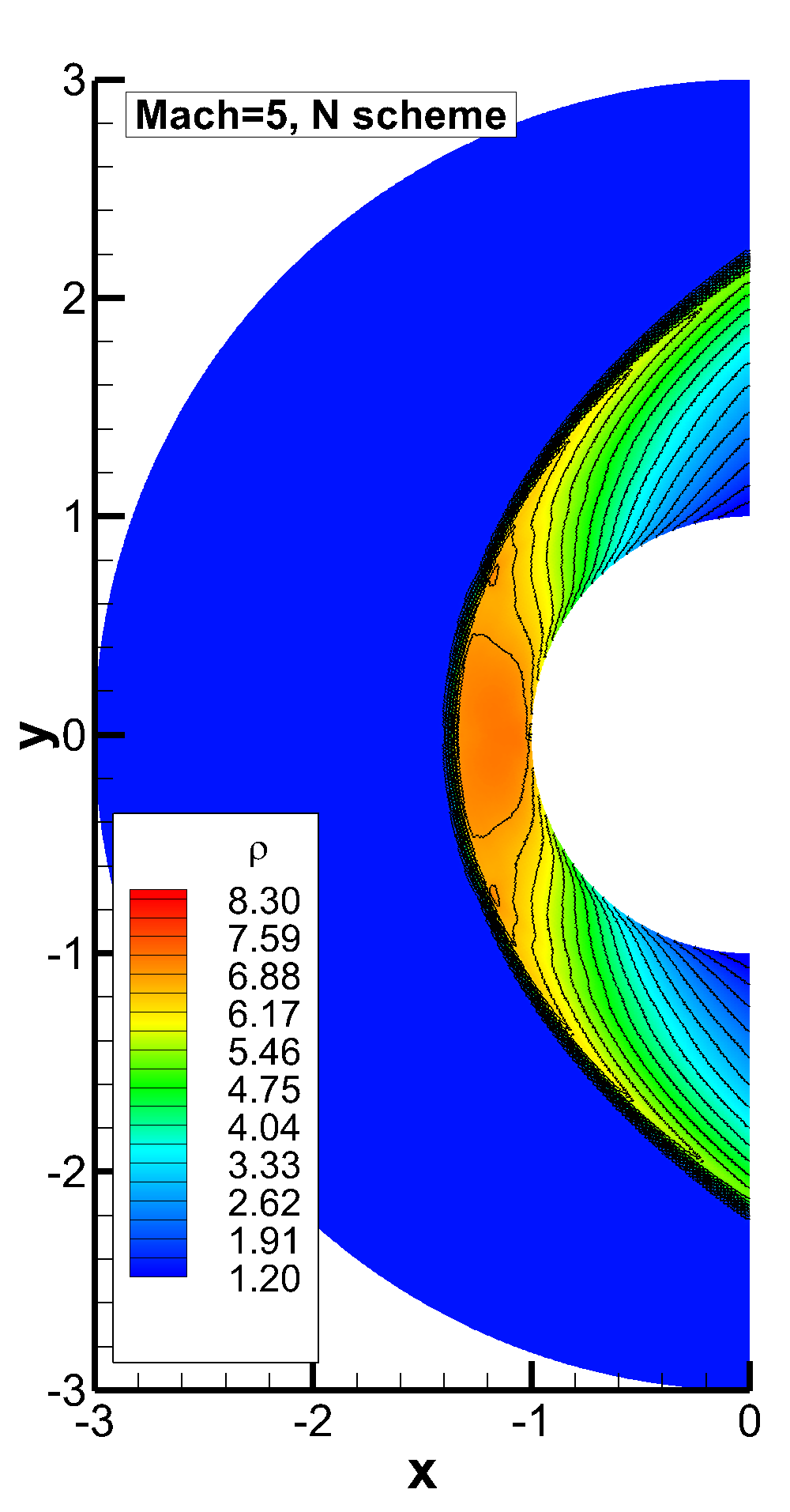}\\
		\includegraphics[width=0.245\linewidth,trim=1 1 1 1,clip]{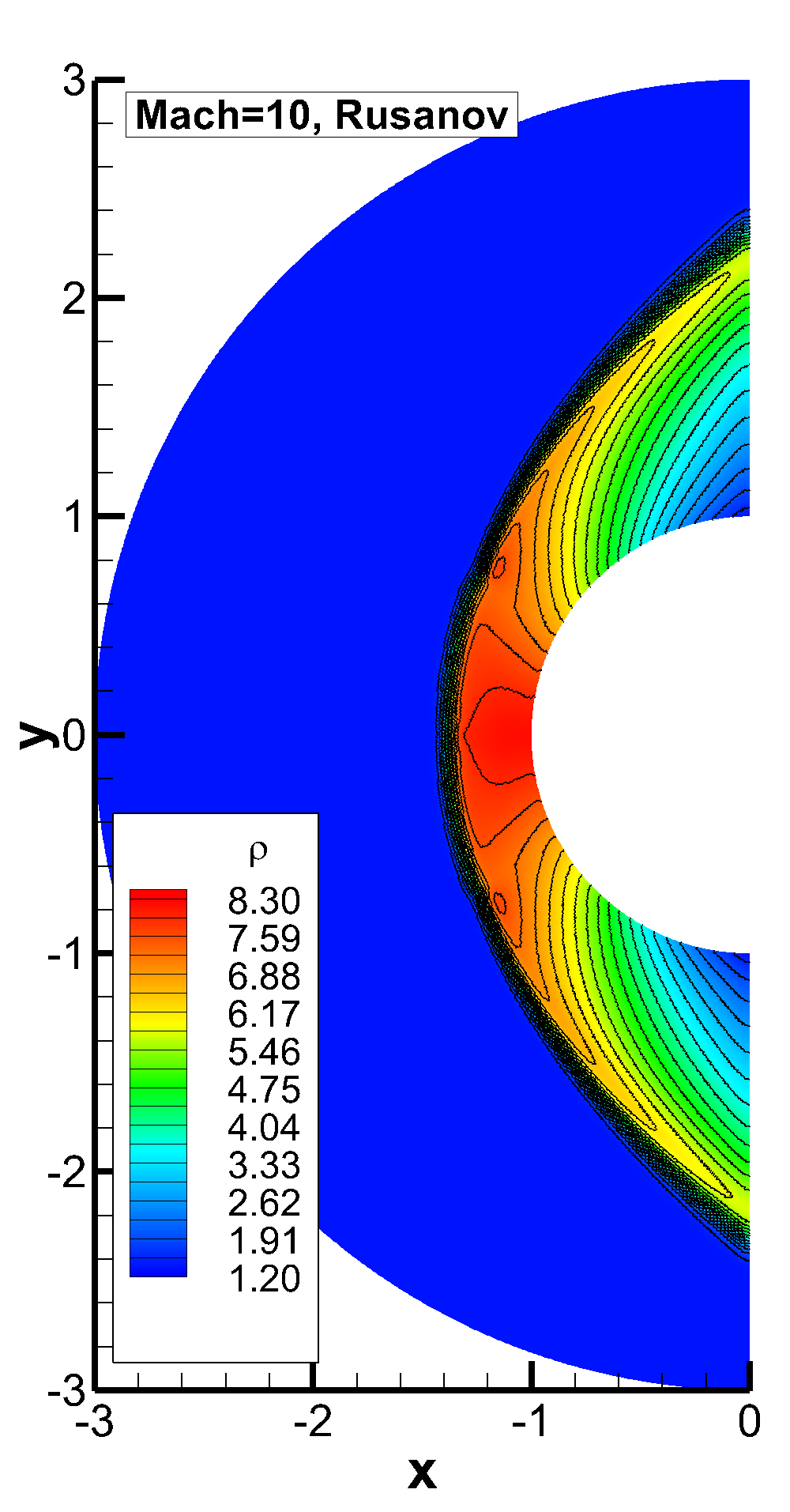}%
		\includegraphics[width=0.245\linewidth,trim=1 1 1 1,clip]{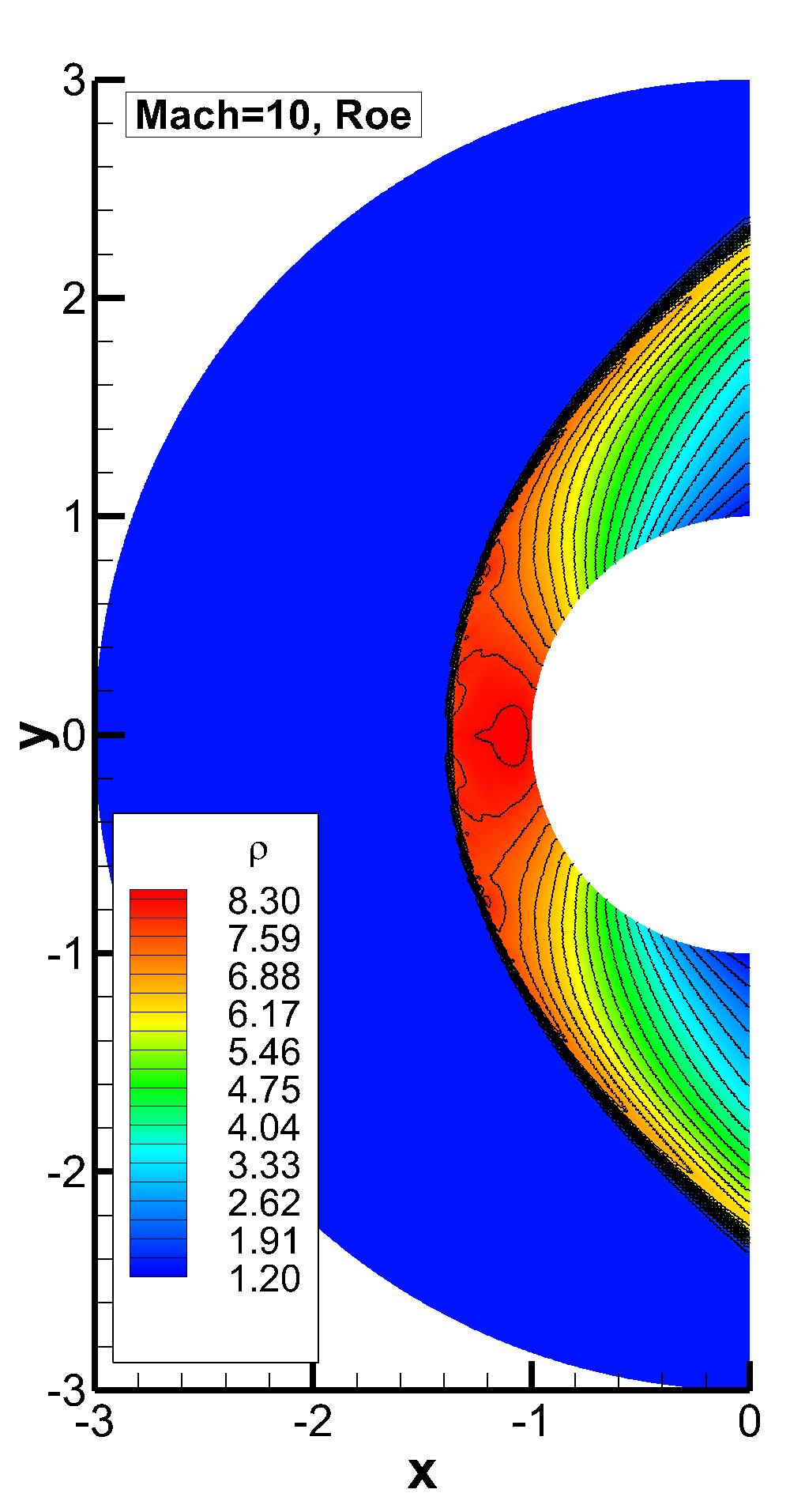}%
		\includegraphics[width=0.245\linewidth,trim=1 1 1 1,clip]{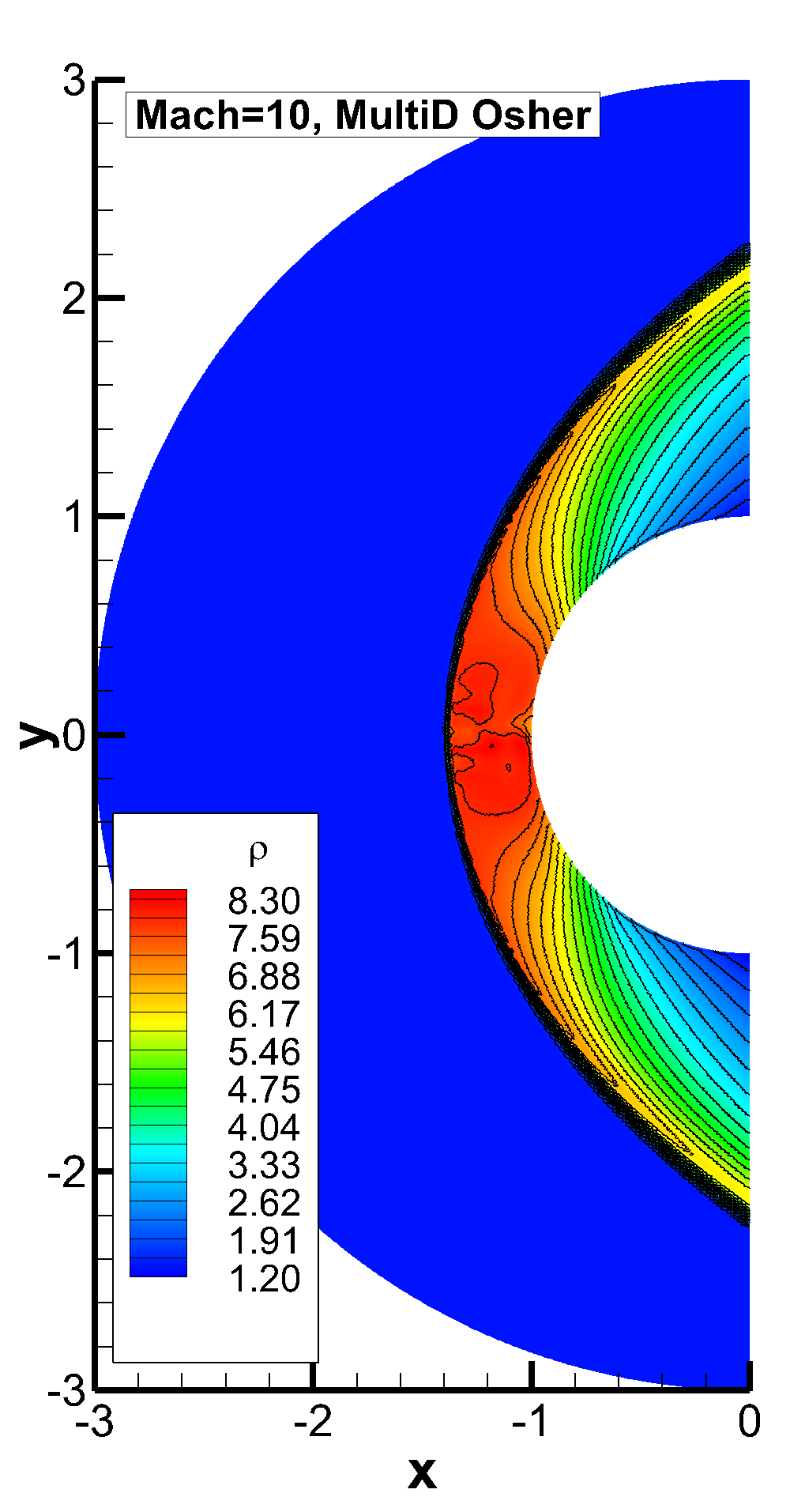}%
		\includegraphics[width=0.245\linewidth,trim=1 1 1 1,clip]{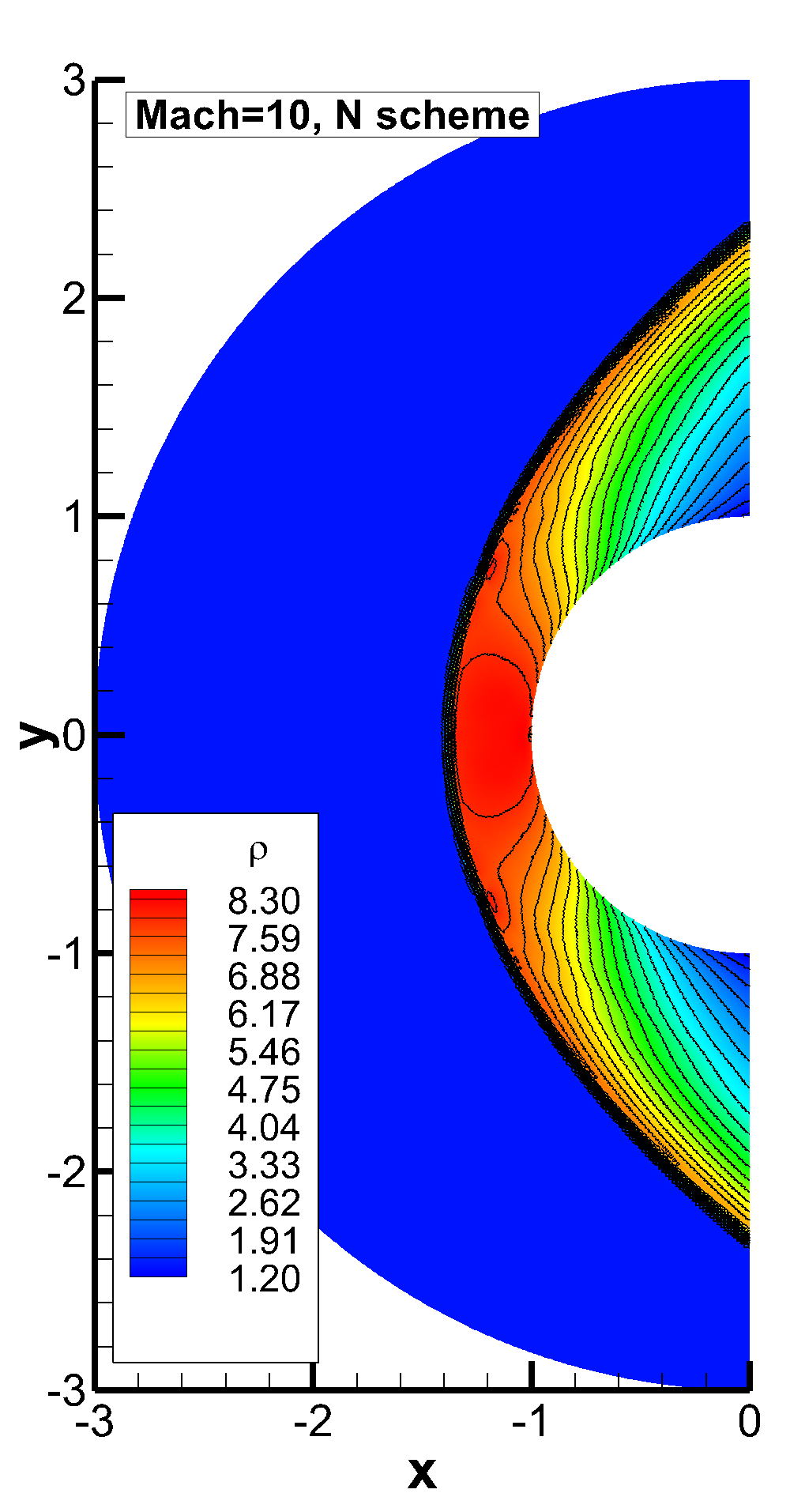}\\
		\caption{Here, we show the density profile of the blunt body test case at Mach=5 (top line) and Mach=10 (bottom line) with 21 contours line in $[1.2,8.3]$.  
			In particular, we compare the behavior of our multidimensional solvers with the classical 1d edge solvers inserted in a first order \RIIcolor{(M=0)} finite volume scheme. 
			We first notice that the Roe scheme, which is usually affected by the carbuncle phenomenon on triangular meshes, here on polygonal meshes does not \RIIcolor{suffer from} this problem.
			The incomplete simple Rusanov scheme is, as usual, carbuncle-free but quite dissipative. 
			Instead, we do not report any results for the 1d Osher scheme because it crashes at the very beginning of the simulation due to the early appearance of the carbuncle phenomenon.
			On the contrary, our multidimensional extension of the Osher scheme, without the \RIIcolor{need for} any trick, provides \RIIcolor{highly} resolved results both at Mach=5 and Mach=10 \RIcolor{with just a small carbuncle effect at high Mach  
				that slightly breaks the symmetry but does not affect the simulation}. 
			Finally, the N scheme, which we recall is implemented with a well-known carbuncle fix, is completely carbuncle free and shows a reduced dissipation despite the employed fix. 
		}
		\label{fig.bluntbodyP0-comparison}
	\end{figure}
	
	We consider now a hypersonic (Mach = 5 and Mach = 10) flow past the forebody of a circular cylinder~\cite{kemm2020simple,PHRaph1,ciallella2023shifted}. As computational domain we take the half ring 
	between two circles of internal radius $r_0=1$ and external $r_1=3$ (centered at the origin) and with $x<0$. We discretize it with a coarse mesh of characteristic size $h\simeq 1/70$.
	The initial condition is given by an incoming hypersonic flow such that 
	$(\rho, u, v, p)(\x) = (\sqrt{\gamma}, \text{Mach}, 0, 1)$ and the boundary conditions are of wall type on $r_0$, of transmissive type on the right vertical boundaries and of Dirichlet type on $r_1$.
	
	We report the obtained results in Figure~\ref{fig.bluntbodyP0-comparison} together with detailed comments regarding in particular the behavior of our schemes w.r.t the well-known appearance of the carbuncle phenomenon.

	\subsection{Circular two-dimensional explosion with boundary reflections}

	\begin{figure}[!b]
		\centering
		\includegraphics[width=0.495\linewidth,trim=1 1 1 1,clip]{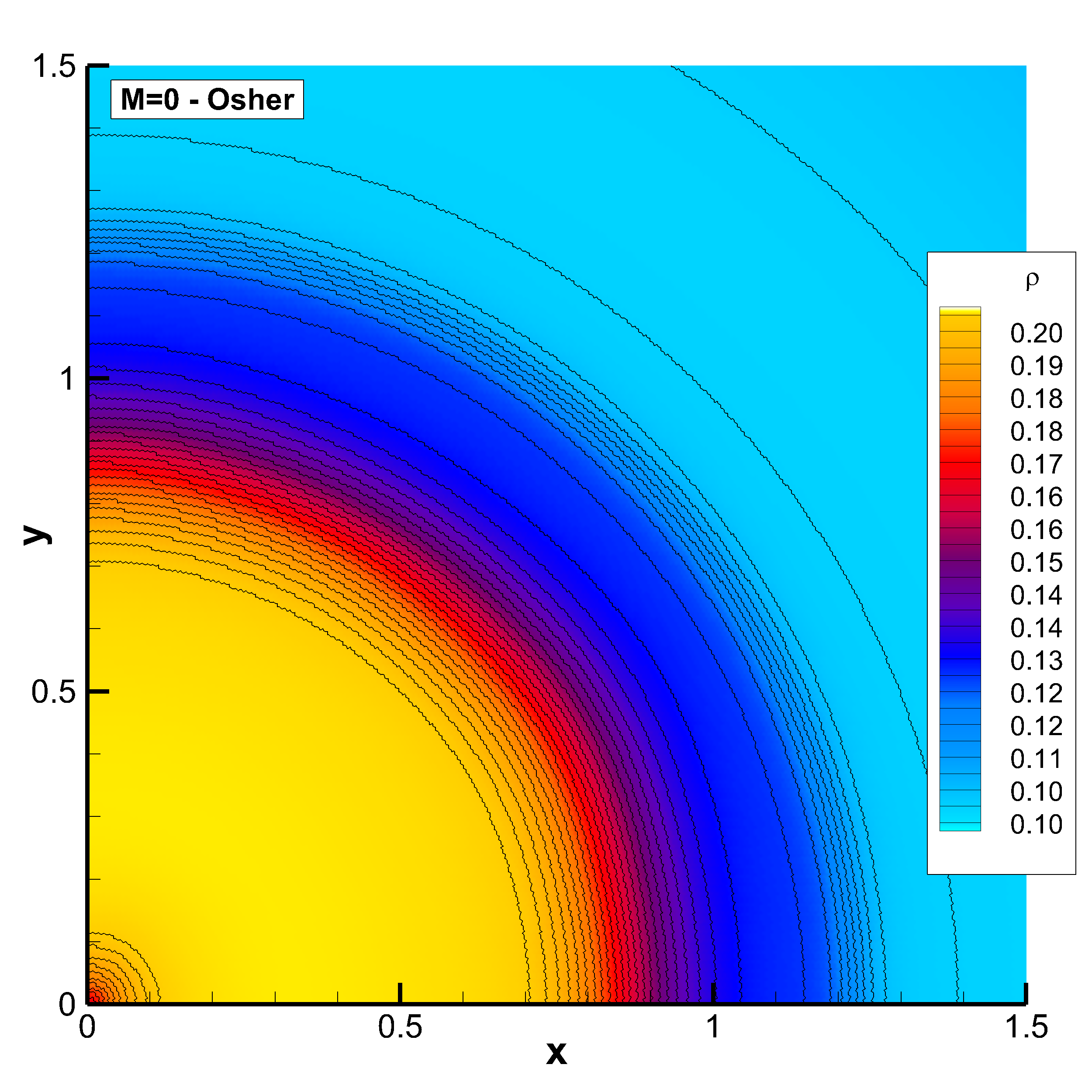}%
		\includegraphics[width=0.495\linewidth,trim=1 1 1 1,clip]{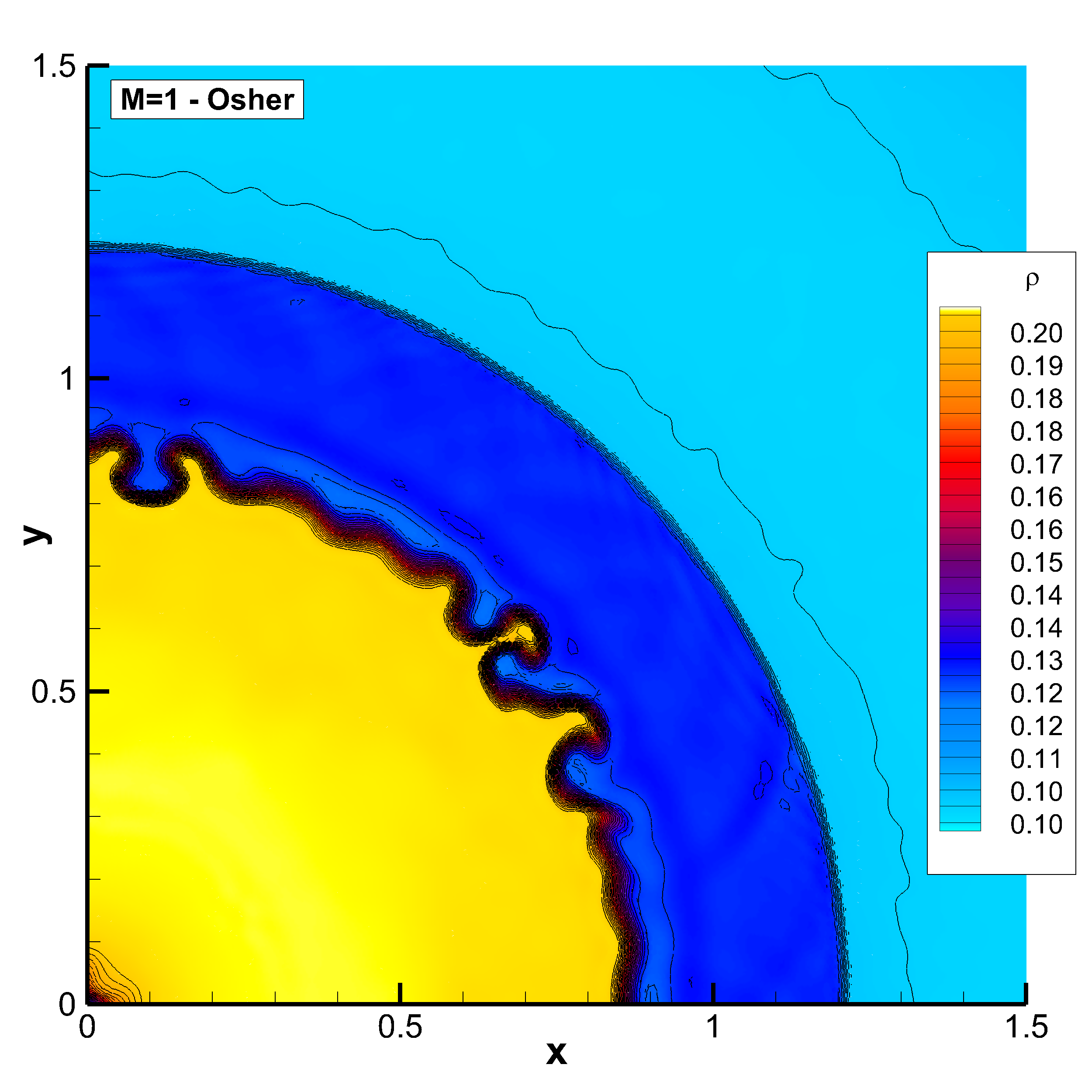}\\
		\includegraphics[width=0.495\linewidth,trim=1 1 1 1,clip]{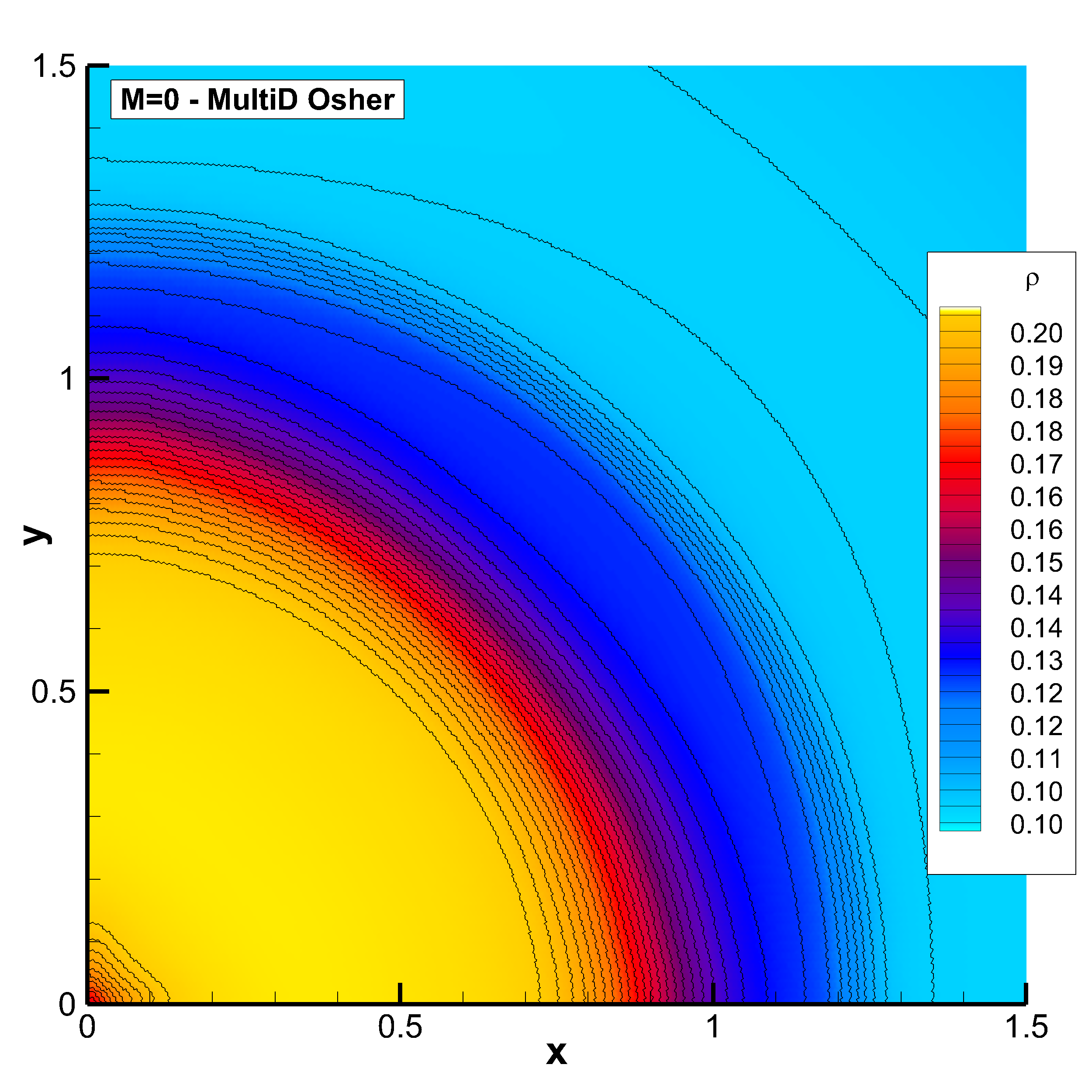}%
		\includegraphics[width=0.495\linewidth,trim=1 1 1 1,clip]{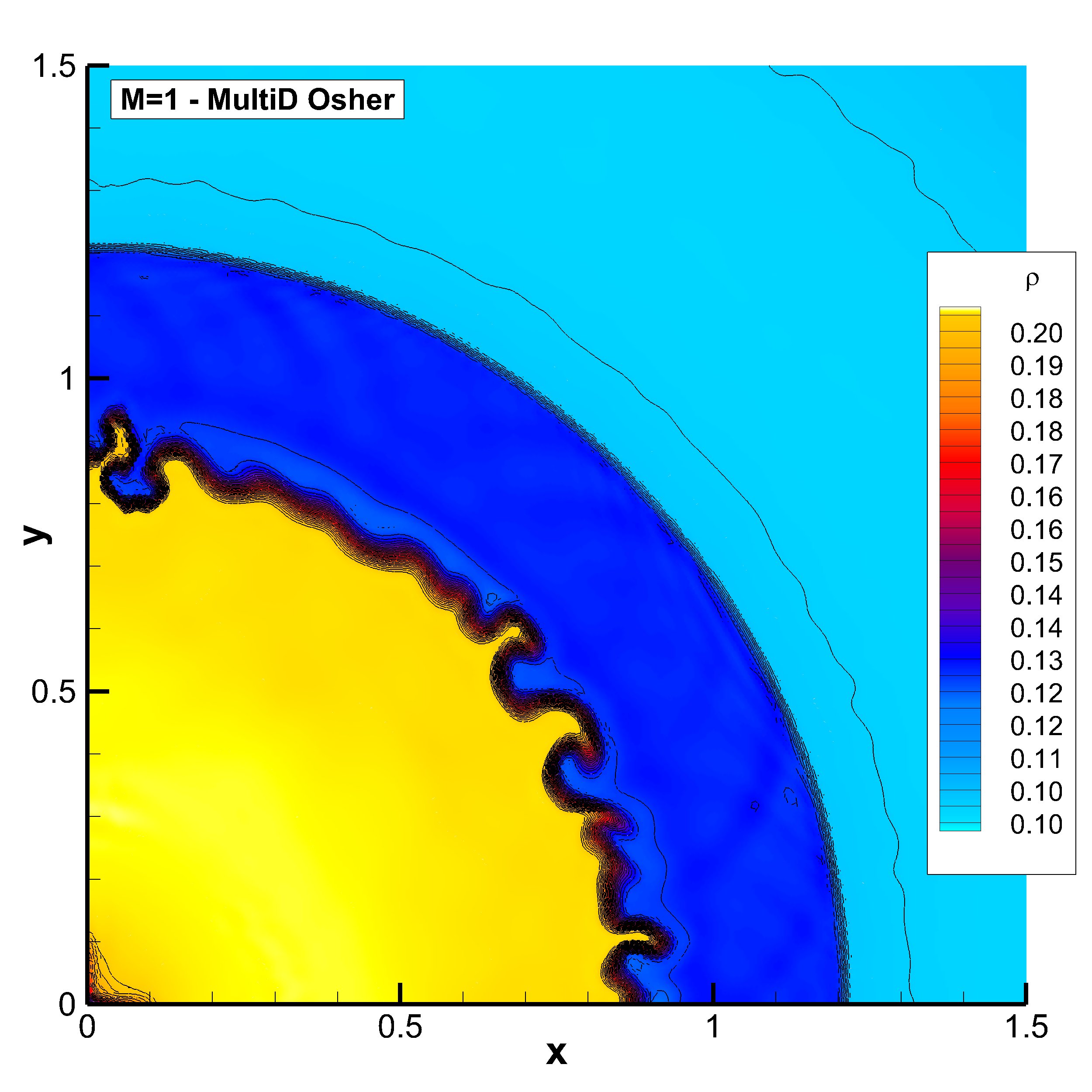}\\
		\caption{\RIIcolor{Circular two-dimensional explosion with inner reflections. We compare the Osher type fluxes at first and second order of accuracy. In particular, we show the density contours at time $t_f = 3.2$ with 16 contour levels in $[0.1,0.2]$.
				On this specific test case, we can only observe small differences between the 1d and the MultiD Osher schemes: in particular, the second has less angular symmetry and more instabilities near the boundaries. This is an indication  of the less dissipative nature of our multidimensional scheme.
		}}
		\label{fig.exp_reflections_Osher}
	\end{figure}
	
	\begin{figure}[!b]
		\centering
		\includegraphics[width=0.495\linewidth,trim=1 1 1 1,clip]{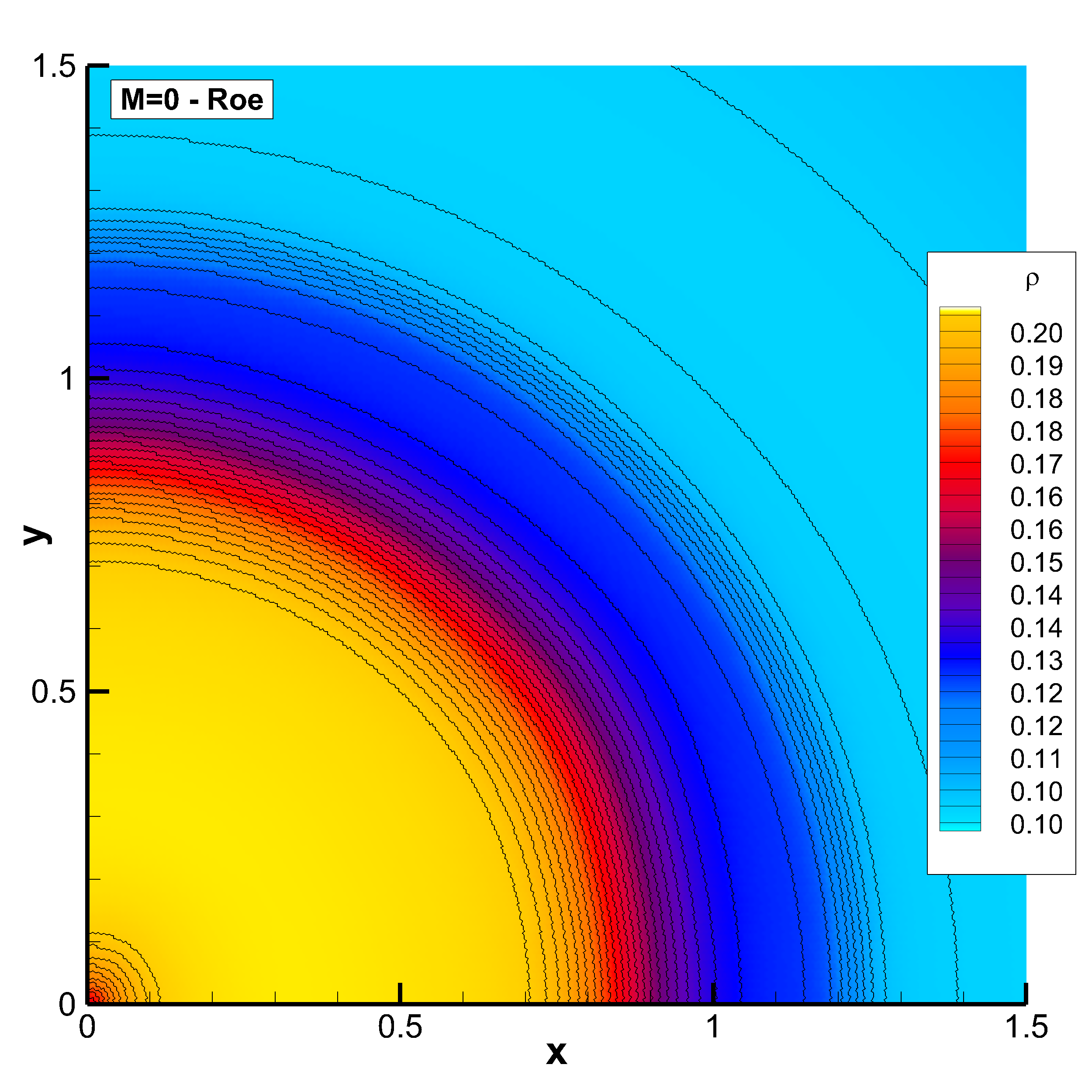}%
		\includegraphics[width=0.495\linewidth,trim=1 1 1 1,clip]{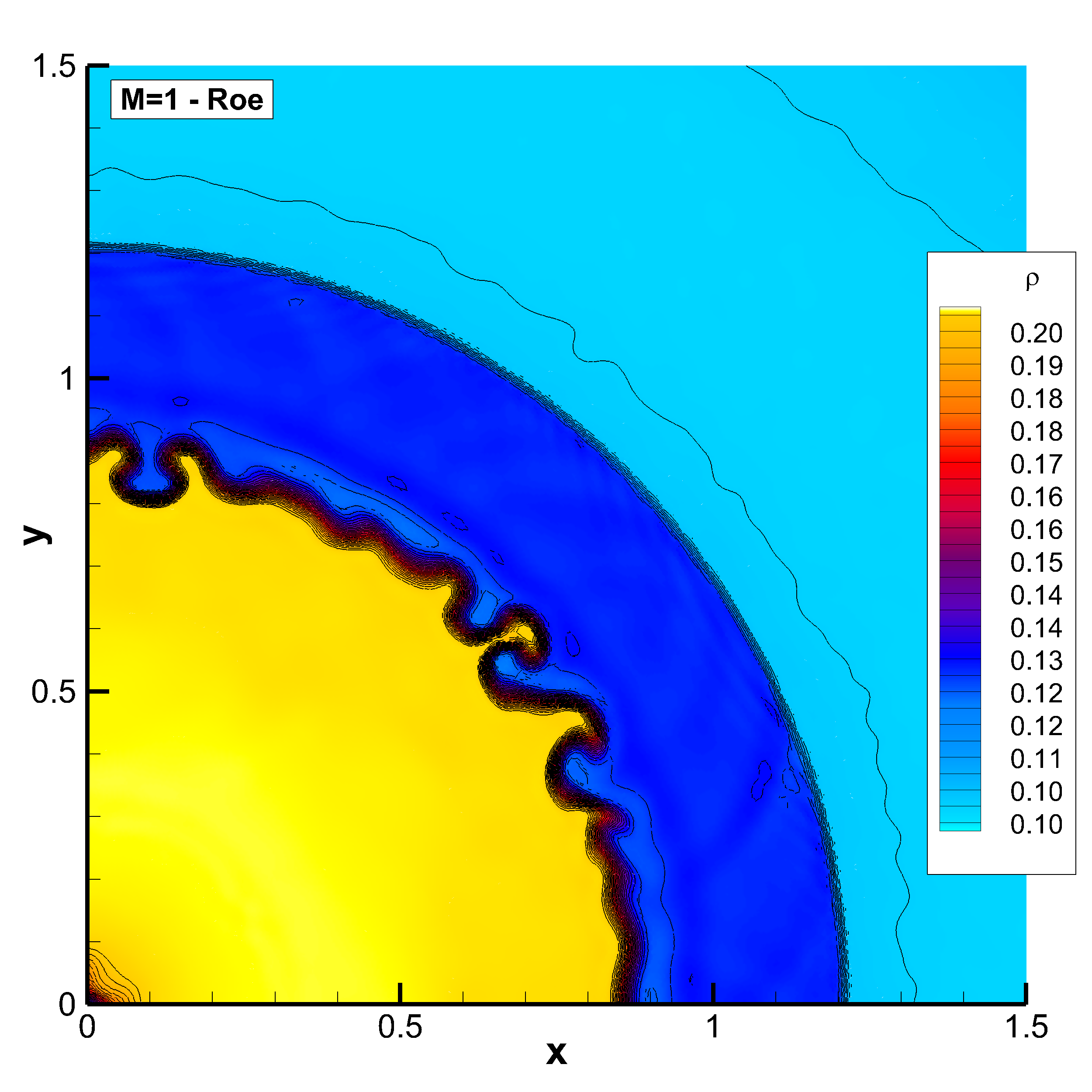}\\
		\includegraphics[width=0.495\linewidth,trim=1 1 1 1,clip]{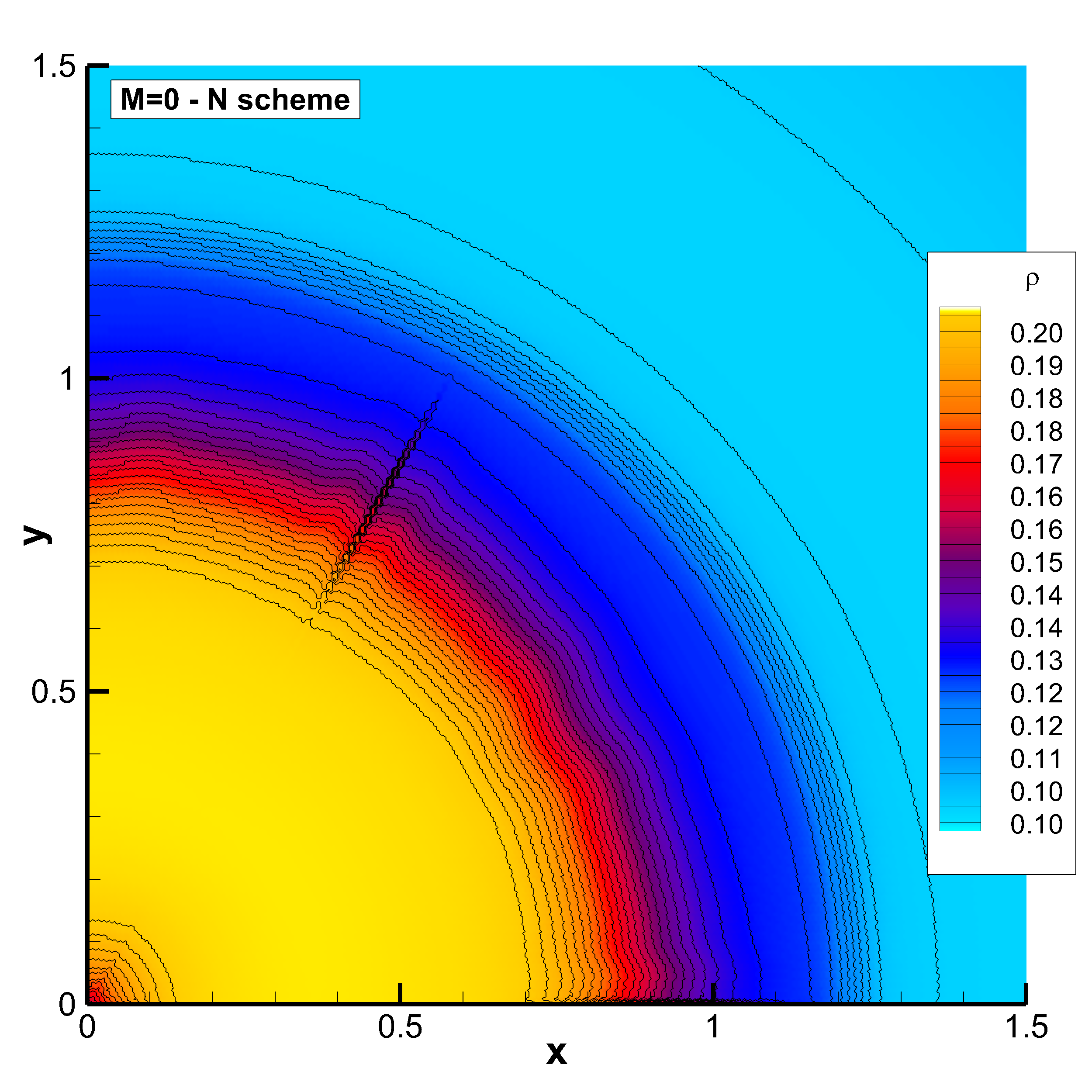}%
		\includegraphics[width=0.495\linewidth,trim=1 1 1 1,clip]{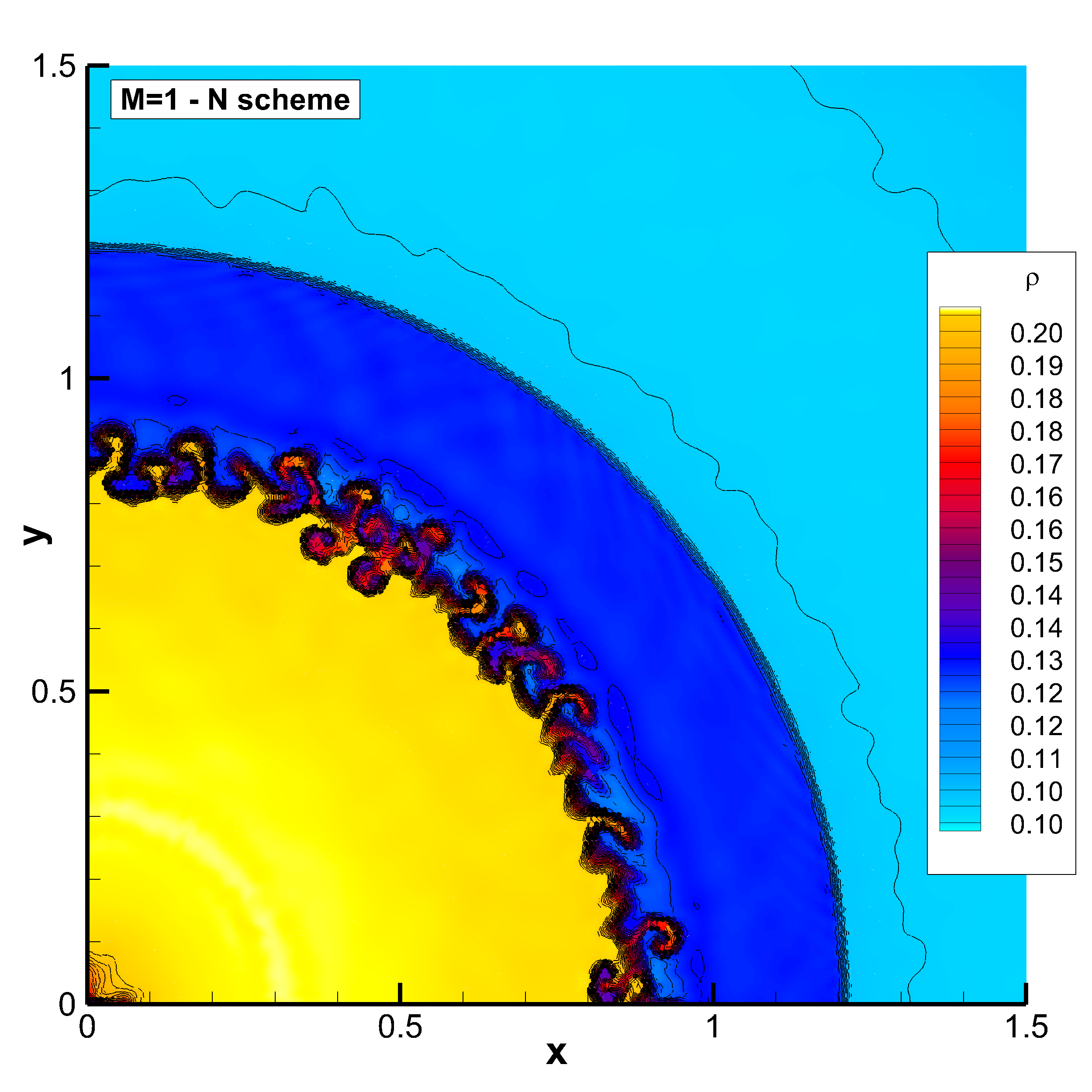}%
		\caption{\RIIcolor{Circular two-dimensional explosion with boundary reflections, comparison of Roe and N scheme solvers at first and second order of accuracy. In particular, we show the density contours at time $t_f = 3.2$ with 16 contours level in $[0.1,0.2]$.
				We note that the N scheme shows a significantly higher resolution compared to all other schemes, highlighting all the instabilities that develop along the discontinuity located at $r \simeq 1.75$.
		}}
		\label{fig.exp_reflections_Roe-Nscheme}
	\end{figure}
	
	We conclude our sets of benchmarks with an explosion problem which has an initial configuration similar to the Sod test case but it is performed on a larger domain and for a longer time, so to see \RIIcolor{many more reflections and interactions around the origin}. The final computational time is \RIIcolor{$t_f=3.2$}.
	The initial conditions for this problem are 
	\begin{equation}
		(\rho, u, v, p)(\x) = 
		\begin{cases} 
			(1.0  , 0.0, 0.0,  1.0   )   &  \ \text{ if } \  x > 0.4, \\
			(0.125,  0.0, 0.0, 0.1 )  &   \ \text{ if } \  x \le 0.4, \\	
		\end{cases}
	\end{equation}
	and \RIIcolor{the domain of interest is $\Omega = [0,1.5]\times[0,1.5]$. 
		(We have solved the problem on a larger domain $\Omega = [-3.5,3.5]\times[-3.5,3.5]$ to avoid the effect of external boundary conditions).}
	
	We have simulated this test case with our first and second order finite volume schemes, and we report in Figure~\ref{fig.exp_reflections_Osher} the results obtained with the Osher solvers and in Figure~\ref{fig.exp_reflections_Roe-Nscheme} those obtained with the Roe scheme and the N scheme.

	\section{Conclusions and outlook to future work}
	\label{sec.conclusions}
	
	In this paper we have designed new high order accurate fully-discrete one-step finite volume schemes that are cell centered and employ nodal fluxes. High order of accuracy in space is achieved via the CWENO reconstruction procedure and high order of accuracy in time is obtained at the aid of the ADER approach, making use of a weak formulation of the PDE system in space-time. With respect to existing high order cell-centered finite volume schemes that employ vertex-based numerical fluxes this paper introduces two novelties: 
	\begin{enumerate} 
		\item the formulation of a new genuinely multidimensional Osher-type Riemann solver for general nonlinear systems of hyperbolic conservation laws, see~\eqref{eqn.multid.osher.final} for the algebraic expression of the numerical flux and~\eqref{eqn.multid.path} for the generalization of the simple 1d integration path used in~\cite{OsherNC,OsherUniversal} to a manifold in multiple space dimensions. The very peculiar feature of this new vertex-based Riemann solver is the fact that it produces an entire numerical flux tensor rather than just a numerical flux projected into a particular direction. The mathematical structure of the flux is a simple and direct extension of the one-dimensional Osher-type solvers introduced in~\cite{OsherNC,OsherUniversal}. It consists in a central part of the flux in which simply the arithmetic averages of the flux tensors surrounding a vertex is computed. The multidimensional numerical dissipation is then achieved by computing a surface integral of the matrix absolute value operator applied to the flux Jacobians in each spatial direction over a virtual simplex element constituted by the cell centers around each vertex. The integral is computed numerically via suitable quadrature rules. These integrals are then multiplied by a discrete multidimensional gradient that is simply obtained via the Gauss theorem, as already proposed in~\cite{BarthJespersen} to obtain discrete gradients on general unstructured meshes.       
		\item the careful rewriting of well-known fluctuations from the residual distribution (RD) framework allows \RIIcolor{constructing} suitable approximate multidimensional Riemann solvers that can then be incorporated inside classical high order cell-centered WENO finite volume schemes. It is therefore very interesting to note that RD fluctuations are not limited to the RD context, but can be used in a much more general setting. The most prominent example of genuinely multidimensional upwinding that accounts for all characteristic fields present in the nonlinear hyperbolic system is perhaps the N scheme. Its rewriting as vertex flux in the context of multidimensional Riemann solvers has been achieved in Equations~\eqref{eq.N-rd}-\eqref{eq.N-rd3} in combination with the general definition~\eqref{eq:point_flux_phi} that relates existing RD fluctuations to numerical point fluxes. 
	\end{enumerate}
	\RIIcolor{Future work will first focus on a more in-depth analysis of our approach, with the aim of providing a formulation that, for schemes of order higher than two, differs from the one proposed here and relies solely on multidimensional solvers. Indeed, our current formulation for orders higher than two mixes the use of 1d and multidimensional solvers, which allows reaching arbitrarily high accuracy but introduces 1d dissipation that we aim to avoid. 
		Additionally, we plan to explore the application of these new multidimensional Riemann solvers in the context of exactly divergence-free schemes for MHD and more complex hyperbolic PDE systems, as well as their extension to non-conservative hyperbolic systems.}

	\section{Acknowledgments}
	
	E.~Gaburro gratefully acknowledges the support received from the European Union 
	with the ERC Starting Grant \textit{ALcHyMiA} (grant agreement No. 101114995).
	
	M. Dumbser was funded by the Italian Ministry of Education, University and Research (MIUR) in the framework of the PRIN 2022 project \textit{High order structure-preserving semi-implicit schemes for hyperbolic equations}, via the  Departments of Excellence  Initiative 2018--2027 attributed to DICAM of the University of Trento (grant L. 232/2016) and  
	by the European Union Next Generation EU projects PNRR Spoke 7 CN HPC and SMART. M. Dumbser has also received funding by the Fondazione Caritro via the project SOPHOS. 
	
	Views and opinions expressed are however those of the authors only and do not necessarily reflect those of the European Union or the European Research Council Executive Agency. 
	Neither the European Union nor the granting authority can be held responsible for them.
	
	E.~Gaburro and M.~Dumbser are members of the INdAM GNCS group in Italy. 
	M.~Ricchiuto is member of the CARDAMOM team at the Inria center of the University of Bordeaux. 
	
	\noindent \section*{In memoriam}
	
	\noindent This paper is dedicated to the memory of Prof. Arturo Hidalgo L\'opez ($^*$July 3\textsuperscript{rd} 1966 - $\dagger$August 26\textsuperscript{th} 2024) of the Universidad Politecnica de Madrid, organizer of HONOM 2019, active participant in many other editions of HONOM and, most of all, a very dear friend. 
	Our thoughts and wishes go to his wife Lourdes and his sister Mar\'ia Jes\'us, whom he left behind.
	
	\appendix
	
	\section{Entropy conservation:  Tadmor's shuffle relation and point entropy fluctuation\label{app.entropy}}

	We now consider the entropy evolution for the  multidimensional scheme with  corner fluxes 
	$$
	|\Omega_c|\dfrac{d\Q_c}{dt} + \sum\limits_p \hat\F_p \cdot\n_{pc} =0.
	$$
	Let $(\eta(\Q), \mathbf{G}(\Q) )$ denote the entropy/entropy flux pair, with as usual
	$$
	\mathbf{G}(\Q)  =  \mathbf{W}^t \F(\Q)  - \boldsymbol{\Psi}(\Q)
	$$
	with $ \mathbf{W}^t  = \nabla_{\Q}\eta$ the entropy variables, and $\boldsymbol{\Psi}$ the entropy potential vector.\\
	
	We first generalize Tadmor's shuffle condition. To this end, we compute
	$$
	\sum_c|\Omega_c| \mathbf{W}^t_c\dfrac{d\Q_c}{dt} + \sum_c\sum\limits_{p\in \mathcal{P}_c} \mathbf{W}^t_c\hat\F_p \cdot\n_{pc} =0,
	$$
	with $\mathbf{W}_c=\mathbf{W}(\Q_c)$ the entropy variables vector. By definition we have
	$$
	\sum_c|\Omega_c|  \dfrac{d\eta(\Q_c)}{dt} + \sum_c\sum\limits_{p\in \mathcal{P}_c} \mathbf{W}^t_c\hat\F_p \cdot\n_{pc} =0.
	$$
	Following~\cite{tadmor} we now write at a given $p\in \mathcal{P}_c$
	$$
	\mathbf{W}_c = \overline{\mathbf{W}}_p + (  \mathbf{W}_c - \overline{\mathbf{W}}_p) =\overline{\mathbf{W}}_p + \dfrac{1}{d+1} \sum_{\ell \in \mathcal{C}_p}(  \mathbf{W}_c - \mathbf{W}_{\ell}), 
	$$
	where $\overline{\mathbf{W}}_p = \sum_{\ell \in \mathcal{C}_p} \mathbf{W}_{\ell}/(d+1)$. As a consequence we have 
	$$
	\sum_c|\Omega_c|  \dfrac{d\eta(\Q_c)}{dt} + \sum_c\sum\limits_{p\in \mathcal{P}_c} \overline{\mathbf{W}}_p^t \hat\F_p \cdot\n_{pc}
	+  \sum_c\sum\limits_{p\in \mathcal{P}_c}  (   \mathbf{W}_c - \overline{\mathbf{W}}_p)^t\hat\F_p \cdot\n_{pc} =0.
	$$
	Assume now that  $\forall \,c\in\mathcal{C}_p$ we have
	\begin{equation} \label{eq.shuffle}
		\sum_{c\in\mathcal{C}_p}   \mathbf{W}_c \hat\F_p \cdot\n_{pc}  =   \sum_{c\in\mathcal{C}_p}    \boldsymbol{\Psi}_c  \cdot\n_{pc},
	\end{equation}
	which is equivalent to
	\begin{equation} \label{eq.shuffle1}
		\sum_{c\in\mathcal{C}_p}  ( \mathbf{W}_c - \overline{\mathbf{W}}_p)^t  \hat\F_p \cdot\n_{pc}  =   \sum_{c\in\mathcal{C}_p}   ( \boldsymbol{\Psi}_c -\overline{\boldsymbol{\Psi}}_p )  \cdot\n_{pc}.
	\end{equation}
	Under this condition  we can write
	$$
	\sum_c|\Omega_c|  \dfrac{d\eta(\Q_c)}{dt} + \sum_c\sum\limits_{p\in \mathcal{P}_c}( \overline{\mathbf{W}}_p^t \hat\F_p -  \overline{\boldsymbol{\Psi}}_p) \cdot\n_{pc}
	=  -\sum_c   \boldsymbol{\Psi}_c  \cdot\sum\limits_{p\in \mathcal{P}_c} \n_{pc} =0,
	$$
	which shows that under the generalized shuffle condition~\eqref{eq.shuffle}  entropy conservation is verified with the multidimensional numerical entropy flux    
	$$
	\hat {\mathbf{G}}_p     := \overline{\mathbf{W}}_p^t \hat\F_p   - \overline{\boldsymbol{\Psi}}_p .
	$$
	As argued in~\cite{tadmor} any scheme with more numerical viscosity than an entropy conservative scheme will be entropy stable.\\
	
	As an alternative approach we now consider the schemes in fluctuation form. In absence of reconstruction, and proceeding as before this gives
	$$
	\sum_c|\Omega_c|  \dfrac{d\eta(\Q_c)}{dt} +   \sum_c\sum\limits_{p\in \mathcal{P}_c}  \mathbf{W}_c^t\boldsymbol{\phi}_{pc}   =0.
	$$
	We can add the sum of the internal entropy value flux  along the nodal normals, which is null due to the properties of the normals:  
	$$
	\sum_c|\Omega_c|  \dfrac{d\eta(\Q_c)}{dt} +   \sum_c\sum\limits_{p\in \mathcal{P}_c}  \mathbf{G}_c\cdot \n_{pc} +  \sum_c\sum\limits_{p\in \mathcal{P}_c}  \mathbf{W}_c^t\boldsymbol{\phi}_{pc}   =0,
	$$
	with $ \mathbf{G}_c = \mathbf{G}(\Q_c)$. We can easily see that a scheme is entropy conservative if for any corner $p$ we have 
	\begin{equation} \label{eq.entropy_fluct}
		\sum\limits_{c\in \mathcal{C}_p}  \mathbf{W}_{c}^t\boldsymbol{\phi}_{pc}  = \boldsymbol{\Phi}_{\eta}:=\sum\limits_{c \in \mathcal{C}_p}  \mathbf{G}_{c}\cdot\n_{c p},
	\end{equation}
	with associated numerical  entropy flux verifying 
	$$
	\hat{\mathbf{G}}_{p}\cdot \n_{pc} = {\mathbf{G}}_{c}\cdot \n_{pc}  + \mathbf{W}_c^t\boldsymbol{\phi}_{pc}.  
	$$
	\begin{proposition}[Equivalence of the entropy condition in potential/fluctuation  form] The potential form of the shuffle condition~\eqref{eq.shuffle} and the fluctuation  condition~\eqref{eq.entropy_fluct}
		are equivalent.
	\end{proposition}
	\begin{proof}
		We prove that one condition implies the other and  vice-versa.
		This first implication be easily shown by  a direct computation, summing up the shuffle conditions around  node $p$:
		$$
		\begin{aligned}
			\sum_{c\in\mathcal{C}_p}  (\mathbf{W}_c - \overline{\mathbf{W}}_p)^t \F_c\cdot \n_{pc} +  \sum_{c\in\mathcal{C}_p} (\mathbf{W}_c - \overline{\mathbf{W}}_p)^t\boldsymbol{\phi}_{pc} =  \sum_{c\in\mathcal{C}_p}( \boldsymbol{\Psi}_c - \bar{\boldsymbol{\Psi}}_p) \cdot n_{pc}\\[10pt]
			\sum_{c\in\mathcal{C}_p}  \mathbf{W}_c^t   \F_c\cdot \n_{pc} +   \overline{\mathbf{W}}_p^t  \boldsymbol{\phi}_p  +  \sum_{c\in\mathcal{C}_p} \mathbf{W}_c^t \boldsymbol{\phi}_{pc} -  \overline{\mathbf{W}}_p^t\boldsymbol{\phi}_{p} =  \sum_{c\in\mathcal{C}_p} \boldsymbol{\Psi}_c  \cdot n_{pc},\\
		\end{aligned}
		$$
		which readily shows that 
		$$
		\begin{aligned}
			\sum_{c\in\mathcal{C}_p} \mathbf{W}_c^t \boldsymbol{\phi}_{pc} -
			=  \sum_{c\in\mathcal{C}_p}  \mathbf{W}_c^t   \F_c\cdot \n_{cp} - \sum_{c\in\mathcal{C}_p} \boldsymbol{\Psi}_c  \cdot n_{cp}=
			\sum_{c\in\mathcal{C}_p}  \mathbf{G}_c \cdot \n_{cp}.
		\end{aligned}
		$$
		The reverse implication is also  easy to show using the definition of the $\boldsymbol{\phi}_{pc}$:
		$$
		\begin{aligned}
			\sum_{c\in\mathcal{C}_p} \mathbf{W}_c^t   ( \hat\F_p \cdot  \n_{pc}  - \F_c \cdot  \n_{pc} )    =       \sum_{c\in\mathcal{C}_p}  \mathbf{G}_c\cdot \n_{cp} \\ 
			\sum_{c\in\mathcal{C}_p} \mathbf{W}_c^t     ( \hat \F_p \cdot  \n_{pc}  - \F_c \cdot  \n_{pc} )    =    - \sum_{c\in\mathcal{C}_p}  (\mathbf{W}_c^t \F_c  \cdot  \n_{pc}  - \boldsymbol{\Psi}_c  \cdot  \n_{pc}) \\
			\sum_{c\in\mathcal{C}_p} \mathbf{W}_c^t      \hat \F_p \cdot  \n_{pc}      =     \sum_{c\in\mathcal{C}_p}    \boldsymbol{\Psi}_c  \cdot  \n_{pc} 
		\end{aligned}
		$$
	\end{proof}
	
	\bibliographystyle{plain}
	\bibliography{./referencesMultiD-RS-black.bib}

\begin{thebibliography}{100}

\bibitem{abgrall_eno}
R.~Abgrall.
\newblock On essentially non-oscillatory schemes on unstructured meshes:
  analysis and implementation.
\newblock {\em J. Comput. Phys.}, 144:45--58, 1994.

\bibitem{abg2017}
R.~Abgrall.
\newblock High order schemes for hyperbolic problems using globally continuous
  approximation and avoiding mass matrices.
\newblock {\em J. Sci. Comput.}, 73:461--494, 2017.

\bibitem{abgrall2018general}
R.~Abgrall.
\newblock {A general framework to construct schemes satisfying additional
  conservation relations. Application to entropy conservative and entropy
  dissipative schemes}.
\newblock {\em J. Comput. Phys.}, 372:640--666, 2018.

\bibitem{Abgrall23}
R.~Abgrall.
\newblock A personal discussion on conservation, and how to formulate it.
\newblock In Emmanuel Franck, J{\"u}rgen Fuhrmann, Victor Michel-Dansac, and
  Laurent Navoret, editors, {\em Finite Volumes for Complex Applications
  X---Volume 1, Elliptic and Parabolic Problems}, pages 3--19, Cham, 2023.
  Springer Nature Switzerland.

\bibitem{abgrall2023extensions}
R.~Abgrall and W.~Barsukow.
\newblock Extensions of active flux to arbitrary order of accuracy.
\newblock {\em ESAIM: Mathematical Modelling and Numerical Analysis},
  57(2):991--1027, 2023.

\bibitem{abg_bar_klin_2025semi}
R.~Abgrall, W.~Barsukow, and C.~Klingenberg.
\newblock {A Semi-discrete Active Flux Method for the Euler Equations on
  Cartesian Grids}.
\newblock {\em Journal of Scientific Computing}, 102(2):36, 2025.

\bibitem{barth02}
R.~Abgrall and T.J. Barth.
\newblock Residual distribution schemes for conservation laws via adaptive
  quadrature.
\newblock {\em {SIAM} J. Sci. Comput.}, 24(3):732--769, 2002.

\bibitem{alr11}
R.~Abgrall, A.~Larat, and M.~Ricchiuto.
\newblock Construction of very high order residual distribution schemes for
  steady inviscid flow problems on hybrid unstructured meshes.
\newblock {\em Journal of Computational Physics}, 230(11):4103 -- 4136, 2011.

\bibitem{abg2001c}
R.~Abgrall and M.~Mezine.
\newblock Construction of second-order accurate monotone and stable residual
  distribution schemes for unsteady flow problems.
\newblock {\em J. Comput. Phys.}, 188:16--55, 2003.

\bibitem{abgrall2022reinterpretation}
R.~Abgrall, P.~{\"O}ffner, and H.~Ranocha.
\newblock {Reinterpretation and extension of entropy correction terms for
  residual distribution and discontinuous Galerkin schemes: Application to
  structure preserving discretization}.
\newblock {\em Journal of Computational Physics}, 453:110955, 2022.

\bibitem{RD-ency2}
R.~Abgrall and M.~Ricchiuto.
\newblock High-order methods for cfd.
\newblock In {\em Encyclopedia of Computational Mechanics Second Edition},
  pages 1--54. John Wiley \& Sons, Ltd, 2017.

\bibitem{abgrall2022hyperbolic}
R.~Abgrall and M.~Ricchiuto.
\newblock Hyperbolic balance laws: residual distribution, local and global
  fluxes.
\newblock {\em Numerical Fluid Dynamics}, pages 177--222, 2022.

\bibitem{abg2003}
R.~Abgrall and P.L. Roe.
\newblock High-order fluctuation schemes on triangular meshes.
\newblock {\em J. Sci. Comput.}, 19(3):3--36, 2003.

\bibitem{balsarahlle2d}
D.S. Balsara.
\newblock {Multidimensional HLLE Riemann solver: Application to Euler and
  magnetohydrodynamic flows}.
\newblock {\em J. Comput. Phys.}, 229:1970--1993, 2010.

\bibitem{balsarahllc2d}
D.S. Balsara.
\newblock {A two-dimensional HLLC Riemann solver for conservation laws:
  Application to Euler and magnetohydrodynamic flows}.
\newblock {\em J. Comput. Phys.}, 231:7476--7503, 2012.

\bibitem{MUSIC1}
D.S. Balsara.
\newblock {Multidimensional Riemann problem with self-similar internal
  structure – Part I – Application to hyperbolic conservation laws on
  structured meshes}.
\newblock {\em J. Comput. Phys.}, 277:163--200, 2014.

\bibitem{ADERdivB}
D.S. Balsara and M.~Dumbser.
\newblock {Divergence-free MHD on unstructured meshes using high order finite
  volume schemes based on multidimensional {R}iemann solvers}.
\newblock {\em Journal of Computational Physics}, 299:687--715, 2015.

\bibitem{MUSIC2}
D.S. Balsara and M.~Dumbser.
\newblock {Multidimensional Riemann problem with self-similar internal
  structure – Part II – Application to hyperbolic conservation laws on
  unstructured meshes}.
\newblock {\em J. Comput. Phys.}, 287:269--292, 2015.

\bibitem{BalsaraMultiDRS}
D.S. Balsara, M.~Dumbser, and R.~Abgrall.
\newblock {Multidimensional HLLC Riemann Solver for Unstructured Meshes - With
  Application to Euler and MHD Flows.}
\newblock {\em J. Comput. Phys.}, 261:172--208, 2014.

\bibitem{barsukow2019active}
W.~Barsukow, J.~Hohm, C.~Klingenberg, and P.~Roe.
\newblock {The active flux scheme on Cartesian grids and its low Mach number
  limit}.
\newblock {\em Journal of Scientific Computing}, 81:594--622, 2019.

\bibitem{barthlsq}
T.J. Barth and P.O. Frederickson.
\newblock Higher order solution of the {Euler} equations on unstructured grids
  using quadratic reconstruction.
\newblock {\em AIAA paper no. 90-0013}, 28th Aerospace Sciences Meeting January
  1990.

\bibitem{BarthJespersen}
T.J. Barth and D.C. Jespersen.
\newblock {The design and application of upwind schemes on unstructured
  meshes}.
\newblock {\em AIAA Paper 89-0366}, pages 1--12, 1989.

\bibitem{LagrangeMHD}
W.~Boscheri, M.~Dumbser, and D.S. Balsara.
\newblock {High Order Lagrangian ADER-WENO Schemes on Unstructured Meshes --
  Application of Several Node Solvers to Hydrodynamics and
  Magnetohydrodynamics}.
\newblock {\em International Journal for Numerical Methods in Fluids},
  76:737--778, 2014.

\bibitem{frontiers2020}
S.~Busto, S.~Chiocchetti, M.~Dumbser, E.~Gaburro, and I.~Peshkov.
\newblock {High order ADER schemes for continuum mechanics}.
\newblock {\em Frontiers in Physics}, 8:32, 2020.

\bibitem{Despres2009}
G.~Carr\'e, S.~Del Pino, B.~Despr\'es, and E.~Labourasse.
\newblock {A cell-centered Lagrangian hydrodynamics scheme on general
  unstructured meshes in arbitrary dimension.}
\newblock {\em Journal of Computational Physics}, 228:5160--5183, 2009.

\bibitem{ciallella2023shifted}
M.~Ciallella, E.~Gaburro, M.~Lorini, and M.~Ricchiuto.
\newblock Shifted boundary polynomial corrections for compressible flows: high
  order on curved domains using linear meshes.
\newblock {\em Applied Mathematics and Computation}, 441:127698, 2023.

\bibitem{ColellaMultiD}
P.~Colella.
\newblock {Multidimensional upwind methods for hyperbolic conservation laws}.
\newblock {\em Journal of Computational Physics}, 87:171--200, 1990.

\bibitem{cravero2018cweno}
I.~Cravero, G.~Puppo, M.~Semplice, and G.~Visconti.
\newblock Cweno: uniformly accurate reconstructions for balance laws.
\newblock {\em Math. Comput.}, 87(312):1689--1719, 2018.

\bibitem{crd02}
\'A. Cs\'{\i}k, M.~Ricchiuto, and H.~Deconinck.
\newblock {A Conservative Formulation of the Multidimensional Upwind Residual
  Distribution Schemes for General Nonlinear Conservation Laws}.
\newblock {\em J. Comput. Phys}, 179(2):286--312, 2002.

\bibitem{Dec93}
H.~Deconinck.
\newblock Beyond the riemann problem, part ii.
\newblock In M.~Y. Hussaini, A.~Kumar, and M.~D. Salas, editors, {\em
  Algorithmic Trends in Computational Fluid Dynamics}, pages 341--367, New
  York, NY, 1993. Springer New York.

\bibitem{RD-ency}
H.~Deconinck and M.~Ricchiuto.
\newblock Residual distribution schemes: Foundations and analysis.
\newblock In {\em Encyclopedia of Computational Mechanics Second Edition},
  pages 1--53. John Wiley \& Sons, Ltd, 2017.

\bibitem{Deconinck1993}
H.~Deconinck, P.L. Roe, and R.~Struijs.
\newblock {A multidimensional generalization of Roe's flux difference splitter
  for the Euler equations}.
\newblock {\em Computers and Fluids}, 22:215--222, 1993.

\bibitem{deconinck93}
H.~Deconinck, P.L. Roe, and R.~Struijs.
\newblock A multidimensional generalization of {Roe's} difference splitter for
  the{ Euler} equations.
\newblock {\em Computers and Fluids}, 22(2/3):215--222, 1993.

\bibitem{PHRaph2}
A.~{Del Grosso}, M.J. Castro, A.~Chan, G.~Gallice, R.~Loub\`ere, and P.H.
  Maire.
\newblock A well-balanced, positive, entropy--stable, and
  multi--dimensional--aware finite volume scheme for 2d shallow--water
  equations with unstructured grids.
\newblock {\em Journal of Computational Physics}, 503:112829, 2024.

\bibitem{Despres2005}
B.~Despr\'es and C.~Mazeran.
\newblock {Lagrangian gas dynamics in two-dimensions and Lagrangian systems}.
\newblock {\em Archive for Rational Mechanics and Analysis}, 178:327--372,
  2005.

\bibitem{hllem}
M.~Dumbser and D.~S. Balsara.
\newblock A new efficient formulation of the {HLLEM} {R}iemann solver for
  general conservative and non-conservative hyperbolic systems.
\newblock {\em Journal of Computational Physics}, 304(C):275--319, 2016.

\bibitem{SIMHD}
M.~Dumbser, D.S. Balsara, M.~Tavelli, and F.~Fambri.
\newblock {A divergence-free semi-implicit finite volume scheme for ideal,
  viscous and resistive magnetohydrodynamics}.
\newblock {\em International Journal for Numerical Methods in Fluids},
  89:16--42, 2019.

\bibitem{dumbser2008unified}
M.~Dumbser, D.S. Balsara, E.F. Toro, and C.-D. Munz.
\newblock {A unified framework for the construction of one-step finite volume
  and discontinuous Galerkin schemes on unstructured meshes}.
\newblock {\em Journal of Computational Physics}, 227(18):8209--8253, 2008.

\bibitem{ADERCWENO}
M.~Dumbser, W.~Boscheri, M.~Semplice, and G.~Russo.
\newblock Central weighted eno schemes for hyperbolic conservation laws on
  fixed and moving unstructured meshes.
\newblock {\em SIAM J. Sci. Comput.}, 39(6):A2564--A2591, 2017.

\bibitem{DumbserEnauxToro}
M.~Dumbser, C.~Enaux, and E.F. Toro.
\newblock {Finite Volume Schemes of Very High Order of Accuracy for Stiff
  Hyperbolic Balance Laws}.
\newblock {\em J. Comput. Phys.}, 227:3971--4001, 2008.

\bibitem{DumbserKaeser06b}
M.~Dumbser and M.~K\"aser.
\newblock {Arbitrary high order non-oscillatory Finite Volume schemes on
  unstructured meshes for linear hyperbolic systems}.
\newblock {\em J. Comput. Phys.}, 221:693--723, 2007.

\bibitem{DumbserKaeser07}
M.~Dumbser, M.~K\"aser, V.A Titarev, and E.F. Toro.
\newblock {Quadrature-Free Non-Oscillatory Finite Volume Schemes on
  Unstructured Meshes for Nonlinear Hyperbolic Systems}.
\newblock {\em J. Comput. Phys.}, 226:204--243, 2007.

\bibitem{OsherUniversal}
M.~Dumbser and E.~F. Toro.
\newblock {On Universal {Osher}--Type Schemes for General Nonlinear Hyperbolic
  Conservation Laws}.
\newblock {\em Communications in Computational Physics}, 10:635--671, 2011.

\bibitem{OsherNC}
M.~Dumbser and E.~F. Toro.
\newblock A simple extension of the {Osher} {Riemann} solver to
  non-conservative hyperbolic systems.
\newblock {\em J. Sci. Comput.}, 48:70--88, 2011.

\bibitem{HybridHexa1}
F.~Fambri, E.~Zampa, S.~Busto, L.~R\'{\i}o-Mart\'{\i}n, F.~Hindenlang,
  E.~Sonnendr\"{u}cker, and M.~Dumbser.
\newblock A well-balanced and exactly divergence-free staggered semi-implicit
  hybrid finite volume / finite element scheme for the incompressible {MHD}
  equations.
\newblock {\em Journal of Computational Physics}, 493:112493, 2023.

\bibitem{gaburro2021unified}
E.~Gaburro.
\newblock {A unified framework for the solution of hyperbolic PDE systems using
  high order direct Arbitrary-Lagrangian--Eulerian schemes on moving
  unstructured meshes with topology change}.
\newblock {\em Archives of Computational Methods in Engineering},
  28(3):1249--1321, 2021.

\bibitem{gaburro2025high}
E.~Gaburro.
\newblock High order well-balanced arbitrary-lagrangian-eulerian ader
  discontinuous galerkin schemes on general polygonal moving meshes.
\newblock {\em Computers \& Fluids}, page 106764, 2025.

\bibitem{gaburro2020high}
E.~Gaburro, W.~Boscheri, S.~Chiocchetti, C.~Klingenberg, V.~Springel, and
  M.~Dumbser.
\newblock High order direct arbitrary-lagrangian-eulerian schemes on moving
  voronoi meshes with topology changes.
\newblock {\em Journal of Computational Physics}, 407:109167, 2020.

\bibitem{GaburroRicchiutoPrimitive}
E.~Gaburro, W.~Boscheri, S.~Chiocchetti, and M.~Ricchiuto.
\newblock {Discontinuous Galerkin schemes for hyperbolic systems in
  non-conservative variables: quasi-conservative formulation with subcell
  finite volume corrections}.
\newblock {\em Computer Methods in Applied Mechanics and Engineering},
  431:117311, 2024.

\bibitem{PHRaph1}
G.~Gallice, A.~Chan, R.~Loub\`ere, and P.H. Maire.
\newblock Entropy stable and positivity preserving godunov-type schemes for
  multidimensional hyperbolic systems on unstructured grid.
\newblock {\em Journal of Computational Physics}, 468:111493, 2022.

\bibitem{GARICANOMENA201643}
J.~Garicano-Mena, A.~Lani, and H.~Deconinck.
\newblock {An energy-dissipative remedy against carbuncle: Application to
  hypersonic flows around blunt bodies}.
\newblock {\em Computers and Fluids}, 133:43--54, 2016.

\bibitem{godunov}
{S.K.} Godunov.
\newblock {Finite Difference Methods for the Computation of Discontinuous
  Solutions of the Equations of Fluid Dynamics}.
\newblock {\em Mathematics of the USSR: Sbornik}, 47:271--306, 1959.

\bibitem{harten_eno}
{A.} Harten, {B.} Engquist, {S.} Osher, and {S.R.} Chakravarthy.
\newblock {Uniformly High Order Accurate Essentially Non--oscillatory Schemes
  {III}}.
\newblock {\em J. Comput. Phys.}, 71:231--303, 1987.

\bibitem{hll}
A.~Harten, {P. D.} Lax, and B.~{van Leer}.
\newblock On upstream differencing and {Godunov}-type schemes for hyperbolic
  conservation laws.
\newblock {\em SIAM Review}, 25(1):35--61, 1983.

\bibitem{hu1999weighted}
C.~Hu and C.-W. Shu.
\newblock {Weighted essentially non-oscillatory schemes on triangular meshes}.
\newblock {\em Journal of Computational Physics}, 150(1):97--127, 1999.

\bibitem{HuShuVortex1999}
C.~Hu and {C.W.} Shu.
\newblock {A high-order WENO finite difference scheme for the equations of
  ideal magnetohydrodynamics.}
\newblock {\em Journal of Computational Physics}, 150:561 -- 594, 1999.

\bibitem{jackson2017eigenvalues}
H.~Jackson.
\newblock {On the eigenvalues of the ADER-WENO Galerkin predictor}.
\newblock {\em J. Comput. Phys.}, 333:409--413, 2017.

\bibitem{shu_efficient_weno}
{G.S.} Jiang and {C.W.} Shu.
\newblock Efficient implementation of weighted {ENO} schemes.
\newblock {\em J. Comput. Phys.}, 126:202--228, 1996.

\bibitem{kemm2020simple}
F.~Kemm, E.~Gaburro, F.~Thein, and M.~Dumbser.
\newblock {A simple diffuse interface approach for compressible flows around
  moving solids of arbitrary shape based on a reduced Baer--Nunziato model}.
\newblock {\em Computers \& fluids}, 204:104536, 2020.

\bibitem{kepler}
J.~Kepler.
\newblock {\em {Nova stereometria doliorium vinariorum, in primis Austriaci}}.
\newblock Plancus, Lincii, Austria, 1615.

\bibitem{CiCP-13-479}
A.~Lani, M.~Panesi, and H.~Deconinck.
\newblock {Conservative Residual Distribution Method for Viscous Double Cone
  Flows in Thermochemical Nonequilibrium}.
\newblock {\em Communications in Computational Physics}, 13(2):479--501, 2013.

\bibitem{laxliu98}
P.~Lax and X.D. Liu.
\newblock {Solution of two-dimensional Riemann problems of gas dynamics by
  positive schemes}.
\newblock {\em SIAM Journal of Scientific Computing}, 19:319--340, 1998.

\bibitem{lax2}
P.D. Lax.
\newblock {Weak solutions of nonlinear hyperbolic equations and their numerical
  approximation}.
\newblock {\em Comm. Pure Appl. Math.}, 7:159--193, 1954.

\bibitem{CWENO1}
D.~Levy, G.~Puppo, and G.~Russo.
\newblock Central {WENO} schemes for hyperbolic systems of conservation laws.
\newblock {\em M2AN Math. Model. Numer. Anal.}, 33(3):547--571, 1999.

\bibitem{CWENO2}
D.~Levy, G.~Puppo, and G.~Russo.
\newblock {A third order central {WENO} scheme for 2D conservation laws}.
\newblock {\em Applied Numerical Mathematics}, 33:415--421, 2000.

\bibitem{CWENO3}
D.~Levy, G.~Puppo, and G.~Russo.
\newblock {A fourth-order central {WENO} scheme for multidimensional hyperbolic
  systems of conservation laws}.
\newblock {\em SIAM J. Sci. Comput.}, 24:480--506, 2002.

\bibitem{FVEG}
M.~Lukacovaá-Medvidova, G.~Warnecke, and Y.~Zahaykah.
\newblock {Finite volume evolution Galerkin (FVEG) methods for
  three-dimensional wave equation system}.
\newblock {\em Applied Numerical Mathematics}, 57:1050--1064, 2007.

\bibitem{Maire2009}
{P.-H.} Maire.
\newblock A high-order cell-centered {Lagrangian} scheme for two-dimensional
  compressible fluid flows on unstructured meshes.
\newblock {\em J. Comput. Phys.}, 228:2391--2425, 2009.

\bibitem{Maire2007}
P.-H. Maire, R.~Abgrall, J.~Breil, and J.~Ovadia.
\newblock A cell-centered lagrangian scheme for two-dimensional compressible
  flow problems.
\newblock {\em SIAM J. Sci. Comput.}, 29:1781--1824, 2007.

\bibitem{Maire2008}
P.-H. Maire and J.~Breil.
\newblock A second-order cell-centered lagrangian scheme for two-dimensional
  compressible flow problems.
\newblock {\em Int. J. Numer. Methods Fluids}, 56:1417--1423, 2007.

\bibitem{Maire2011}
P.H. Maire.
\newblock A high-order one-step sub-cell force-based discretization for
  cell-centered {Lagrangian} hydrodynamics on polygonal grids.
\newblock {\em Computers and Fluids}, 46(1):341--347, 2011.

\bibitem{Maire2010}
P.H. Maire.
\newblock A unified sub-cell force-based discretization for cell-centered
  {Lagrangian} hydrodynamics on polygonal grids.
\newblock {\em International Journal for Numerical Methods in Fluids},
  65:1281--1294, 2011.

\bibitem{hiro3}
H.~Nishikawa, M.~Rad, and P.L. Roe.
\newblock {A third-order fluctuation splitting scheme that preserves potential
  flow}.
\newblock 15th {AIAA} Computational Fluid Dynamics Conference, Anaheim, {CA},
  {USA}, June 2001.

\bibitem{osherandsolomon}
S.~Osher and F.~Solomon.
\newblock {Upwind Difference Schemes for Hyperbolic Conservation Laws}.
\newblock {\em Math. Comput.}, 38:339--374, 1982.

\bibitem{mr11}
M.~Ricchiuto.
\newblock {\em {Contributions to the development of residual discretizations
  for hyperbolic conservation laws with application to shallow water flows}}.
\newblock Habilitation {\`a} diriger des recherches, {Universit{\'e} Sciences
  et Technologies - Bordeaux I}, December 2011.
\newblock PDF available at {\tt https://theses.hal.science/tel-00651688/}.

\bibitem{HDR-MR}
M.~Ricchiuto.
\newblock {\em {Contributions to the development of residual discretizations
  for hyperbolic conservation laws with application to shallow water flows}}.
\newblock Habilitation {\`a} diriger des recherches, {Universit{\'e} Sciences
  et Technologies - Bordeaux I}, December 2011.
\newblock PDF available at {\tt https://theses.hal.science/tel-00651688/}.

\bibitem{ra10}
M.~Ricchiuto and R.~Abgrall.
\newblock Explicit runge-kutta residual distribution schemes for time dependent
  problems: Second order case.
\newblock {\em Journal of Computational Physics}, 229(16):5653 -- 5691, 2010.

\bibitem{rad07}
M.~Ricchiuto, R.~Abgrall, and H.~Deconinck.
\newblock Application of conservative residual distribution schemes to the
  solution of the shallow water equations on unstructured meshes,.
\newblock {\em J. Comput. Phys.}, 222:287--331, 2007.

\bibitem{rcd05}
M.~Ricchiuto, \'A. Cs\'{\i}k, and H.~Deconinck.
\newblock {Residual distribution for general time dependent conservation laws}.
\newblock {\em J. Comput. Phys}, 209(1):249--289, 2005.

\bibitem{RODIONOV2017308}
A.V. Rodionov.
\newblock Artificial viscosity in godunov-type schemes to cure the carbuncle
  phenomenon.
\newblock {\em Journal of Computational Physics}, 345:308--329, 2017.

\bibitem{Roe:87}
P.~L. Roe.
\newblock Linear advection schemes on triangular meshes.
\newblock Technical Report CoA 8720, Cranfield Institute of Technology, 1987.

\bibitem{Roe:90}
P.~L. Roe.
\newblock {``Optimum'' upwind advection on a triangular mesh}.
\newblock Technical Report ICASE 90-75, ICASE, 1990.

\bibitem{Roe92}
P.~L. Roe and D.~Sidilkover.
\newblock Optimum positive linear schemes for advection in two and three
  dimensions.
\newblock {\em SIAM Journal on Numerical Analysis}, 29(6):1542--1568, 1992.

\bibitem{roe1d}
{P.L.} Roe.
\newblock {Approximate {Riemann} Solvers, Parameter vectors, and Difference
  schemes}.
\newblock {\em J. Comput. Phys.}, 43:357--372, 1981.

\bibitem{Roe82}
P.L. Roe.
\newblock Fluctuations and signals - a framework for numerical evolution
  problems.
\newblock In K.W. Morton and M.J. Baines, editors, {\em Numerical Methods for
  Fluids Dynamics}, pages 219--257. Academic Press, 1982.

\bibitem{RoeMultiD}
P.L. Roe.
\newblock {Discrete models for the numerical analysis of time--dependent
  multidimensional gas dynamics}.
\newblock {\em Journal of Computational Physics}, 63:458--476, 1986.

\bibitem{Roe93}
P.L. Roe.
\newblock Beyond the riemann problem, part i.
\newblock In M.~Y. Hussaini, A.~Kumar, and M.~D. Salas, editors, {\em
  Algorithmic Trends in Computational Fluid Dynamics}, pages 341--367, New
  York, NY, 1993. Springer New York.

\bibitem{Roe24}
P.L. Roe.
\newblock Musings of a computational philosopher.
\newblock In James~C. Tyacke and Nagabhushana~Rao Vadlamani, editors, {\em
  Proceedings of the Cambridge Unsteady Flow Symposium 2024}, pages 1--35,
  Cham, 2025. Springer Nature Switzerland.

\bibitem{RumseyLeerRoe}
C.L. Rumsey, B.~van Leer, and P.L. Roe.
\newblock {A multidimensional flux function with applications to the Euler and
  Navier-Stokes equations}.
\newblock {\em Journal of Computational Physics}, 105:306--323, 1993.

\bibitem{SaltzmanMultiD}
J.~Saltzman.
\newblock {An unsplit 3D upwind method for hyperbolic conservation laws}.
\newblock {\em Journal of Computational Physics}, 115:153--168, 1994.

\bibitem{schneider2021multidimensional}
K.~A Schneider, J.M. Gallardo, D.S. Balsara, B.~Nkonga, and C.~Par{\'e}s.
\newblock {Multidimensional approximate Riemann solvers for hyperbolic
  nonconservative systems. Applications to shallow water systems}.
\newblock {\em Journal of Computational Physics}, 444:110547, 2021.

\bibitem{schulzrinne}
C.~W. Schulz-Rinne.
\newblock Classification of the {Riemann} problem for two-dimensional gas
  dynamics.
\newblock {\em SIAM J. Math. Anal.}, 24:76--88, 1993.

\bibitem{RinneCollinsGlaz}
C.W. Schulz-Rinne, J.P. Collins, and H.M. Glaz.
\newblock {Numerical solution of the Riemann problem for two-dimensional
  gasdynamics}.
\newblock {\em SIAM Journal on Scientific Computing}, 14:1394--1414, 1993.

\bibitem{rd-efix}
K.~Sermeus and H.~Deconinck.
\newblock An entropy fix for multidimensional upwind residual distribution
  schemes.
\newblock {\em Computers and Fluids}, 34(4):617--640, 2005.

\bibitem{stroud}
{A.H.} Stroud.
\newblock {\em {Approximate Calculation of Multiple Integrals}}.
\newblock Prentice-Hall Inc., Englewood Cliffs, New Jersey, 1971.

\bibitem{tadmor}
E.~Tadmor.
\newblock Entropy stability theory for difference approximations of nonlinear
  conservation laws and related time-dependent problems.
\newblock {\em Acta Numerica}, 12:451--512, 2003.

\bibitem{toro3}
{V.A.} Titarev and {E.F.} Toro.
\newblock {ADER}: Arbitrary high order {Godunov} approach.
\newblock {\em Journal of Scientific Computing}, 17(1-4):609--618, December
  2002.

\bibitem{MixedWENO2D}
V.A. Titarev, P.~Tsoutsanis, and D.~Drikakis.
\newblock {WENO schemes for mixed--element unstructured meshes}.
\newblock {\em Communications in Computational Physics}, 8:585--609, 2010.

\bibitem{torobook}
E.~F. Toro.
\newblock {\em {Riemann Solvers and Numerical Methods for Fluid Dynamics. A
  Practical Introduction, Third edition}}.
\newblock Springer-Verlag, Berlin, 2009.

\bibitem{Torohllc}
E.F. Toro, M.~Spruce, and W.~Speares.
\newblock Restoration of the contact surface in the {Harten-Lax-van Leer
  Riemann} solver.
\newblock {\em Journal of Shock Waves}, 4:25--34, 1994.

\bibitem{toro2005ader}
E.F. Toro and V.A. Titarev.
\newblock {ADER schemes for scalar non-linear hyperbolic conservation laws with
  source terms in three-space dimensions}.
\newblock {\em Journal of Computational Physics}, 202(1):196--215, 2005.

\bibitem{MixedWENO3D}
P.~Tsoutsanis, V.A. Titarev, and D.~Drikakis.
\newblock {WENO schemes on arbitrary mixed-element unstructured meshes in three
  space dimensions}.
\newblock {\em J. Comput. Phys.}, 230:1585--1601, 2011.

\bibitem{vanderweide}
E.~van~der Weide and H.~Deconinck.
\newblock {Positive matrix distribution schemes for hyperbolic systems}.
\newblock In {\em Computational Fluid Dynamics}, pages 747--753, New York,
  1996. Wiley.

\bibitem{CompatibleMHD}
E.~Zampa, S.~Busto, and M.~Dumbser.
\newblock A divergence-free hybrid finite volume / finite element scheme for
  the incompressible mhd equations based on compatible finite element spaces.
\newblock {\em Applied Numerical Mathematics}, 198:346--374, 2024.

\bibitem{ZhangShu3D}
{Y.T.} Zhang and {C.W.} Shu.
\newblock {Third Order {WENO} Scheme on Three Dimensional Tetrahedral Meshes }.
\newblock {\em Communications in Computational Physics}, 5:836--848, 2009.

\end{thebibliography}

\end{document}